\documentclass[12pt]{article}

\input{epsf}
\usepackage{amsfonts,amsmath}

\usepackage{enumitem}
\usepackage{amsthm}

\usepackage[pdftex]{graphicx}

\usepackage{xcolor,soul}

\def\cA{\mathcal{A}}
\def\cB{\mathcal{B}}

\def\cD{\mathcal{D}}
\def\cE{\mathcal{E}}
\def\cF{\mathcal{F}}
\def\cG{\mathcal{G}}
\def\cI{\mathcal{I}}
\def\cJ{\mathcal{J}}
\def\cM{\mathcal{M}}
\def\cP{\mathcal{P}}

\newtheorem{theorem}{Theorem}[subsection]
\newtheorem{lemma}[theorem]{Lemma}
\newtheorem{corollary}[theorem]{Corollary}
\newtheorem{definition}[theorem]{Definition}
\newtheorem{proposition}[theorem]{Proposition}
\newtheorem{example}[theorem]{Example}
\newtheorem{remark}[theorem]{Remark}

\def\bit{\begin{itemize}}
\def\eit{\end{itemize}}

\def\bc{\begin{center}}
\def\ec{\end{center}}

\def\bthm{\begin{theorem}}
\def\ethm{\end{theorem}}

\def\bcor{\begin{corollary}}
\def\ecor{\end{corollary}}

\def\bprop{\begin{proposition}}
\def\eprop{\end{proposition}}

\def\blem{\begin{lemma}}
\def\elem{\end{lemma}}

\def\bex{\begin{example}}
\def\eex{\end{example}}

\def\brem{\begin{remark}}
\def\erem{\end{remark}}

\def\prf{\noindent{\bf Proof. }}

\def\bdes{\begin{description}}
\def\edes{\end{description}}

\def\ita{\item[(a)]}
\def\itb{\item[(b)]}
\def\itc{\item[(c)]}
\def\itd{\item[(d)]}
\def\ite{\item[(e)]}

\def\iti{\item[(i)]}
\def\itii{\item[(ii)]}
\def\itiii{\item[(iii)]}

\def\beq{\begin{equation}}
\def\eeq{\end{equation}}

\def\ben{\begin{enumerate}}
\def\een{\end{enumerate}}

\def\beqar{\begin{eqnarray}}
\def\eeqar{\end{eqnarray}}

\def\beqarr{\begin{eqnarray*}}
\def\eeqarr{\end{eqnarray*}}

\def \non{{\nonumber}}

\def\RR{{\mathbb R}}  

\def\SS{{\mathbb S}}

\def\qed{\vbox{\hrule\hbox{\vrule height 1.5 ex\kern 1 ex\vrule}\hrule}}

\def\P{{\mathsf P}} 
\def\Q{{\mathsf Q}}
\def\E{{\mathsf E}} 

\def\ZZ{{\mathbb Z}}       
\def\NN{{\mathbb N}}       

\def\eps{\varepsilon}

\def\a{\alpha}

\def\p{\varphi}

\def\o{\omega}
\def\O{\Omega}

\def\la{\langle}
\def\ra{\rangle}

\def\part{\partial}

\begin{document}

\title{Flows, coalescence and noise }
\author{Yves Le Jan and Olivier Raimond}
\date{}
\maketitle
\begin{center}Universit\'e Paris-Sud\ and Universit\'e Paris Nanterre\\e-mail: yves.lejan@math.u-psud.fr and oraimond@parisnanterre.fr\end{center}

\vskip50pt
\noindent{\bf Summary.} We are interested in stationary ``fluid''
random evolutions with independent increments. Under some mild
assumptions, we show they are solutions of a stochastic differential
equation (SDE). There are situations where these evolutions are not
described by flows of diffeomorphisms, but by coalescing flows or by
flows of probability kernels.

In an intermediate phase, for which there exists a coalescing flow and
a flow of kernels solution of the SDE, a classification is
given~: All solutions of the SDE can be obtained by filtering a
coalescing motion with respect to a sub-noise containing the Gaussian
part of its noise. Thus, the coalescing motion cannot be described by
a white noise.

\vskip50pt
\noindent{\bf Keywords~:} Stochastic differential equations, strong
solution, stochastic flow, stochastic flow of kernels, Sobolev flow, isotropic Brownian flow,
coalescing flow, noise, Feller convolution semigroup.

\medskip
\noindent{\bf 2000 MSC classifications~:} 60H10, 60H40, 60G51, 76F05.

\newpage
\tableofcontents
\bibliographystyle{apalike}

\setcounter{equation}{0}
\section*{Introduction.}
A stationary motion on the real line with independent increments is
described by a Levy process, or equivalently by a convolution
semigroup of probability measures. This naturally extends to ``rigid''
motions represented by Levy processes on Lie groups. If one
assumes the continuity of the paths, a convolution semigroup on a Lie
group $G$ is determined by an element of the Lie algebra
$\mathfrak{g}$ (the drift) and a scalar product on $\mathfrak{g}$ (the
diffusion matrix) (see for example \cite{McKean}). We call them
the local characteristics of the convolution semigroup.

We will be interested in stationary ``fluid'' random evolutions which have
independent increments. Strong solutions of stochastic
differential equations (SDEs) driven by smooth vector fields define
such evolutions. Those are of a regular type, namely 
\bdes \ita The probability that two points thrown in the fluid at the
same time and at distance $\eps$ separate at distance one in one unit
of time tends to 0 as $\eps$ tends to 0.
\itb Such points will never hit each other. \edes
Their laws can be viewed as convolution semigroups of probability
measures on the group of diffeomorphisms. 

On a compact manifold, let $V_0,V_1,\dots,V_n$ be
vector fields and $B^1,\dots,B^n$ be independent Brownian
motions. Consider the SDE
\beq dX_t = \sum_{k=1}^n V_k(X_t) \circ dB^k_t + V_0(X_t)~dt,\eeq
which equivalently can be written
\beq df(X_t)=\sum_{k=1}^n V_kf(X_t)dB^k_t+\frac{1}{2}Af(X_t)dt \eeq
for all smooth function $f$ and $Af=\sum_{k=1}^n V_k(V_kf)+V_0f$. Note
that $Af^2-2fAf=\sum_{k=1}^n (V_kf)^2$. Then, strong solutions (when they exist), as defined for
example in \cite{R-W}, of this SDE produce a flow of maps $\p_t$, such
that for all $x$, $\p_t(x)$ is a strong solution of the SDE with
$\p_0(x)=x$, which means that $\p_t$ is a function of the Brownian
paths $B^1,\dots,B^n$ up to time $t$. When the vector fields are
smooth, strong solutions are known to exist, and to be unique. The
framework can be extended to include flows of maps driven by vector
fields valued Brownian motions, which means essentially that
$n=\infty$ (see for example \cite{baxendale,iw,ku,lj,lw}). 

In a previous work \cite{ljr}, this was extended again to include flows
of Markovian operators $S_t$ solutions of the SPDE
\beq dS_tf=\sum_{k=1}^\infty
S_t(V_kf)dB^k_t+\frac{1}{2}S_t(Af)dt, \eeq
assuming the covariance function $C=\sum_{k=1}^\infty V_k\otimes V_k$
of the Brownian vector field $\sum_{k=1}^\infty V_k B^k$ is compatible
with $A$, namely that
\beq Af^2-2fAf \leq \sum_{k=1}^\infty (V_kf)^2.\eeq
Existence and uniqueness of a flow of Markovian operators $S_t$, which is a strong
solution of the previous SPDE in the sense that $S_t$ is a function of
the Brownian paths $(B^i)_{i\geq 1}$ up to time $t$ holds
under rather weak assumptions. However it is assumed in \cite{ljr}
that $A$ is self-adjoint with respect to a measure $m$ and the
Markovian operators act on $L^2(m)$ only.

The local characteristics of these flows are given by $A$ 
and the covariance function $C$, and they determine the SDE or the SPDE. 
But it was shown in \cite{ljr} that covariance functions which are not 
smooth on the diagonal (e.g. covariance associated with Sobolev norms
of order between $d/2$ and $(d+2)/2$, $d$ being the dimension of the
space) can produce strong solutions, which define random evolutions of
different type~:
\bdes \item - turbulent evolutions where {\bf (a)} is not
satisfied, which means that two points thrown initially at the same
place separate, though there is no pure diffusion, i.e. that
$Af^2-2fAf = \sum_{k=1}^\infty (V_kf)^2$.
\item - coalescing evolutions where {\bf (b)} does not hold. \edes

In this paper, we adopt a different approach based on consistent
systems of $n$ point Markovian Feller semigroups which can be viewed
as determining the law of the motion of $n$ unsecable points thrown
into the fluid. Regular and coalescing evolutions are represented by
flows of maps. Turbulent evolutions by flows of probability kernels
$K_{s,t}(x,dy)$ describing how a point mass (made of a continuum of
unsecable points) in $x$ at time $s$ is spread at time $t$. (Note that
in that case, the motion of an unsecable point is not fully determined
by the flow).

Among turbulent evolutions, we can distinguish the intermediate ones
where two points thrown in the fluid at the same place separate but
can meet after, i.e. where {\bf (a)} and {\bf (b)} are both not satisfied.

In the intermediate phase, it has been shown in \cite{Darl} (for gradient fields) and (at a
physical level) in \cite{E,VdE,G-V} that 
a coalescing solution of the SDE can be defined in law,
i.e. in the sense of the martingale problems for the $n$-point
motions. We present a construction of a coalescing flow in the 
intermediate phase. This flow obviously differs from the strong
solution $(S_{s,t},~s\leq t)$ and corresponds to an absorbing boundary
condition on the diagonal for the two-point motion.

This flow generates a vector field valued white noise $W$ and we can
identify the strong solution to the coalescing flow $(\p_{s,t},~s\leq
t)$ filtered by the velocity field $\sigma(W)$.
The noise, in Tsirelson sense (see \cite{Tsirelson}), associated to
the coalescing flow, is not linearizable, i.e. cannot be generated by
a white noise though it contains $W$.

A classification of the solutions of the SDE (or of the SPDE) can be given~: They are
obtained by filtering a coalescing motion defined on an extended
probability space with respect to a
sub-noise containing the Gaussian part of its noise.

\medskip
Let us explain in more details the contents of the paper.
We give in section \ref{sfm} and \ref{sfk} construction results, which
generalize a theorem by De Finetti on exchangeable variables (see for
example \cite{kal}). 
A stochastic flow of kernels $K$ is associated
with a general compatible family $(\P^{(n)}_t,~n\geq 1)$ of Feller
semigroups. The flow $K$ is induced by a flow of measurable mappings
when 
$$\P^{(2)}_tf^{\otimes 2}(x,x)=\P_tf^2(x),$$
for all $f\in C(M)$, $x\in M$ and $t\geq 0$. The Markov process
associated with $\P^{(n)}_t$ represents the motion of $n$ unsecable
points thrown in the fluid. The two notions are shown to be
equivalent: the law of a stochastic flow of kernels is uniquely
determined by the compatible system of $n$-point motions.
This construction is related to a recent result of Ma and Xiang
\cite{Ma} where an associated measure valued process was constructed
in a special case (the flow can actually be viewed as giving the
genealogy of this process, i.e. as its ``historical  process'') and to
a result of Darling \cite{Darl}. Note however that Darling did not get
flows of measurable maps except in very special cases. See also
Tsirelson \cite{TsirelsonFlour} for an alternative approach to this
construction.

In section \ref{secnoise}, we define the noise associated with $K$ and
introduce the notion of ``filtering with respect to a sub-noise''.

In section \ref{coal}, coalescing flows are constructed and briefly
studied. They can be obtained from any flow whose two-point motion
hits the diagonal. Then the original flow is shown to be recovered by
filtering the coalescing flow with respect to a sub-noise.

In section \ref{secW} we restrict our attention to diffusion
generators. We define the vector field valued white noise $W$
associated with the stochastic flow of kernels $K$ and prove that the
flow solves the SDE driven by the white noise $W$.

In section \ref{noise}, under some off diagonal uniqueness assumption
for the law of the $n$-point motion, we show there is only one strong
solution of the SDE. In the intermediate phase described above, the
classification of other solutions by filtering of the coalescing
solution is established. Then we identify the linear part of the noise
generated by these solutions to the noise generated by $W$.

The examples related to our previous work (see \cite{ljr}) are presented
in section \ref{isotrope}, with an emphasis on the verification of the
Feller property for the semigroups $\P^{(n)}_t$, the classification
of the solutions and the appearance of non-classical noise,
i.e. predictable noises which cannot be generated by white noises.

\renewcommand {\theequation}{\arabic{section}.\arabic{equation}}
\setcounter{equation}{0}
\section{Stochastic flow of measurable mappings.}\label{sfm}  
\subsection{Compatible family of Feller semigroups.}\label{hypsg}
Let $M$ be a separable compact metric space and $d$ a distance on $M$.

\begin{definition} Let $(\P^{(n)}_t,~n\geq 1)$ be a family of Feller
semigroups, respectively defined on $M^n$ and acting on
$C(M^n)$. We say that this family is compatible as soon as for all
$k\leq n$,
\beq \P^{(k)}_tf(x_1,\dots,x_k)=\P^{(n)}_tg(y_1,\dots,y_n) \eeq 
where $f$ and $g$ are any continuous functions such that
\beq g(y_1,\dots,y_n)=f(y_{i_1},\dots,y_{i_k}) \eeq with
$\{i_1,\dots,i_k\}\subset\{1,\dots,n\}$ and
$(x_1,\dots,x_k)=(y_{i_1},\dots,y_{i_k})$. 

We will denote by $\P^{(n)}_{(x_1,\dots,x_n)}$ the law of the Markov
process associated with $\P^{(n)}_t$ starting from
$(x_1,\dots,x_n)$. This Markov process will be called the $n$-point
motion of this family of semigroups. It is defined on the set of
c\`adl\`ag paths on $M^n$. \end{definition} 

\brem
$\P^{(n)}_t$ is a Feller semigroup on $M^n$ if
and only if $\P^{(n)}_t$ is positive (i.e. $\P^{(n)}_tf\geq 0$ for all
$f\geq 0$), $\P^{(n)}_t1=1$ and for all continuous function $f$,
$\lim_{t\to 0}\P^{(n)}_t f(x)=f(x)$ which implies the uniform
convergence of $\P^{(n)}_t f$ towards $f$ (see theorem 9.4 in chapter
I of \cite{BluGet}).
\erem

\subsection{Convolution semigroups on the space of measurable mappings.}
We equip $M$ with its Borel $\sigma$-field $\cB(M)$.
Let $(F,\cF)$ be the space of measurable
mappings on $M$ equipped with the $\sigma$-field generated by the
mappings $\p\mapsto \p(x)$ for all $x\in M$.

\begin{definition}\label{def:measpres} A probability measure $\Q$ on $(F,\cF)$ is called
regular if there exists a measurable mapping $\cJ~:~(F,\cF)\to
(F,\cF)$ such that
$$ \begin{array}{lll}
(M\times F,\cB(M)\otimes \cF) &\to& (M,\cB(M))\\
(x,\p) &\mapsto& \cJ(\p)(x)
\end{array}$$
is measurable and for all $x\in M$,
\beq \Q(d\p)-a.s.,\qquad \cJ(\p)(x)=\p(x),\eeq
i.e. $\cJ$ is a measurable modification of the identity mapping
on $(F,\cF,\Q)$. We call it a measurable presentation of
$\Q$. \end{definition}

\brem \label{rem:12}{\em If $\cJ$ and $\cJ'$ are two measurable presentations of a regular probability measure $\Q$ and if $\mu\in \mathcal{P}(M)$, then (using Fubini's theorem), 
$$\mu(dx)\otimes \Q(d\varphi)-a.s., \qquad \cJ(\p)(x)=\cJ'(\p)(x).$$
This remark will be used to prove (\ref{eq:1.4}) below.}\erem

\bprop \label{compos} Let $\Q_1$ and $\Q_2$ be two probability
measures on $(F,\cF)$. Assume $\Q_1$ is regular. Let $\cJ$ be a
measurable presentation of $\Q_1$. Then the mapping
$$\begin{array}{lll} (F^2,\cF^{\otimes 2}) & \to & (F,\cF)\\
 (\p_1,\p_2) &\mapsto& \cJ(\p_1)\circ\p_2 \end{array}$$
is measurable. Moreover, if $\cJ'$ is another
measurable presentation of $\Q_1$, then for all $x\in M$
\beq \label{eq:1.4} \Q_1(d\p_1)\otimes \Q_2(d\p_2)-a.s.,\qquad 
\cJ(\p_1)\circ \p_2(x)=\cJ'(\p_1)\circ \p_2(x).\eeq
\eprop
\brem  {\em \bdes \iti $(\p_1,\p_2)\mapsto \cJ(\p_1)\circ \p_2$ is
measurable but $(\p_1,\p_2)\mapsto \p_1\circ \p_2$ is not measurable.
\itii The law of $\cJ(\p_1)\circ \p_2$ does not
depend of the chosen presentation $\cJ$. \edes} \erem
\noindent{\bf Proof of proposition \ref{compos}.} Let $\cJ$ be a
measurable presentation of $\Q_1$. For all $x\in
M$, the mapping $(\p_1,\p_2) \mapsto \cJ(\p_1)\circ\p_2(x)$ is
measurable since it is the composition of the measurable mappings
$(\p_1,\p_2)\mapsto (\p_1,\p_2(x))$ and $(\p_1,y)\mapsto
\cJ(\p_1)(y)$. By definition of $\cF$, the mapping  $(\p_1,\p_2)
\mapsto \cJ(\p_1)\circ\p_2$ is measurable.

For $x\in M$, we have
$$\Q_1(d\p_1)-a.s.,\qquad \cJ(\p_1)(x)=\p_1(x).$$
Thus, for all $x\in M$ and all $\p_2\in F$,
$$\Q_1(d\p_1)-a.s.,\qquad \cJ(\p_1)\circ\p_2(x)
=\p_1\circ \p_2(x) = \cJ'(\p_1)\circ\p_2(x) .$$
Therefore, using Fubini's theorem,
$$\Q_1(d\p_1)\otimes\Q_2(d\p_2)-a.s.,\qquad
\cJ(\p_1)\circ\p_2(x)=\cJ'(\p_1)\circ \p_2(x). \quad \qed $$

\begin{definition} We denote $\Q_1 * \Q_2$, and we call the
convolution product of $\Q_1$ and $\Q_2$, the law of the random
variable $(\p_1,\p_2)\mapsto \cJ(\p_1)\circ \p_2$ defined on the
probability space $(F^2,\cF^{\otimes 2},\Q_1\otimes \Q_2)$.
\end{definition}

\begin{definition} A convolution semigroup on $(F,\cF)$ is a family
$(\Q_t)_{t\geq 0}$ of regular probability measures on $(F,\cF)$ such
that for all nonnegative $s$ and $t$, $\Q_{s+t}=\Q_s * \Q_t$.
\end{definition}

\begin{definition} \label{defconv}
A convolution semigroup $(\Q_t)_{t\geq 0}$ on $(F,\cF)$ is called
Feller if
\bdes \iti $\forall f\in C(M)$, $\lim_{t\to 0} \sup_{x\in M}
\int(f\circ \p(x)-f(x))^2\Q_t(d\p) = 0$.
\itii $\forall f\in C(M)$, $\forall t\geq 0$, $\lim_{d(x,y)\to 0} 
\int(f\circ \p(x)-f\circ \p(y))^2\Q_t(d\p) = 0$. \edes \end{definition}

\bprop \label{QP}Let $(\Q_t)_{t\geq 0}$  be a Feller convolution semigroup on
$(F,\cF)$. For all $n\geq 1$, $f\in C(M^n)$ and $x\in M^n$, set
\beq \label{eq:p(n)m}\P^{(n)}_tf(x)=\int f\circ \p^{\otimes n}(x)~\Q_t(d\p).\eeq
Then $(\P^{(n)}_t,~n\geq 1)$ is a compatible family of Feller
semigroups on $M$ satisfying 
\beq \label{mapcond}\P^{(2)}_tf^{\otimes 2}(x,x)=\P_tf^2(x),\eeq
for all $f\in C(M)$, $x\in M$ and $t\geq 0$. \eprop
\prf It is easy to see that this family is compatible and that for all
$n\geq 1$ and all $t\geq 0$, $\P^{(n)}_t$ is Markovian. Let $s$ and
$t$ in $\RR^+$, $f\in C(M^n)$ and $x\in M^n$, then
\beqarr
\P^{(n)}_{s+t}f(x)
&=& \int f\circ \p^{\otimes n}(x)~\Q_{s+t}(d\p)\\
&=& \int f\circ \cJ(\p_1)^{\otimes n} \circ 
\p_2^{\otimes n}(x)~\Q_t(d\p_1)\otimes \Q_s(d\p_2)\\
&=& \int \P^{(n)}_t f \circ \p_2^{\otimes n} (x)~\Q_s(d\p_2)\\
&=& \P^{(n)}_{s}\P^{(n)}_{t}f(x) \eeqarr
where $\cJ$ is a measurable presentation of $\Q_t$. This proves that
$\P^{(n)}_t$ is a semigroup.

\smallskip
Let us now prove the Feller property. Let $h\in C(M^n)$ be in the
form $f_1\otimes\cdots\otimes f_n$, $x=(x_1,\dots,x_n)$ and
$y=(y_1,\dots,y_n)$. We have for $M$ large enough 
\beq |\P^{(n)}_t h(y) - \P^{(n)}_t h(x)| \leq M\sum_{k=1}^n
\left(\int (f_k\circ\p(y_k)-
f_k\circ \p(x_k))^2\Q_t(d\p)\right)^{\frac{1}{2}} \eeq
which converges towards 0 as $d(x,y)$ goes to 0 since {\bf (ii)} in
definition \ref{defconv} is satisfied. We also have
\beq |\P^{(n)}_t h(x) - h(x)| \leq M\sum_{k=1}^n
\left(\int(f_k\circ \p(x_k)-f_k(x_k))^2 
\Q_t(d\p)\right)^\frac{1}{2} \eeq
which converges towards 0 as $t$ goes to 0 since {\bf (i)} in
definition \ref{defconv} is satisfied. These properties
extend to all function $h$ in $C(M^n)$ by an approximation
argument. This proves the Feller property of the Markovian semigroups
$\P^{(n)}_t$.

\smallskip
It remains to prove (\ref{mapcond}). This follows from
\beqarr
\P^{(2)}_tf^{\otimes 2}(x,x) 
&=& \int f^{\otimes 2}\circ \p^{\otimes 2} (x,x)~\Q_t(d\p)\\
&=& \int f^{2}\circ \p(x) ~\Q_t(d\p) 
\quad=\quad \P^{(1)}_t f^2(x). \qquad \qed
\eeqarr

\brem
{\em Let $(\Q_t)_{t\geq 0}$  be a Feller convolution semigroup on $(F,\cF)$.
\begin{itemize}
\item The semigroup $(\Q_t)_{t\geq 0}$ is uniquely determined by
$(\P^{(n)}_t,~n\geq 1)$.
\item Let $X$ and $\p$ be independent random variables respectively in $M^n$ and in $F$.
Denote by $\mu$ the law of $X$ and suppose that the law of $\p$ is $\Q_t$. 
Then Fubini's theorem implies that for all  measurable presentation $\cJ$ of $\Q_t$, the random variable $\cJ(\p)^{\otimes n}(X)$ is distributed as $\mu \P^{(n)}_t$, where $\P^{(n)}_t$ is defined by \eqref{eq:p(n)m}.
\item For all $x\in M$, $\Q_0(d\p)$-almost surely, $\p(x)=x$.
\end{itemize}
}\erem

\subsection{Stochastic flows of mappings.}
\begin{definition}\label{defsfm} 
Let $(\Omega,{\cal A},\P)$ be a
probability space and let $\p=(\p_{s,t},~s\leq t)$ be a family of $(F,\cF)$-valued random variables such that for all $x\in M$ and all $t\in\mathbb{R}$, $\P$-a.s. $\p_{t,t}(x)=x$. For $t\ge 0$, denote by $\Q_t$ the law of $\p_{0,t}$.
The family $\p$ is called a stochastic flow of mappings if for all $t\ge 0$, $\Q_t$ is regular and if the following properties are satisfied by $\p$ 
\bdes \ita For all $s\le u\le t$, all $x\in M$ and all measurable presentation $\cJ_{t-u}$ of $\Q_{t-u}$,\\ $\P$-almost surely,
$\p_{s,t}(x)=\cJ_{t-u}(\p_{u,t})\circ\p_{s,u}(x)$. (cocycle  property)
\itb For all $s\leq t$, the law of  $\p_{s,t}$ is $\Q_{t-s}$. (Stationarity)
\itc The flow has independent increments, i.e. for all
$t_1\le t_2\le \cdots\le t_n$, the family $\{\p_{t_i,t_{i+1}},~1\leq i\leq n-1\}$
is independent.
\itd For all $f\in C(M)$ and all $s\leq t$, $\lim_{(u,v)\to (s,t)}\sup_{x\in M}
\E[(f\circ \p_{s,t}(x)-f\circ\p_{u,v}(x)) ^2]=0$.
\ite For all $f\in C(M)$ and all $s\leq t$,
$\lim_{d(x,y)\to 0}\E[(f\circ\p_{s,t}(x)-f\circ\p_{s,t}(y))^2]=0$.
\edes
\end{definition}
\brem
\begin{itemize}
\item Item (a) holds for all set of measurable presentations as soon as (a) holds for one of them.
\item If $\psi$ is equal in law to a stochastic flow of mappings $\p$, then $\psi$ is also a stochastic flow of mappings. Indeed, it is straightforward to check that $\psi$ satisfies (b), (c), (d) and (e), and after having remarked that for all $x\in M$, $(\p_1,\p_2,\p_3)\mapsto (\p_3(x),\cJ_{t-u}(\p_2)\circ \p_1(x))$ is measurable, we prove that $\psi$ satisfies (a).
\end{itemize}
\erem
\bprop \label{prop:modppres}
Let $\p$ be a stochastic flow of mappings, and for $t\ge 0$, let $\cJ_t$ be a measurable presentation of the law of $\p_{0,t}$. Then $\p'=(\cJ_{t-s}(\p_{s,t}),\, s\le t)$ is a stochastic flow of mappings satisfying
\begin{itemize}
\item[(i)] For all $s\le t$ and all $x\in M$, a.s. $\p'_{s,t}(x)=\p_{s,t}(x)$.
\item[(ii)] For all $s\le u\le t$ and all $x\in M$, $\P$-almost surely,
$\p'_{s,t}(x)=\p'_{u,t}\circ\p'_{s,u}(x)$.
\end{itemize}
\eprop 
\prf Item (i) is a consequence of the fact that $\cJ_{t-s}$ is a measurable presentation of the law of $\p_{s,t}$. Then $\p$ and $\p'$ share the same law and $\p'$ is a stochastic flow of mappings.
Let us now prove (ii). For $s\le u\le t$ and $x\in M$, it holds that $\P$-almost surely,
\begin{align*}
\p'_{s,t}(x)&=\cJ_{t-u}(\p'_{u,t})\circ\p'_{s,u}(x)\\
&= \cJ_{t-u}\circ \cJ_{t-u}(\p_{u,t})\circ\p'_{s,u}(x).
\end{align*}
Since $\cJ_{t-u}\circ \cJ_{t-u}$ is also a measurable presentation of $\Q_{t-u}$, using remark \ref{rem:12}, $\P$-almost surely,
\begin{align*}
\cJ_{t-u}\circ \cJ_{t-u}(\p_{u,t})\circ\p'_{s,u}(x)=\cJ_{t-u}(\p_{u,t})\circ\p'_{s,u}(x)= \p'_{u,t}\circ\p'_{s,u}(x) .
\end{align*}
This proves (ii). \qed\begin{definition}
\begin{itemize}
\item The stochastic flow of mappings $\p'$ defined in proposition \ref{prop:modppres} will be called a measurable modification of $\p$.
\item  A stochastic flow of mappings which is a measurable modification of a stochastic flow of mappings is called a measurable stochastic flow of mappings.
\end{itemize}
\end{definition}
\brem \label{rem:15} {\em In the proof of theorem \ref{mainthm} below, the measurable presentations we construct satisfy $\cJ_t\circ\cJ_t=\cJ_t$ for all $t\ge 0$, and the measurable stochastic flow $\p$ we construct satisfies $\cJ_{t-s}(\p_{s,t})=\p_{s,t}$ for all $s\le t$. 
But, we will see that this property doesn't hold for stochastic flows of kernels studied in section 2.}
\erem

\bprop \label{SFMP}
Let $\p=(\p_{s,t},~s\leq t)$  be a stochastic flow of
mappings. For all $n\geq 1$, $f\in C(M^n)$ and $x\in M^n$, set 
\beq \P^{(n)}_tf(x)=\E[f\circ \p_{0,t}^{\otimes n}(x)].\eeq
Then $(\P^{(n)}_t,~n\geq 1)$ is a compatible family of Feller
semigroups on $M$ satisfying (\ref{mapcond}). \eprop
\prf For $t\ge 0$, denote by $\Q_t$ the law of $\p_{0,t}$. Then $\Q_t$ is regular and there is $\cJ_t$ a measurable presentation of $\Q_t$. With (a), (b) and (c), we show that for all $(s,t)\in\mathbb{R}_+^2$,  $\Q_s*\Q_t=\Q_{s+t}$, i.e. $(\Q_t)_{t\ge 0}$ is a convolution semigroup. Finally (d) and (e) imply that it is Feller. To conclude, we apply proposition \ref{QP}. \qed

\subsection{Construction and characterization.}\label{sec:conschar}
In this section, we present a theorem stating that to any compatible
family $(\P^{(n)}_t,~n\geq 1)$ of Feller semigroups, one can associate
a Feller convolution semigroup on $(F,\cF)$ and a stochastic
flow of mappings.

Let $(\Omega^0,{\cal A}^0)$ denote the measurable space $(\prod_{s\leq t}F ,
\otimes_{s\leq t} \cF)$. For $s\leq t$, let $\p^0_{s,t}$ denote
the random variable $\omega\mapsto \omega(s,t)$. Let  $\p^0$ be the
random variable $(\p^0_{s,t},~s\leq t)$. Then $\p^0(\omega)=\omega$.
Let $(T_h)_{h\in\RR}$ be the one-parametric group of transformations
of $\Omega^0$ defined by $T_h(\omega)(s,t)=\omega(s+h,t+h)$, for all
$s\leq t$, $h\in \RR$ and $\omega\in\Omega^0$. 

\begin{definition} A probability space $(\Omega,{\cal A},\P)$ is said
separable if the Hilbert space $L^2(\Omega,{\cal A},\P)$ is
separable. 
(Note that this implies that for all $1\leq p <\infty$, $L^p(\Omega,{\cal A},\P)$ is
separable.)
\end{definition}

\bthm\label{mainthm} {\em (i)} Let $(\P^{(n)}_t,~n\geq 1)$ be a
compatible family of Feller semigroups on $M$ satisfying
\beq \P^{(2)}_tf^{\otimes 2}(x,x)=\P_tf^2(x),\eeq
for all $f\in C(M)$, $x\in M$ and $t\geq 0$. Then
there exists a unique Feller convolution semigroup $(\Q_t)_{t\geq 0}$
on $(F,\cF)$ such that for all $n\geq 1$, $t\geq 0$, $f\in C(M^n)$
and $x\in M^n$, 
\beq \label{relconv} \P^{(n)}_tf(x)=\int f\circ\p^{\otimes n}(x)~\Q_t(d\p). \eeq

\medskip
{\em (ii)} For all Feller convolution semigroup $\Q=(\Q_t)_{t\geq 0}$
on $(F,\cF)$, there exists a unique $(T_h)_{h\in\RR}$-invariant
probability measure $\P_\Q$ on $(\Omega^0,{\cal A}^0)$ such that
$(\Omega^0,{\cal A}^0,\P_\Q)$ is separable, the family of random
variables $\p^0=(\p^0_{s,t},~s\leq t)$ is a stochastic flow of
mappings and for all $s\leq t$, the law of $\p^0_{s,t}$ is
$\Q_{t-s}$. 
Every measurable modification $\p'$
of $\p^0$ satisfies $\p'_{s+h,t+h}=\p'_{s,t}\circ T_h$ for all $s\le t$ and all $h\in \mathbb{R}$.

The flow $\p^0$ is called the canonical stochastic flow of
mappings associated with $\Q$ (or equivalently with
$(\P^{(n)}_t,~n\geq 1)$). 
\ethm

\brem  Theorem \ref{mainthm} is also satisfied when $M$ is a
locally compact separable metric space. In this case,
$(\P^{(n)}_t,~n\geq 1)$ is a compatible family of Markovian semigroups
acting continuously on $C_0(M^n)$, the set of continuous functions on
$M^n$ converging towards 0 at $\infty$ (we call them Feller
semigroups). In the previous definitions (\ref{defconv} and
\ref{defsfm}) and in the statement of the theorem the function $f$ has
to be taken in $C_0(M)$ or in $C_0(M^n)$. Moreover {\bf (ii)} of
definition \ref{defconv} must be modified by: for all $x\in M$, $f\in
C_0(M)$ and $t\geq 0$,
\beq \left\{ \begin{array}{lll} 
\lim_{y\to x} \int (f\circ\p(y)-f\circ\p(x))^2~\Q_t(d\p) &=& 0 \\
\hbox{and} \qquad \lim_{y\to \infty} \int (f\circ\p(y))^2~\Q_t(d\p) &=& 0. 
\end{array}\right. \eeq
In definition \ref{defsfm}, {\bf (e)} must be modified by: for all $x\in
M$ and $s\leq t$,
\beq \left\{ \begin{array}{lll} 
\lim_{y\to x} \E[(f\circ\p_{s,t}(y)-f\circ\p_{s,t}(x)))^2] &=& 0\\
\hbox{and} \qquad \lim_{y\to \infty} \E[(f\circ\p_{s,t}(y))^2] &=& 0
\end{array} \right. \eeq 
\erem 

\prf In order to prove this remark, note that the one-point
compactification of $M$, $\hat{M}=M\cup \{\infty\}$, is a  separable
compact metric space. On $\hat{M}$, we define the compatible family
of Feller semigroups, $({\hat{\P}}^{(n)}_t,~n\geq 1)$,  by the
following relations:

for all $n\geq 2$ and all family of continuous functions on $\hat{M}$,
$\{f_i,~i\geq 1\}$, 
\beqar
{\hat{\P}}^{(n)}_tf_1\otimes \cdots \otimes f_n
&=& \P^{(n)}_t g_1\otimes \cdots \otimes g_n\\
&+& \sum_{i=1}^n f_i(\infty){\hat{\P}}^{(n-1)}_tf_1\otimes \cdots
\otimes f_{i-1}\otimes g_{i+1} \otimes \cdots \otimes
g_n\non \eeqar
and
\beq {\hat{\P}}^{(1)}_t f_1 = f_1(\infty) + \P^{(1)}_tg_1, \eeq
where $g_i=f_i-f_i(\infty)\in C_0(M)$ and with the convention 
$\P^{(n)}_tg_1\otimes \cdots \otimes g_n(x_1,\dots,x_n)=0$ if there
exists $i$ such that $x_i=\infty$. We apply theorem \ref{mainthm} to
$\hat{M}$ and to the family $(\hat{\P}^{(n)}_t,~n\geq 1)$ to construct
a Feller convolution semigroup $\hat{\Q}$ and a stochastic flow of
mappings $(\hat{\p}_{s,t},~s\leq t)$ on $\hat{M}$. This stochastic flow
of mappings satisfies 
\bdes
\iti $\hat{\p}_{s,t}(\infty)=\infty$ for all $s\leq t$ and
\itii $\hat{\p}_{s,t}(x)\neq \infty$ for all $x\in M$ and $s\leq t$.
\edes
\noindent{\bf Proof of (i).} For all $f\in C(\hat{M})$,
\beqarr \E[(f\circ\hat{\p}_{s,t}(\infty)-f(\infty))^2] &=&
\hat{\P}^{(2)}_{t-s}f^{\otimes 2}(\infty,\infty) -
2f(\infty)\hat{\P}^{(1)}_{t-s}f(\infty)+f(\infty)^2\\
&=& 0 \eeqarr
since $\hat{\P}^{(2)}_{t-s}f^{\otimes 2}(\infty,\infty)=f(\infty)^2$
and $\hat{\P}^{(1)}_{t-s}f(\infty)=f(\infty)$. This implies {\bf
(i)}. \qed 

\smallskip
\noindent{\bf Proof of (ii).} Let $g_n$ be a sequence in $C_0(M)$ such
that $g_n\in [0,1]$ and simply converging towards 1. Then
$f_n=1-g_n\in C(\hat{M})$ is such that $f_n(\infty)=1$ and for all
$x\in M$
$$\E[(f_n\circ\hat{\p}_{s,t}(x))^2] = \hat{\P}^{(2)}_{t-s}g_n^{\otimes 2}(x,x)
+1 - 2\hat{\P}^{(1)}_{t-s}g_n(x).$$
This implies that
$\lim_{n\to\infty}\E[(f_n\circ\hat{\p}_{s,t}(x))^2]=0$. Assertion {\bf (ii)}
follows since $1_{\{\hat{\p}_{s,t}(x)=\infty\}} =
\lim_{n\to\infty}f_n\circ\hat{\p}_{s,t}(x)$. \qed

\medskip
For all $x\in M$, let us denote $\hat{\p}_{s,t}(x)$ by $\p_{s,t}(x)$.
Assertions {\bf (i)} and {\bf (ii)} implies that $\p_{s,t}\in F$ and
that $(\p_{s,t},~s\leq t)$ is a stochastic flow of mappings on $M$. In
a similar way, one can show that $\hat{\Q}$ induces a Feller
convolution semigroup on $(F,\cF)$. \qed 

\medskip
Let us explain briefly the method we use to prove theorem
\ref{mainthm}. We first suppose we are given a compatible family of
Feller semigroups satisfying (\ref{mapcond}). Then we define a
convolution semigroup  $(\Q_t,~t\geq 0)$ on measurable mappings on
$M$. For all $t$, to define $\Q_t$, we define $\P^{(\infty)}_t$, the law
of $(\p(z_l),~l\in\NN)$, where the law of $\p$ is $\Q_t$, for some dense
family $(z_l,~l\in\NN)$ in $M$ and get $\Q_t$ by an
approximation. Hence $\Q_t$ is defined as the law of a random
variable, which takes its values in the ``bad'' space $F$,  but is
defined on a ``nice'' space $M^\NN$.

The approximation used to construct this convolution semigroup
allows us to define a stochastic flow of mappings on $M$ in
such a way that these mappings are measurable, defining it first on
the dyadic numbers. We get a measurable flow defined on a ``nice''
space. Note that a difficulty to get this measurability comes from the
fact that the composition of mappings from $M$ onto $M$ is not
measurable with respect to the natural $\sigma$-field.

\subsection{Proof of the first part of theorem \ref{mainthm}.}
In the following we assume we are given $(\P^{(n)}_t,~n\geq 1)$, a
compatible family of Feller semigroups satisfying (\ref{mapcond}). 
And we intend to construct a Feller convolution semigroup
$(\Q_t)_{t\geq 0}$ on $(F,\cF)$ satisfying (\ref{relconv}). The
uniqueness of such a convolution semigroup is immediate since
(\ref{relconv}) characterizes $\Q_t$.

\subsubsection{A measurable choice of limit points in $M$.}
It is known that, as a separable compact metric space, $M$ is
homeomorphic to a closed subset of $[0,1]^\NN$ (see corollaire 1 �6.1
of chapter 9 in \cite{bbki}). A point $y$ can be represented by a
sequence $(y^n)_{n\in\NN}\in [0,1]^\NN$. Let $y=(y_i)_{i\in\NN}$ be a
sequence of elements of $M$.

Let $y^1=\limsup_{i\to\infty} y^1_i$. Let $i^1_k=\inf\{i,~|y^1-y^1_i|<1/k\}$.
By induction, for all integer $j$, we construct $y^j$ and
$\{i_k^j,~k\in\NN\}$ by the relations
$$y^j=\limsup_{k\to\infty} y^j_{i^{j-1}_k}
\quad\hbox{and}\quad i^j_k = 
\inf\{i\in\{i^{j-1}_k,~k\in\NN\},~|y^j-y^j_i|<1/k\}.$$
We denote $(y^n)_{n\in\NN}$ by $l(y)$. Note that
$l(y)^j=\lim_{n\to\infty}y_{i_n^n}^j$. Hence
$l(y)$ belongs to $M$.
It is easy to see that $l$ satisfies the following lemma.

\blem\label{lmes} $l:M^\NN\to M$ is a measurable
mapping, $M$ being equipped with the Borel $\sigma$-field
${\cal B}(M)$ and $M^\NN$ with the product
$\sigma$-field ${\cal B}(M)^{\otimes\NN}$. Moreover
$l((y_i)_{i\in\NN})=y_\infty$ when $y_i$ converges towards $y_\infty$. 
\elem

\subsubsection{Notation and definitions.}
Let $\{z_l,~l\in \NN\}$ be a dense family in $M$, which
will be fixed in the following. We wish to define a measurable mapping
$i~:~M^\NN \to F$ such that
$i((y_j)_{j\in\NN})(z_l)=y_l$ for all integer $l$.  

Let $(\eps_k)_{k\in\NN}$ be a positive sequence decreasing
towards 0 (this sequence will be fixed later). Let $i~:~M^\NN\to F$ be
the injective mapping defined by
\beq\label{defi}i(y)(x) = l((y_{n_k^x})_{k\in\NN}) \eeq
where
\beq \label{defn}  n_k^x = \inf\{n,~d(z_n,x)\leq
\eps_k\},\eeq
for $(y,x)\in M^\NN\times M$. Note that $i(y)$
defined this way is a measurable mapping since $l$ is measurable and
$x\mapsto (y_{n_k^x})_{k\in\NN}$ is measurable. Note also that the
relation $i(y)(z_l)=y_l$ is satisfied for all integer $l$. 
\blem \label{lem:12} For $n\geq 1$, the mappings $\Phi_n:(M^\NN)^n\to F$ and
$\Psi_n:M\times (M^\NN)^n\to M$, defined by
$$\begin{array}{l} \Phi_n(y^1,\dots,y^n) \quad=\quad 
i(y^n)\circ i(y^{n-1})\circ\cdots\circ i(y^1)\\
\Psi_n(x,y^1,\dots,y^n) \quad=\quad \Phi_n(y^1,\dots,y^n)(x)
\end{array}$$
are measurable. ($(M^\NN)^n$ and $M\times (M^\NN)^n$ are equipped with
the product $\sigma$-field.) In particular, $i$ is measurable.
\elem
\prf  Note that $\Psi_1$ is the composition of the mappings $l$ and
$(x,y)\mapsto (y_{n_k^x})_{k\in\NN}$. Since these mappings are
measurable, $\Psi_1$ is measurable. 
By induction, we prove that $\Psi_n$ is measurable since, for $n\geq 2$,
$$\Psi_n(x,y^1,\dots,y^n)=\Psi_1(\Psi_{n-1}(x,y^1,\dots,y^{n-1}),y^n).$$

For all $A\in \cB(M)$ and $x\in M$,
$$\Phi_n^{-1}(\{\p\in F,~\p(x)\in A\})=\{y\in
(M^\NN)^n,~(x,y)\in \Psi_n^{-1}(A)\}.$$
 This event belongs to $({\cal B(M)}^{\otimes \NN})^{\otimes
n}$ since $\Psi_n$ is measurable. This shows the measurability of
$\Phi_n$. \qed

\medskip We need to introduce $\Phi_n$ because the composition
application $F^n\to F$, $(\p_1,\dots,\p_n)\mapsto
\p_n\circ\cdots\circ\p_1$ is not $\cF^{\otimes n}$-measurable in
general.

\medskip
Let $j:F\to M^{\NN}$ be the mapping defined by 
\beq j(\p)=(\p(z_l))_{l\in\NN}.\eeq
\blem \label{lem:13} The mapping $j$ is measurable and satisfies $j\circ i(y)=y$ for all $y\in M^\NN$.\elem
\prf We have for all $A\in \cB(M)^{\otimes n}$,
$$j^{-1}(\{y\in M^\NN,~(y_1,\dots,y_n)\in A\}) =
\{\p\in F,~(\p(z_1),\dots,\p(z_n))\in A\}.$$ 
 This set belongs to $\cF$. \qed

\medskip
Note that for all $l\in\NN$ and $\p\in F$, $i\circ j(\p)(z_l)=\p(z_l)$. 

\brem {\em Set $\cJ=i\circ j$. Lemmas \ref{lem:12} and  \ref{lem:13} imply that the mappings $(\p_1,\dots,\p_n)\mapsto \cJ(\p_n)\circ \cdots \circ \cJ(\p_1)$ and $(x,\p_1,\dots,\p_n)\mapsto \cJ(\p_n)\circ \cdots \circ\cJ(\p_1)(x)$ are measurable.} \erem

\subsubsection{Constructions of probabilities on $M^\NN$ and on $F$.} \label{defprob}
By Kolmogorov's theorem, we construct on $M^\NN$ a
probability measure $\P^{(\infty)}_{t}$ such that
$\P^{(\infty)}_t(A\times M^\NN) = \P^{(n)}_t
1_A(z_1,\dots,z_n)$ for any $A\in {\cal B}(M)^{\otimes n}$. 
We now prove useful lemmas satisfied by $\P^{(\infty)}_t$~:
\blem\label{baslem}
For all positive $T$, there exists a positive function $\eps_T(r)$
converging towards 0 as $r$ goes to 0 such that
\beq \sup_{t\in [0,T]} \E^{(2)}_{(x,y)} [(d(X_t,Y_t))^2]
\leq \eps_T(d(x,y)). \eeq \elem
\prf For all continuous function $f$, we have
\beqarr
\E^{(2)}_{(x,y)}[(f(X_t)-f(Y_t))^2]
&=&\P_tf^{2}(x)+\P_tf^{2}(y)-2\P^{(2)}_tf^{\otimes 2}(x,y).\\
&=&\P^{(2)}_tf^{\otimes 2}(x,x)+\P^{(2)}_tf^{\otimes
2}(y,y)-2\P^{(2)}_tf^{\otimes 2}(x,y),
\eeqarr
since (\ref{mapcond}) is satisfied. Let $(f_n)_{n\geq 1}$ be a dense
sequence in $\{f\in C(M),~\|f\|_\infty\leq 1\}$. Then
$d'(x,y)=\left(\sum_{n\geq 1}
2^{-n}(f_n(x)-f_n(y))^2\right)^{\frac{1}{2}}$ is a distance equivalent
to $d$ and we have
$$\E^{(2)}_{(x,y)}[(d'(X_t,Y_t))^2] = \P^{(2)}_t h(x,x) + 
\P^{(2)}_th(y,y)-2\P^{(2)}_t h(x,y),$$
where $h$ is the continuous function $\sum_{n\geq 1} 2^{-n} f_n\otimes
f_n$. We conclude the lemma after remarking that this function is
uniformly continuous in $(t,x,y)$ on $[0,T]\times M^2$. \qed

\medskip
From now on we fix $T$ and define the sequence $(\eps_k)_{k\in\NN}$ (which
defines the sequence $(n_k^x)_{k\in\NN}$ for all $x\in M$ by
equation (\ref{defn})) such that $0\leq r\leq 2\eps_k$ implies
$\eps_T(r)\leq 2^{-3k}$. The sequence $(\eps_k)_{k\in\NN}$ is well
defined since $\lim_{r\to 0}\eps_T(r)=0$. Since $i$ depends on $T$, we
now denote $i$ by $i_T$, $\Phi_n$ by $\Phi_n^T$ and $\Psi_n$ by $\Psi_n^T$.

\brem {\em One can construct $i_t$ such that $i_t=i_1$ for all $t\le 1$.}\erem

\blem\label{lemps}
For all $t\in [0,T]$ and 
for any independent random variables $X$ and $Y$ respectively in
$M$ and  $M^\NN$, such that the law of $Y$ is
$\P^{(\infty)}_t$,  then $Y_{n_k^X}$ converges almost surely towards
$l((Y_{n_k^X})_{k\in\NN})=i_T(Y)(X)$ as $k$ goes to $\infty$. \elem
\prf Note that $(Y_{n_k^X})_{k\in\NN}$ is a random variable (the mapping
$(x,y)\mapsto(y_{n_k^x})_{k\in\NN}$ is measurable). For all
integer $k$, $d(z_{n_k^x},z_{n_{k+1}^x})\leq 2\eps_k$ and
\beq\label{2-k} \P[d(Y_{n_k^X},Y_{n_{k+1}^X})>2^{-k}]\leq
2^{2k}\E[\eps_T(d(z_{n_k^X},z_{n_{k+1}^X}))]\leq 2^{-k}. \eeq
Using Borel-Cantelli's lemma, we prove that almost surely,
$(Y_{n_k^X})_{k\in\NN}$ is a Cauchy sequence and therefore
converges. Its limit can only be $l((Y_{n_k^X})_{k\in\NN})$. \qed

\blem\label{lemprob} Let $(X_n)_{n\in\NN}$ be a sequence of random variables
in $M$ converging in probability towards a random variable $X$. Let $Y$
a random variable in $M^\NN$ of law $\P^{(\infty)}_t$
independent of $(X_n)_{n\in\NN}$. Then $i_T(Y)(X_n)=l((Y_{n_k^{X_n}})_{k\in\NN})$
converges in probability towards $i_T(Y)(X)=l((Y_{n_k^{X}})_{k\in\NN})$ as $n$
tends to $\infty$.
\elem
\prf Let $Z_n=l((Y_{n_k^{X_n}})_{k\in\NN})$ and
$Z=l((Y_{n_k^{X}})_{k\in\NN})$. For all integer $k$, we have
\beqarr
\P[d(Z_n,Z)>\eps]&\leq&
\P[d(Z_n,Y_{n_k^{X_n}})>\eps/3]+\P[d(Y_{n_k^{X_n}},Y_{n_k^{X}})>\eps/3]\\
&&+\quad\P[d(Y_{n_k^{X}},Z)>\eps/3].
\eeqarr
Lemma \ref{lemps} implies that the first and last terms of the right
hand side of the preceding equation converge towards 0 as $k$ goes to
$\infty$. The second term is lower than
$\frac{9}{\eps^2}\E[\eps_T(d(z_{n_k^{X_n}},z_{n_k^{X}}))]$.
Since for all positive $\a$, there exists a positive $\eta$ such that
$|r|<\eta$ implies $\frac{9}{\eps^2}|\eps_T(r)|<\a$, we get
\beqarr
\P[d(Y_{n_k^{X_n}},Y_{n_k^{X}})>\eps/3] &\leq& \a+C\P[d(z_{n_k^{X_n}},z_{n_k^{X}})>\eta]\\
&\leq& \a+C\P[d(X_n,X)>\eta-2\eps_k],
\eeqarr
where $C=9D^2/\eps^2$, where $D$ is the diameter of $M$ (one can choose
$\eps_T$ such that $\eps_T(r)\leq D^2$ for all $r$).
Therefore, we get $\P[d(Z_n,Z)>\eps]\leq \a+C\P[d(X_n,X)\geq \eta]$
and for all positive $\a$, $\limsup_{n\to\infty} \P[d(Z_n,Z)>\eps]\leq \a.$
Thus we prove that $Z_n$ converges in probability towards $Z$. \qed

\medskip
For all $t\in [0,T]$, set $\Q_t=i_T^*(\P^{(\infty)}_t)$. It is a
probability measure on $(F,\cF)$ and it satisfies the
following proposition. 

\bprop \label{unique}$\Q_t$ is the unique probability measure on
$(F,\cF)$ such that for any continuous function $f$ on
$M^n$ and any $x\in M^n$, 
\beq\label{Qn} 
\int_F f\circ \p^{\otimes n}(x)~\Q_t(d\p)
=\P^{(n)}_tf(x).\eeq
Moreover, $j^*(\Q_t)=\P^{(\infty)}_t$
and $(i_T\circ j)^*(\Q_t)=i_T^*(\P^{(\infty)}_t)=\Q_t$.
\eprop
\prf The unicity is obvious since (\ref{Qn}) characterizes $\Q_t$. Let
us check that $\Q_t=i_T^*(\P^{(\infty)}_t)$ satisfies (\ref{Qn}). Let
$Y$ be a random variable of law $\P^{(\infty)}_t$ then for all $f\in
C(M^n)$ and all $x\in M^n$,
\beqarr
\int_F f\circ\p^{\otimes n}(x)\Q_t(d\p)
&=& \E[f(i_T(Y)(x_1),\dots,i_T(Y)(x_n))]\\
&=& \lim_{k\to\infty}\E[f(Y_{n_k^{x_1}},\dots,Y_{n_k^{x_n}})]\\
&=& \lim_{k\to\infty}\P^{(n)}_tf(z_{n_k^{x_1}},\dots,z_{n_k^{x_n}})
\quad = \quad \P^{(n)}_tf(x),
\eeqarr
using first dominated convergence theorem and lemma \ref{lemps}, then
the definition of $\P^{(\infty)}_t$ and the fact that $\P^{(n)}_t$
is Feller. \qed 

\brem  Since $T$ can be taken arbitrarily large, we can define $\Q_t$
for all positive $t$ and the definition of $\Q_t$ is independent of
the chosen $T$, since $\Q_t$ satisfies proposition \ref{unique}. 
\erem

\subsubsection{A convolution semigroup on $(F,\cF)$.}
\blem\label{iojps} 
For all $t\geq 0$, $\Q_t$ is regular. And for all $T\geq t$,
$i_T\circ j$ is a measurable presentation of $\Q_t$.\elem
\prf Let $0\leq t\leq T$. For all $x\in M$ and $\p\in F$,  
$i_T\circ j(\p)(x)=\Psi^T_1(x,j(\p))$. Since $\Psi_1^T$ and $j$ are
measurable, the mapping $(x,\p)\mapsto i_T\circ j(\p)(x)$ is
measurable.

Let $x\in M$. Since $\Q_t=i_T^*(\P^{(\infty)}_t)$, if $Y$ is
a random variable of law $\P^{(\infty)}_t$,
\beqarr \Q_t[d(\p(z_{n_k^x}),\p(x))>2^{-k}] &=&
\P[d(Y_{n_k^x},i_T(Y)(x))\geq 2^{-k}]\\
&=& \lim_{l\to\infty}\P[d(Y_{n_k^x},Y_{n_l^x})\geq 2^{-k}]\leq 2^{-k}
\eeqarr
since for all $l\geq k$, $d(z_{n_k^x},z_{n_l^x})\leq 2\eps_k$ (see
equation (\ref{2-k})). Using Borel-Cantelli's lemma, we prove that
$\p(z_{n_k^x})$ converges almost surely towards $\p(x)$. Therefore,
$$\Q_t(d\p)-a.s.,\qquad i_T\circ j(\p)(x)=\p(x).$$
This proves the lemma. \qed 

\medskip In the published version of this paper \cite{FCN}, we have made the following incorrect remark (Fubini's theorem cannot be applied here since $(x,\p)\mapsto \p(x)$ is not measurable):
\brem  \label{remiojps} Let $\p$ and $X$ be independent random variables
respectively $F$-valued and $M$-valued. Then, if the law of $\p$ is
$\Q_t$ and if $M\times \Omega\ni (x,\o)\mapsto\p(x,\o)\in M$ is measurable,
Fubini's theorem implies that for all $T\geq t$,
\beq \label{eqiojps} \P-\hbox{a.s.},\qquad i_T\circ j(\p)(X)=\p(X). \eeq
\erem
\noindent\textbf{Counterexample to remark \ref{remiojps} (labelled 1.7 in \cite{FCN})} This counterexample was communicated to us by G. Riabov (see \cite{riabov}).
Let $\p$ be a random variable in $F$ of law $\Q$ such that 
$M\times \Omega\ni (x,\o)\mapsto\p(x,\o)\in M$ is measurable. 
Suppose that 
$\Q$ is regular and let $\cJ$ be a regular presentation of $\Q$. 
Let $X$ be a random variable in $M$ independent of $\p$.
Out of $\p$ and $X$, define $\psi\in F$ by $\psi(x)=\p(x)$ is $x\ne X$ and $\psi(x)=X$ is $x=X$. Then $M\times \Omega\ni (x,\o)\mapsto\psi(x,\o)\in M$ is measurable. 
Suppose also that the law of $X$ has no atoms, then (reminding the definition of $\cF$) $\psi$ and $X$ are independent and the law of $\psi$ is $\Q$. Note that $\psi(X)=X$ and (except for very special cases) we won't have that a.s. $\cJ(\psi)(X)=\psi(X)=X$.
\brem {\em Let $\p$ and $X$ be independent random variables
respectively $F$-valued and $M$-valued.
If $\Q$, the law of $\p$, is regular and if $\cJ$ and $\cJ'$ are two measurable presentations of $\Q$.
Then $\P$-a.s., $\cJ(\p)(X)=\cJ'(\p)(X)$.
Lemma \ref{lemps} also shows that if $\Q=\Q_t$ ($t\le T$), $\p(z_{n_k^X})$ converges a.s. towards $\cJ(\p)(X)$.}\erem

\blem\label{convolution} For all $t_1,\dots,t_n$ in $[0,T]$,
\beq\label{eqconv}
(\Phi_n^T)^*(\P^{(\infty)}_{t_1}\otimes\cdots\otimes~
\P^{(\infty)}_{t_n})=\Q_{t_1+\cdots+t_n}.\eeq \elem
\prf Let us prove that
$(\Phi_n^T)^*(\P^{(\infty)}_{t_1}\otimes\cdots\otimes\P^{(\infty)}_{t_n})$
satisfies (\ref{Qn}) for all $f\in C(M^k)$, all
$x\in M^k$ and $t=t_1+\cdots+t_n$. To simplify we prove this
for $k=1$. Let $f\in C(M)$ and $x\in M$, then
applying Fubini's theorem,
\beqarr
\int_F f(\p(x))~(\Phi_n^T)^*(\P^{(\infty)}_{t_1}
\otimes\cdots\otimes  ~\P^{(\infty)}_{t_n})(d\p)\hskip-200pt&&\\
&=& \int f(i_T(y^n)\circ i_T(y^{n-1})\circ\cdots\circ
i_T(y^1)(x))~\P^{(\infty)}_{t_1}(dy^1)\otimes\cdots\otimes
\P^{(\infty)}_{t_n}(dy^n)\\ 
&=& \int \P^{(1)}_{t_n}f(i(y^{n-1})\circ\cdots\circ
i_T(y^1)(x))~\P^{(\infty)}_{t_1}(dy^1)\otimes\cdots\otimes
\P^{(\infty)}_{t_{n-1}}(dy^{n-1})\\ 
&=&\cdots~=~\P^{(1)}_{t_1+\cdots+t_n}f(x).
\eeqarr
The proof is similar for $f\in C(M)^k$ and 
$x\in M^k$. We conclude using proposition \ref{unique}. \qed 

\bprop $(\Q_t)_{t\geq 0}$ is a Feller convolution semigroup on
$(F,\cF)$.\eprop
\prf For all nonnegative $s$ and $t$, $\Phi_2^T\circ j^{\otimes 2}$ is
measurable. Proposition \ref{unique} and lemma \ref{convolution}
implies that
$(\Phi_2^T\circ j^{\otimes 2})^*(\Q_s\otimes \Q_t)=\Q_{s+t}$. Since
$(\Phi_2^T\circ j^{\otimes 2})(\p_1,\p_2)=(i_T\circ j)(\p_1)\circ
(i_T\circ j)(\p_2)$, we have easily that $\Q_s * \Q_t=\Q_{s+t}$. 
The Feller property for $\Q$ is easy to prove. \qed

\medskip This proves the first part of theorem \ref{mainthm}. \qed

\subsection{Proof of the second part of theorem \ref{mainthm}.}

We now assume we are given a Feller convolution semigroup
$\Q=(\Q_t)_{t\geq 0}$. With $\Q$, we associate a compatible family of
Feller semigroups $(\P^{(n)}_t,~n\geq 1)$ and construct
$\P^{(\infty)}_t$ like in section \ref{defprob}.

\subsubsection{Construction of a probability space.}\label{constprob}
For all $n\in\NN$, let $D_n=\{j2^{-n},~j\in\ZZ\}$ and
$D=\cup_{n\in\NN}D_n$ the set of the dyadic numbers. We take $T=1$ and
set $i=i_1$ and $\Phi_n=\Phi_n^1$.

\medskip
For all integer $n\geq 1$, let $(S_n,{\cal S}_n,\P_n)$ denote the
probability space\\ $(M^\NN,\cB(M)^{\otimes \NN}, 
\P_{2^{-n}}^{(\infty)})^{\otimes \ZZ}$.  Let
$\pi_{n-1,n}:S_n\to S_{n-1},\quad \o^n\mapsto \o^{n-1}$,  where 
\beq \o^{n-1}_{i/(2^{n-1})} =
j\circ\Phi_2(\o^{n}_{2i/2^{n}},\o^{n}_{(2i+1)/2^{n}}) =
j(i(\o^{n}_{(2i+1)/2^{n}})\circ
i(\o^{n}_{2i/2^{n}})).\eeq         
From lemma \ref{convolution}, $\pi_{n-1,n}^*(\P_n)=\P_{n-1}$.

Let $\Omega=\{(\o^n)_{n\in\NN}\in\prod S_n,~
\pi_{n-1,n}(\o^n)=\o^{n-1}\}$ and ${\cal A}$ be the $\sigma$-field on
$\Omega$ generated by the mappings $\pi_n:\Omega\to S_n$, with
$\pi_n((\o^k)_{k\in\NN})=\o^n$. Let $\P$ be the unique probability on
$(\Omega,{\cal A})$ such that $\pi_n^*(\P)=\P_n$ (see theorem 3.2 in
\cite{Partha}).

\medskip
For all dyadic numbers $s<t$, let ${\cal F}_{s,t}$ be the
$\sigma$-field generated by the mappings $(\o^k)_{k\in\NN}\mapsto \o^n_u$ for
all $(n,u)$ such that $(s,t)\in D_n^2$  and $u\in D_n\cap [s,t[$.

\subsubsection{A measurable stochastic flow of mappings on $M$.}
For $t\ge 0$, set $\cJ_t=i_t\circ j$. Then, $\cJ_t$ is a measurable presentation of $\Q_s$ for all $s\le t$. Recall that $i_t$ can be chosen such that $i_t=i$ for $t\le 1$, so that $\cJ_t=\cJ:=i\circ j$ if $t\le 1$. Note also that $\cJ_t\circ\cJ_s=i_t\circ j \circ i_s\circ j=\cJ_t$ since $j\circ i_s(y)=y$.
\begin{definition}\label{def} On $(\Omega,\cA,\P)$, we define the
following random variables 
\ben
\item For all $s<t\in D_n$, let
$\p_{s,t}^n((\o^k)_{k\in\NN})=\Phi_{(t-s)2^n}(\o^n_s,\dots,\o^n_{t-2^{-n}})$,
\item 
For all $s<t\in D$, let $\p_{s,t}=\cJ_{t-s}(\p^n_{s,t})$ where
$n=\inf\{k,~(s,t)\in D_k^2\}$.
\item For all $t\in D$, let $\p_{t,t}\in F$ be defined by $\p_{t,t}(x)=x$ for all $x\in M$.
\een
\end{definition}

Let us remark that for all $s<t\in D_n$, the law  of $\p^n_{s,t}$ is $\Q_{t-s}$ (this is a consequence of
lemma \ref{convolution}), and therefore the law of $\p_{s,t}$ is also $\Q_{t-s}$. 
Note also that for all $s\le u\le t\in D_n$,
we have  $\p^n_{s,t}=\p^n_{u,t}\circ \p^n_{s,u}$.

\brem{\em 
\begin{itemize}\item If $t\in D$, since $\p_{t,t}$ is continuous, $\cJ_0(\p_{t,t})=\p_{t,t}$. And if $s<t\in D$, $\cJ_{t-s}(\p_{s,t})=\p_{s,t}$.
\item If $u\in D_n$, then $\p_{u,u+2^{-n}}=\p^n_{u,u+2^{-n}}$. Indeed, using that $j\circ i(y)=y$ for all $y \in M^\NN$,
$\p_{u,u+2^{-n}}(\o)=i\circ j\circ i(\o^n_u)=i(\o^n_u)=\p^n_{u,u+2^{-n}}(\o)$.
\item If $s<t\in D_n$, then by definition, 
\begin{align*}
\p^n_{s,t}&=\p^n_{t-2^{-n},t}\circ \dots \circ \p^n_{s,s+2^{-n}}\\
&=\cJ(\p^n_{t-2^{-n},t})\circ \dots \circ \cJ(\p^n_{s,s+2^{-n}})\\
&=\cJ(\p_{t-2^{-n},t})\circ \dots \circ \cJ(\p_{s,s+2^{-n}})
\end{align*}
and $\p^n_{s,t}$ is a measurable function of $(\p^n_{u,u+2^{-n}})_{u\in D_n}$.
Hence, for all $s<t\in D$, $\p_{s,t}$ is a measurable function of $(\p^n_{u,u+2^{-n}})_{(n,u)\in \cup_{k\in \NN} \{k\}\times D_k}$. 
\item If $u\in D_n$, 
\begin{align*}
\p_{u,u+2^{-n}}=\p^n_{u,u+2^{-n}} &= i(\o^n_u) = i\circ j\big( i(\o^{n+1}_{u+2^{-(n+1)}})\circ i(\o^{n+1}_u)\big)\\
&= \cJ\big(\p^{n+1}_{u,u+2^{-(n+1)}}\circ \p^{n+1}_{u+2^{-(n+1)},u+2^{-n}}\big)\\
&= \cJ\big(\p_{u,u+2^{-(n+1)}}\circ \p_{u+2^{-(n+1)},u+2^{-n}}\big).
\end{align*}
\end{itemize}}\erem
\bprop \label{prop:16} For all $s<t\in D_n$ and all $M$-valued random
variable $X$ independent of ${\cal F}_{s,t}$,
$$\p^n_{s,t}(X)=\p_{s,t}(X)\qquad\hbox{$\P$-almost surely.}$$
\eprop
\prf To prove this proposition we will apply the following lemma:
\blem Let $X$ be an $M$-valued random variable and let $\p_1$ and $\p_2$ be two $F$-valued random variables. Suppose that $X$, $\p_1$ and $\p_2$ are independent and that the laws of $\p_1$ and $\p_2$ are respectively $\Q_s$ and $\Q_t$, with $(s,t)\in\RR_+^2$. Set $\p=\cJ_{s+t}\big(\cJ_t(\p_2)\circ \cJ_s(\p_1)\big)$. Then the law of $\p$ is $\Q_{s+t}$ and a.s., 
$$\p(X)= \cJ_t(\p_2)\circ \cJ_s(\p_1)(X).$$
\elem
\prf First $(x,f,g)\mapsto \cJ_t(g)\circ \cJ_s(f)(x)$ and $(x,f,g)\mapsto \cJ_{s+t}\big(\cJ_t(g)\circ \cJ_s(f)\big)(x)$ are measurable. It holds that the law of $\cJ_t(\p_2)\circ \cJ_s(\p_1)$ is $\Q_s*\Q_t=\Q_{s+t}$ and that for all $x\in M$, a.s. $\p(x)=\cJ_t(\p_2)\circ \cJ_s(\p_1)(x)$. We conclude using Fubini's theorem.    \qed \\ \\
\textbf{Proof of Proposition \ref{prop:16}.} Fubini's theorem implies that for all  $s<t\in D_n$, a.s. $\p^n_{s,t}(X)=\cJ_{t-s}(\p^n_{s,t})(X)$ (using the fact that $\p^n_{s,t}=\Phi_{(t-s)2^n}(\o^n_s,\dots,\o^n_{t-2^{-n}})$ and that $\p^n_{s,t}(x,\o)=\Psi_{(t-s)2^n}(x,\o^n_s,\dots,\o^n_{t-2^{-n}})$). Choosing $n=\inf\{k:\,(s,t)\in D_k\}$, we have that $\p^n_{s,t}(X)=\p_{s,t}(X)$ a.s.

Therefore, it is enough to prove that for all $n\ge 1$ such that $(s,t)\in D_n^2$,
$\p^n_{s,t}(X)=\p^{n+1}_{s,t}(X)$ almost surely. This holds
since 
\beqarr \p^n_{s,t}(X) &=& \cJ(\p^n_{t-2^{-n},t})\circ \dots \circ \cJ(\p^n_{s,s+2^{-n}})(X)\\
&=& \cJ \big(\p^{n+1}_{t-2^{-(n+1)},t}\circ \p^{n+1}_{t-2^{-n},t-2^{-(n+1)}}\big)
\circ\cdots\circ \cJ \big(\p^{n+1}_{s,s+2^{-(n+1)}}\circ \p^{n+1}_{s+2^{-(n+1)},s+2^{-n}}\big)(X).
\eeqarr
Since for $u\in D_{n+1}$, $\p^{n+1}_{u,u+2^{-(n+1)}}=\cJ(\p^{n+1}_{u,u+2^{-(n+1)}})$,
using the independence of the family of random variables
$\{\p^{n+1}_{u,u+2^{-(n+1)}}:\,~ u\in D_{n+1}\cap [s,t[\}$, that $X$ is independent of this family and the lemma $(t-s)2^n$ times, we prove that the last term is a.s. equal to 
$\cJ (\p^{n+1}_{t-2^{-(n+1)},t})\circ \cJ(\p^{n+1}_{t-2^{-n},t-2^{-(n+1)}})
\circ\cdots\circ \cJ (\p^{n+1}_{s,s+2^{-(n+1)}})\circ \cJ(\p^{n+1}_{s+2^{-(n+1)},s+2^{-n}})(X) = \p^{n+1}_{s,t}(X)$. \qed

\brem  {\em The preceding proposition implies that for all $s<u<t\in D$ and
all $M$-valued random variable $X$ independent of 
${\cal F}_{s,t}$,  
\beq \p_{s,t}(X)=\p_{u,t}\circ \p_{s,u}(X) \qquad
\P-\hbox{almost surely}.\eeq }
\erem
\prf There is $n$ such that $(s,t,u)\in D_n^3$, then a.s. $\p_{s,t}(X)=\p^n_{s,t}(X)=\p^n_{u,t}\circ \p^n_{s,u}(X)=\p_{u,t}\circ \p_{s,u}(X)$. \qed.

\medskip
We now intend to define by approximation for all $s<t$ in
$\RR$ a $(F,\cF)$-valued random variable $\p_{s,t}$ of law
$\Q_{t-s}$. In order to do this, we prove the following lemma.

\blem \label{tcont} For all continuous function $f$ on 
$M^2$, the mapping 
\beq(s,t,u,v,x,y)\mapsto
\E[f(\p_{s,t}(x),\p_{u,v}(y))]\eeq
is continuous on $\{(s,t)\in D^2,~s\leq t\}^2\times M^2$. (And
therefore uniformly continuous on every compact.) \elem 
\prf For all $s\leq u\leq t\leq v$ in $D$, using the cocycle property, we have (In order to apply Fubini's theorem, we use that for all $s\le t$, $\p_{s,t}=\cJ_{t-s}(\p_{s,t})$)
\beqarr
\E[f(\p_{s,t}(x),\p_{u,v}(y))]
&=&\E[f(\p_{u,t}\circ \p_{s,u}(x),\p_{t,v}\circ \p_{v,t}(y))]\\
&=&(\P^{(1)}_{u-s}\otimes I)\P^{(2)}_{t-u}(I\otimes\P^{(1)}_{v-t})f(x,y).
\eeqarr
For all $s\leq u\leq v\leq t$ in $D$, using the cocycle property, we have
\beqarr
\E[f(\p_{s,t}(x),\p_{u,v}(y))]
&=&\E[f(\p_{v,t}\circ \p_{u,v}\circ \p_{s,u}(x),\p_{u,v}(y))]\\
&=&(\P^{(1)}_{u-s}\otimes I)\P^{(2)}_{v-u}(\P^{(1)}_{t-v}\otimes I)f(x,y).
\eeqarr
For all $s\leq t\leq u\leq v$ in $D$, 
$$\E[f(\p_{s,t}(x),\p_{u,v}(y))]
=(\P^{(1)}_{t-s}\otimes\P^{(1)}_{v-u})f(x,y).$$
All these functions are continuous and they join. This implies the
lemma. \qed

\medskip
For all real $t$ and all integer $n$, let $t_n=\sup\{u\in D_n,~u\leq
t\}$. For all $s<t\in\RR$, we define the increasing sequences
$(s_n)_{n\in\NN}$ and $(t_n)_{n\in\NN}$. Using lemma \ref{tcont} for
$f(x,y)=d(x,y)$ and the Markov inequality, for all positive $\eps$, we
have
\beq\label{eq}\lim_{n\to\infty}\sup_{k>n}\sup_{x\in M}
\P[d(\p_{s_n,t_n}(x),\p_{s_k,t_k}(x))\geq\eps]=0. \eeq

Set $\cJ^n=\cJ_{t_n-s_n}$ and define the following measurable mappings
$\Phi: F^\NN\to F$ and  $\Psi: M\times F^\NN \to M$ by $\Psi(x,(\p_n)_{n\in \NN})=l\big((\cJ^n(\p_n)(x))_{n\in \NN}\big)$ and $ \Phi((\p_n)_{n\in \NN})(\cdot)=\Psi(\cdot,(\p_n)_{n\in \NN})$. 
We then set for all $s\le t$, $\psi_{s,t}= \Phi((\p_{s_n,t_n})_{n\in \NN})$ and $\p_{s,t}=\cJ_{t-s}(\psi_{s,t})$. 
Note that 
\begin{itemize}
\item for all $x\in M$, a.s. $\psi_{s,t}(x)=l\left((\p_{s_n,t_n}(x))\right)$ and when $s=t$, a.s.  $\psi_{t,t}(x)=x$. 
\item  $\p_{s,t}$ is a measurable function of $(\p^n_{u,u+2^{-n}})_{(n,u)\in \cup_{k\in \NN} \{k\}\times D_k}$.
\end{itemize}
 \blem\label{convprob} For all positive $\eps$ and all $s\leq t$,
\beq\lim_{n\to\infty}\sup_{x\in M}\P[d(\p_{s_n,t_n}(x),\psi_{s,t}(x))\geq\eps]=0.\eeq \elem
\prf Equation (\ref{eq}) implies that $\p_{s_n,t_n}(x)$ converges in
probability towards $\psi_{s,t}(x)$. Thus, for all positive $\eps$,
$$\P[d(\p_{s_n,t_n}(x),\psi_{s,t}(x))\geq\eps] =
\lim_{k\to\infty}\P[d(\p_{s_n,t_n}(x),\p_{s_k,t_k}(x))\geq \eps].$$
Therefore, 
$$\sup_{x\in M} \P[d(\p_{s_n,t_n}(x),\psi_{s,t}(x))\geq\eps]
\leq \sup_{k>n}\sup_{x\in M}
\P[d(\p_{s_n,t_n}(x),\p_{s_k,t_k}(x))\geq\eps],$$
which implies the lemma. \qed 

\bprop For all $s<t\in\RR$, the law of $\p_{s,t}$ is $\Q_{t-s}$. \eprop
\prf For all $k\geq 1$, $f\in C(M^k)$ and $x\in M^k$, lemma
\ref{convprob} and dominated convergence theorem implies that
\beqarr
\E[f\circ\psi_{s,t}^{\otimes k}(x))]
&=&\lim_{n\to\infty}\E[f\circ\p_{s_n,t_n}^{\otimes k}(x)]\\
&=&\lim_{n\to\infty}\P^{(k)}_{t_n-s_n}f(x) 
\quad=\quad \P^{(k)}_{t-s}f(x)
\eeqarr
since $\P^{(k)}_t$ is Feller. This implies that the law of $\psi_{s,t}$
is $\Q_{t-s}$, and therefore that the law of $\p_{s,t}=\cJ_{t-s}(\psi_{s,t})$ is $\Q_{t-s}$. \qed

\medskip
Let us now prove the cocycle property.
\bprop\label{cocycle}  For all $s<u<t$ and all $x\in M$, $\P$-almost
surely, \beq\label{eqcocycle}\p_{s,t}(x)=\p_{u,t}\circ
\p_{s,u}(x).\eeq 
\eprop
\prf
By construction, for all $s\le t$, $\cJ_{t-s}(\p_{s,t})=\p_{s,t}$ (since $\cJ_t\circ \cJ_t=\cJ_t$).

Almost surely, we have
$\p_{s_n,t_n}(x) = \p_{u_n,t_n}\circ
\p_{s_n,u_n}(x)$ since $s_n<u_n<t_n$ belong to $D$. Proposition 1.7 implies that $\p_{s,t}(x)=\psi_{s,t}(x)$ a.s. and Lemma 1.10 implies that $\p_{s_n,t_n}(x)$ converges in probability towards
$\p_{s,t}(x)$.

Let us now show that $\p_{u_n,t_n}\circ
\p_{s_n,u_n}(x)$ converges in probability toward $\p_{u,t}\circ
\p_{s,u}(x)$.
For $n\ge 1$, set $X_n=\p_{s_n,u_n}(x)$ and $X=\p_{s,u}(x)$. Then $X_n$ converges in probability towards $X$ and for $n\ge 1$ and $\varepsilon>0$, 
\begin{align*}
\P[d(\p_{u_n,t_n}&\circ
\p_{s_n,u_n}(x),\p_{u,t} 
\circ \p_{s,u}(x))\geq \eps]\\
= & \qquad \P[d(\cJ_{t_n-u_n}(\p_{u_n,t_n})(X_n),\cJ_{t-u}(\p_{u,t})(X)\geq \eps]\\
\leq & \qquad \P[d(\cJ_{t_n-u_n}(\p_{u_n,t_n})(X_n),\cJ_{t-u}(\p_{u,t})(X_n))\geq \eps/2]\\
& + \quad \P[d(\cJ_{t-u}(\p_{u,t})(X_n),\cJ_{t-u}(\p_{u,t})(X))\geq \eps/2].
\end{align*}
Since $X_n$ is independent of $\cF_{u_n,t}$, the first term is equal to 
\begin{align*}
\int \P[d(\cJ_{t_n-u_n}&(\p_{u_n,t_n})(y),\cJ_{t-u}(\p_{u,t})(y))\geq \eps/2]\, \P_{X_n}(dy)\\ 
&\le\quad \sup_{y\in M} \P[d(\p_{u_n,t_n}(y),\p_{u,t}(y))\geq \eps/2]
\end{align*}
which converges to $0$ as $n\to\infty$ (using Lemma 1.10). We show that the second term converges to $0$ using Lemma 1.6 (since the sequence $(X_n)$ is independent of $\p_{u,t}$). \qed

\medskip
Thus we have constructed a measurable stochastic flow of measurable mappings
on $M$ associated with the compatible family of Feller
semigroups $(\P^{(k)}_t,~k\geq 1)$ and with the Feller convolution
semigroup $(\Q_t,~t\geq 0)$.

\medskip
Let $\p$ be the $(\O^0,\cA^0)$-valued random variable defined
by $\p=(\p_{s,t},~s\leq t)$. Let $\P_\Q = \p^*(\P)$ be the law of
$\p$. Then by a monotone class argument we show that $T_h^*(\P_\Q) =
\P_\Q$ for all $h\in \RR$.

On $(\Omega^0,\mathcal{A}^0,\P_\Q)$, let $\p'$ be a measurable modification of the canonical stochastic flow $\p^0$. Then for all $t\ge 0$, there is a measurable presentation $\cJ'_t$ of $\Q_t$ such that for all $s\le t$, $\p'_{s,t}=\cJ'_{t-s}(\p^0_{s,t})$.
This modification is $(T_h)_{h\in \mathbb{R}}$-invariant, since for all $s\le t$ and all $h\in \mathbb{R}$, $\p'_{s+h,t+h}=\cJ_{t-s}(\p^0_{s+h,t+h})= \cJ_{t-s}(\p^0_{s,t}\circ T_h)=\p'_{s,t}\circ T_h$.  

 The
fact that $(\O^0,\cA^0,\P_\Q)$ is separable is a consequence of the
construction of $\p$. The proof of Theorem \ref{mainthm} is
finished. \qed

\subsection{The example of Lipschitz SDEs.}\label{Lipschitz}

We first show a sufficient condition for a compatible family of
Markovian kernels semigroups to be constituted of Feller semigroups.

\blem\label{Fel2n}
A compatible family $(\P^{(n)}_t,~n\geq 1)$ of semigroups of Markovian
kernels is constituted of Feller semigroups when the following
condition is satisfied  
\bdes \item[(F)] For all $f\in C(M)$ and $x\in M$,
$\lim_{t\to 0}\P^{(1)}_tf(x)=f(x)$ and for all $x\in M$, $\eps>0$ and
$t>0$, $\lim_{y\to x}\P^{(2)}_tf_\eps(x,y)=0$, where
$f_\eps(x,y)=1_{d(x,y)>\eps}$.\edes \elem 
\prf Let $h\in C(M^n)$ be in the form $f_1\otimes\cdots\otimes f_n$ and
$x=(x_1,\dots,x_n)$ in $M^n$. We have for $M$ large enough 
\beq |\P^{(n)}_t h(x) - h(x)| \leq M\sum_{k=1}^n(\P^{(1)}_tf_k^2 +
f_k^2 -2f_k\P^{(1)}_tf_k)^{\frac{1}{2}}(x_k) \eeq
which converges towards 0 as $t$ goes to 0 since for all $f\in C(M)$
and all $x\in M$, $\lim_{t\to 0}\P^{(1)}_tf(x)=f(x)$. We also have for
$y=(y_1,\dots,y_n)$ in $M^n$,
\beq |\P^{(n)}_t h(y) - \P^{(n)}_t h(x)| \leq
M\sum_{k=1}^n \P^{(2)}_t(|1\otimes f_k-f_k\otimes 1|)(y_k,x_k) \eeq
which converges towards 0 as $y$ tends to $x$ since for all $f\in
C(M)$ and $x\in M$, $\lim_{y\to x} \P^{(2)}_t(|1\otimes f-f\otimes 1|)(y,x)=0$. 
Indeed, $\forall \alpha>0$, $\exists \eps>0$ such that $d(x,y)<\eps$
implies $|f(y)-f(x)|<\alpha$. This implies
\beq \P^{(2)}_t(|1\otimes f-f\otimes 1|)(y,x) \leq \alpha +
2\|f\|_\infty \P^{(2)}_t f_\eps (x,y).\eeq
This implies $\limsup_{y\to x} \P^{(2)}_t(|1\otimes f-f\otimes
1|)(y,x)\leq \alpha$ for all $\alpha>0$. \qed

\brem  $\bullet$ The previous result extends to the locally compact
case (using the fact that $C_0(M)$ is constituted of uniformly
continuous functions). 
 
$\bullet$ When {\bf (F)} is satisfied, for all positive $t$, $f\in
C_0(M)$ and $x\in M$, $\P^{(2)}_tf^{\otimes
2}(x,x)=\P^{(1)}_tf^2(x)$. This implies that {\bf (F)} is not a
necessary condition. Theorem \ref{mainthm} shows that a stochastic
flow of mappings is associated with this family of semigroups. \erem

\begin{definition}
A two-parametric family $(W_{s,t},~s\leq t)$ of real random variables
is called a real white noise if 
\bdes \iti for all $s<t$, $W_{s,t}$ is a centered Gaussian
variable with  variance $t-s$,
\itii for all $((s_i,t_i),~1\leq i\leq n)$ with $s_i\leq t_i\leq s_{i+1}$, the
random variables $(W_{s_i,t_i},~1\leq i\leq n)$ are
independent and 
\itiii for all $s\leq t\leq u$, $W_{s,u}=W_{s,t}+W_{t,u}$.
\edes
\end{definition}

Let $V,V_1,\dots,V_k$ be bounded Lipschitz vector fields on a smooth
locally compact manifold $M$. We also assume that $V_1,\dots,V_k$ are
$C^1$. Let $W^1,\dots,W^k$ be $k$ independent
real white noises. We consider the SDE on $M$ 
\beq\label{LipSDE}dX_t=\sum_{i=1}^k V_i(X_t)\circ dW^i_t +
V(X_t)~dt,\qquad t\in\RR.\eeq
From the usual theory of strong solutions of SDEs (see for example
\cite{ku}), it is possible to construct a stochastic flow of
diffeomorphisms $(\p_{s,t},~s\leq t)$ such that for all $x\in M$,
$\p_{s,t}(x)$ is a strong solution of the SDE (\ref{LipSDE}) with
$\p_{s,s}(x)=x$. 

\smallskip Using this stochastic flow, it is possible to construct a
compatible family of Markovian semigroups $(\P^{(n)}_t,~n\geq 1)$ with
\beq \P^{(n)}_th(x_1,\dots,x_n)=\E[h(\p_{0,t}(x_1),\dots,\p_{0,t}(x_n))]\eeq
for $h\in C(M^n)$ and $x_1,\dots,x_n$ in $M$. Using lemma \ref{Fel2n},
it is easy to check that these semigroups are Feller (these properties
were previously observed by P. Baxendale in \cite{baxendale}).

It can easily be shown that the canonical stochastic flow of
maps associated with this family of semigroups is equal in law to
$(\p_{s,t},~s\leq t)$.

\setcounter{equation}{0}
\section{Stochastic flow of kernels.}\label{sfk}
\subsection{Notation and definitions.} 
In order to simplify,we suppose  in this section that $M$ is a compact metric space. But, as it is explained in section \ref{sec:conschar},  all the results of this section extend to locally compact separable metric spaces.
We denote by ${\cal P}(M)$ the space of probability measures on $M$,
equipped with the weak convergence topology. Let $(f_n)_{n\in\NN}$ be
a sequence of functions dense in 
$\{f\in C(M),~\|f\|_\infty\leq 1\}$. We will equip ${\cal P}(M)$ with
the distance $\rho(\mu,\nu)=(\sum_{n}2^{-n}(\int f_n~d\mu-\int f_n~d\nu)^2)^{1/2}$
for all $\mu$ and $\nu$ in ${\cal P}(M)$. Thus ${\cal P}(M)$ is a
separable compact metric space.

Let us recall that a kernel $K$ on $M$ is a measurable
mapping from $M$ into ${\cal P}(M)$, $M$ and ${\cal P}(M)$ being
equipped with their Borel $\sigma$-fields. For all $f\in C(M)$
and $x\in M$, $Kf(x)$ denotes $\int f(y)~K(x,dy)$. For all
$\mu\in\cP(M)$, $\mu K$ denotes the probability measure defined
by $\int f(y)~\mu K(dy) = \int Kf(x)~\mu(dx)$.
We denote by $E$ the space of all kernels on $M$ and we equip $E$
with the $\sigma$-field generated by the mappings $K\mapsto \mu K$,
for all $\mu\in \cP(M)$ ($\cP(M)$ is equipped with its Borel
$\sigma$-field $\cB(\cP(M))$). We denote this $\sigma$-field by $\cE$.

Let $\Gamma$ denote the space of measurable mappings on $\cP(M)$. We
equip $\Gamma$ with the $\sigma$-field generated by the mappings
$\Phi\mapsto\Phi(\mu)$ for all $\mu\in\cP(M)$. Note that
$(\Gamma,\cG)=(F,\cF)$ once we have replaced $M$ by $\cP(M)$.

\subsection{Convolution semigroups on the space of kernels.}
Let $\cI$ denote the measurable mapping from $(E,\cE)$ on
$(\Gamma,\cG)$ defined by $\cI(K)(\mu)=\mu K$. Note that $\cI(E)$ is
not measurable in $\Gamma$ but $\cI$ is measurable.

\begin{definition} \label{defconvE}
\bdes \iti A probability measure $\nu$ on $(E,\cE)$
is called regular if $\cI^*(\nu)$ is a regular probability measure on
$(\Gamma,\cG)$.
\itii A convolution semigroup on $(E,\cE)$ is a family
$(\nu_t)_{t\geq 0}$ of regular probability measures on $(E,\cE)$ such
that $(\cI^*(\nu_t))_{t\geq 0}$ is a convolution semigroup on
$(\Gamma,\cG)$. \edes \end{definition}

Let $\delta~:~\Gamma\to E$ be the mapping defined by
$\delta(\Phi)(x)=\Phi(\delta_x)$. Note that $\delta$ is not measurable
in general.

\bprop \label{Jnu} Let $\Q$ be a regular probability measure on
$(\Gamma,\cG)$ and $\cJ$ a measurable presentation of $\Q$. Then
$\delta\circ \cJ$ is measurable and the probability measure
$\nu=(\delta\circ \cJ)^*(\Q)$ is a regular probability measure on
$(E,\cE)$ if $\cI^*(\nu)=\Q$. \eprop
\prf Let $\Q$ be a regular probability measure on $(\Gamma,\cG)$ and
$\cJ$ a measurable presentation of $\Q$.
The mapping $\cP(M)\times \Gamma\ni (\mu,\Phi)\mapsto
\cJ(\Phi)(\mu)\in \cP(M)$ and $M\ni x\mapsto \delta_x\in\cP(M)$ are
measurable. Thus $M\times \Gamma\ni (x,\Phi)\mapsto
\delta\circ \cJ(\Phi)(x)=\cJ(\Phi)(\delta_x)\in \cP(M)$ is measurable,
which implies that $\delta\circ \cJ$ is measurable. \qed

\brem  {\em The probability measure $\nu$ defined in proposition
  \ref{Jnu} depends only of $\Q$. Indeed, if $\cJ'$ is
  another measurable presentation  of $\Q$, for all $x\in M$,
  $\Q(d\Phi)$-a.s., $\delta\circ
  \cJ(\Phi)(x)=\delta\circ\cJ'(\Phi)(x)$, which implies by Fubini
  theorem that for all $\mu\in\cP(M)$, $\Q(d\Phi)$-a.s.,
  $\mu(\delta\circ\cJ(\Phi)) = \mu(\delta\circ\cJ'(\Phi))$ and then
  that $(\delta\circ\cJ)^*(\Q) = (\delta\circ\cJ')^*(\Q)$. }\erem

\begin{definition} A measurable presentation of a regular probability measure $\nu$ is a measurable mapping $\mathfrak{p}:(E,\cE)\to (E,\cE)$ such that $(x,K)\mapsto\mathfrak{p}(K)(x)$ is measurable and such that for all $x\in M$, $\nu(dK)$-a.s. $\mathfrak{p}(K)(x)=K(x)$.
\end{definition}
\brem{\em
\begin{itemize}
\item If $\mathfrak{p}$ is a measurable presentation of a regular probability measure $\nu$, then $(\mu,K)\mapsto \mu (\mathfrak{p}(K))$ is measurable and for all $\mu\in \cP(M)$, $\nu(dK)$-a.s. $\mu (\mathfrak{p}(K))=\mu K$. 
\item If $\nu$ is a regular probability measure on $(E,\cE)$ and if $\cJ$ is a measurable presentation of $\Q:=\cI^*(\nu)$, then for all $x\in M$, if $\nu(dK)$-a.s. $\delta\circ \cJ \circ \cI(K)(x)=K(x)$. Therefore, the mapping $\mathfrak{p}=\delta\circ \cJ \circ \cI$ is a measurable presentation of $\nu$.
\item Let $(\nu_t)_{t\ge 0}$ be a convolution semigroup on $(E,\cE)$. If $K_1$ and $K_2$ are random kernels with laws $\nu_s$ and $\nu_t$ and if $\mathfrak{p}_t$ is a measurable presentation of $\nu_t$, then $K_1(\mathfrak{p}_t(K_2))$ is a random kernel with law $\nu_{s+t}$ (note that $(K_1,K_2)\mapsto K_1(\mathfrak{p}_t(K_2))$ is measurable). 
\end{itemize}  
}\erem   

\begin{definition} \label{defconvF}
A convolution semigroup $(\nu_t)_{t\geq 0}$ on $(E,\cE)$ is called
Feller if 
\bdes \iti for all $f\in C(M)$, $\lim_{t\to 0} \sup_{x\in
M}\int(Kf(x)-f(x))^2\nu_t(dK) = 0$,  
\itii for all $f\in C(M)$ and all $t\geq 0$, 
$\lim_{d(x,y)\to 0} \int(Kf(x)-Kf(y))^2\nu_t(dK) = 0$. 
\edes  \end{definition}

\bprop \label{nuP}Let $(\nu_t)_{t\geq 0}$  be a Feller convolution
semigroup on $(E,\cE)$. For all $n\geq 1$, $f\in C(M)$ and $x\in M^n$, set
\beq \P^{(n)}_tf(x)=\int K^{\otimes n} f(x)~\nu_t(dK).\eeq
Then $(\P^{(n)}_t,~n\geq 1)$ is a compatible family of Feller
semigroups on $M$. \eprop
\prf This is the same proof as the one of proposition \ref{QP}. \qed

\bprop \label{qnu}
Let $(\Q_t)_{t\geq 0}$ be a convolution semigroup on
 $(\Gamma,\cG)$. Let $\cJ_t$ be a measurable presentation of $\Q_t$
 and $\nu_t=(\delta\circ\cJ_{t})^*(\Q_t)$. If $\Q_t=\cI^*(\nu_t)$,
 $(\nu_t)_{t\geq 0}$ is a
 convolution semigroup on $(E,\cE)$. Then, $(\Q_t)_{t\geq 0}$
 is Feller if and only if $(\nu_t)_{t\geq 0}$ is Feller.
 \eprop
\prf The fact that $(\nu_t)_{t\geq 0}$ is a convolution semigroup
 follows from the definition \ref{defconvE}.

Note that $(\Q_t)_{t\geq 0}$ is Feller if and only if for all $f\in C(M)$,
\beqar \lim_{t\to 0} 
\sup_{\mu\in\cP(M)} \int (\Phi(\mu)f-\mu f)^2\Q_t(d\Phi) &=& 0\label{Q1}\\
\lim_{\rho(\mu,\nu)\to 0} \int (\Phi(\mu)f-\Phi(\nu) f))^2\Q_t(d\Phi)
&=& 0. \label{Q2}\eeqar

We first prove (\ref{Q1}) and {\bf (i)} in definition \ref{defconvF} are
equivalent. Equation (\ref{Q1}) implies {\bf (i)} since 
$\int (Kf(x)-f(x))^2\nu_t(dK) = 
\int(\Phi(\delta_x)f-\delta_xf)^2\Q_t(d\Phi)$.
And {\bf (i)} implies (\ref{Q1}) since
\beqarr \int (\Phi(\mu)f-\mu f)^2\Q_t(d\Phi)
&=& \int (\mu Kf-\mu f)^2\nu_t(dK)\\
&\leq& \int\left(\int (Kf(x)-f(x))^2\nu_t(dK)\right)\mu(dx).
\eeqarr

We now prove (\ref{Q2}) and {\bf (ii)} in definition \ref{defconvF} are
equivalent. Equation (\ref{Q2}) implies {\bf (ii)} since 
$\int (Kf(x)-Kf(y))^2\nu_t(dK)  = 
\int (\Phi(\delta_x)f-\Phi(\delta_y)f)^2\Q_t(d\Phi)$ 
and  $\lim_{d(x,y)\to 0}\rho(\delta_x,\delta_y)=0$. Assume {\bf (ii)}
holds. For $\mu$ and $\nu$ in $\cP(M)$, we have
\beqarr \int (\Phi(\mu)f-\Phi(\nu) f)^2\Q_t(d\Phi)
&=& \int (\mu Kf-\nu Kf)^2\nu_t(dK) \\
&=& (\mu-\nu)^{\otimes 2} \int K^{\otimes 2} f^{\otimes 2}~\nu_t(dK).
\eeqarr
We conclude since $\int K^{\otimes 2} f^{\otimes 2}~\nu_t(dK)$ is a
continuous function (see proposition \ref{nuP}). \qed

\subsection{Stochastic flows of kernels.}

\begin{definition}\label{defsfk} Let $(\Omega,{\cal A},\P)$ be a
probability space and let $K=(K_{s,t},~s\leq t)$ be a family of $(E,\cE)$-valued random variables such that for all $x\in M$ and all $t\in\mathbb{R}$, $\P$-a.s. $K_{t,t}(x)=\delta_x$. For $t\ge 0$, denote by $\nu_t$ the law of $K_{0,t}$.
The family $K$ is called a stochastic flow of kernels if for all $t\ge 0$, $\nu_t$ is regular and if the following properties are satisfied by $K$ 
\bdes \ita For all $s\le u\le t$, for all $x\in M$, for all $f\in C(M)$ and for all measurable presentation $\mathfrak{p}_{t-u}$ of $\nu_{t-u}$, $\P$-almost surely,
$K_{s,t}f(x)=K_{s,u}(\mathfrak{p}_{t-u}(K_{u,t}) f)(x)$. (cocycle  property)
\itb For all $s\leq t$, the law of  $K_{s,t}$ is $\nu_{t-s}$. (Stationarity)
\itc The flow has independent increments, i.e. for all
$t_1\le t_2\le \cdots\le t_n$, the family $\{K_{t_i,t_{i+1}},~1\leq i\leq n-1\}$
is independent.
\itd For all $f\in C(M)$ and all $s\leq t$, $\lim_{(u,v)\to (s,t)}\sup_{x\in M}
\E[(K_{s,t}f(x)-K_{u,v}f(x)) ^2]=0$.
\ite For all $f\in C(M)$ and all $s\leq t$,
$\lim_{d(x,y)\to 0}\E[(K_{s,t}f(x)-K_{s,t}f(y))^2]=0$.
\edes
\end{definition}
\brem {\em
\begin{itemize}
\item Item (a) holds for all set of measurable presentations as soon as (a) holds for one of them.
\item If $K'$ is equal in law to a stochastic flow of kernels $K$, then $K'$ is also a stochastic flow of kernels.
\end{itemize}
}\erem
\bprop \label{prop:233}
Let $K$ be a stochastic flow of kernels, and for $t\ge 0$, let $\mathfrak{p}_t$ be a measurable presentation of the law of $K_{0,t}$. Then $K'=(\mathfrak{p}_{t-s}(K_{s,t}),\, s\le t)$ is a stochastic flow of kernels satisfying
\begin{itemize}
\item[(i)] For all $s\le t$ and $\mu\in \cP(M)$, a.s. $\mu K'_{s,t}(x)=\mu K_{s,t}(x)$.
\item[(ii)] For all $s\le u\le t$ and for all $\mu\in \cP(M)$, $\P$-almost surely,
$\mu K'_{s,t}= \mu K'_{s,u} K'_{u,t}$.
\end{itemize} \eprop
\prf Follow the proof of proposition at page 10 in section 2. \qed
\begin{definition}
\begin{itemize}
\item The stochastic flow of kernels $K'$ defined in the proposition \ref{prop:233} will be called a measurable modification of $K$.
\item  A stochastic flow of kernels which is a measurable modification of a stochastic flow of kernels is called a measurable stochastic flow of kernels.
\end{itemize}
\end{definition}

\bprop \label{SFKP}
Let $(K_{s,t},~s\leq t)$  be a stochastic flow of
kernels. For all $n\geq 1$, $f\in C(M^n)$ and $x\in M^n$, set 
\beq \P^{(n)}_tf(x)=\E[K^{\otimes n}f(x)].\eeq
Then $(\P^{(n)}_t,~n\geq 1)$ is a compatible family of Feller
semigroups on $M$. \eprop
\prf This is the same proof as the one to prove proposition
\ref{QP}. \qed

\subsection{Construction and characterization.}
Let $(\Omega^0,{\cal A}^0)$ denote the measurable space
$(\prod_{s\leq t}E ,\otimes_{s\leq t} {\cal E})$.
For $s\leq t$, let $K^0_{s,t}$ denote the random variable
$\omega\mapsto \omega(s,t)$. Let also $K^0$ be the random variable
$(K^0_{s,t},~s\leq t)$. Then $K^0(\omega)=\omega$.
Let $(T_h)_{h\in\RR}$ be the one-parametric group of transformations
of $\Omega^0$ defined by $T_h(\omega)(s,t)=\omega(s+h,t+h)$, for all
$s\leq t$, $h\in \RR$ and $\omega\in\Omega^0$. 

\bthm\label{mainthmk} {\bf 1-} For all compatible family
$(\P^{(n)}_t,~n\geq 1)$  of Feller semigroups on $M$, there
exists a unique Feller convolution semigroup $(\nu_t)_{t\geq 0}$ on
$(E,{\cal E})$ such that for all $n\geq 1$, $t\geq 0$, $f\in C(M^n)$
and $x\in M^n$, 
\beq \label{eqrelconv} \P^{(n)}_tf(x)=\int K^{\otimes n}
f(x)~\nu_t(dK). \eeq

\medskip
{\bf 2-} For all Feller convolution semigroup $\nu=(\nu_t)_{t\geq 0}$
on $(E,\cE)$, there exists a unique $(T_h)_{h\in\RR}$-invariant
probability measure $\P_\nu$ on $(\Omega^0,\cA^0)$ such that
$(\Omega^0,\cA^0,\P_\nu)$ is separable, the family of random
variables $(K^0_{s,t},~s\leq t)$ is a stochastic flow of kernels and for
all $s\leq t$, the law of $K^0_{s,t}$ is $\nu_{t-s}$. 
Every measurable modification $K'$
of $K^0$ satisfies $K'_{s+h,t+h}=K'_{s,t}\circ T_h$ for all $s\le t$ and all $h\in \mathbb{R}$.

The flow $K^0$ is called the canonical stochastic flow of
kernels associated with $\nu$ (or equivalently with
$(\P^{(n)}_t,~n\geq 1)$).
 \ethm

\brem {\em
In the case (\ref{mapcond}) is satisfied, $(\P^{(n)}_t,~n\geq 1)$ is associated to  a  stochastic flow of mappings $\p$ by Theorem 1.1. 
Set $\delta_\p=(\delta_{\p_{s,t}},~s\leq t)$. 
Then, for all $s\le t$, $\delta_{\p_{s,t}}$ is a random kernel of law $\nu_{s,t}$ (since $\E[\delta_{\p_{s,t}}^{\otimes n} f(x)]=\E[f\circ \p_{s,t}^{\otimes n} (x)] = \P^{(n)}_{t-s}f(x)$) and one can check that $\delta_\p$ is a stochastic flow of kernels of law $\P_{\nu}$.
}\erem

\subsection{Proof of theorem \ref{mainthmk}.} 
Let $(\P^{(n)}_t,~n\geq 1)$ be a compatible family of Feller
semigroups on $M$. Starting with this family of semigroups, we intend
to construct a Feller convolution semigroup $\nu=(\nu_t)_{t\geq 0}$ on
$(E,\cE)$. The idea is to construct a compatible family of Feller
semigroups on $\cP(M)$, then to apply theorem \ref{mainthm} to
construct a Feller convolution semigroup  $\Q=(\Q_t)_{t\geq 0}$ on
$(\Gamma,\cG)$ and to construct $\nu$ using the mappings $\delta\circ
\cJ_t$, where $\cJ_t$ is a measurable presentation of law $\Q_t$.

\subsubsection{Construction of a compatible family of Feller
  semigroups on $\cP(M)$.} 
For all integer $k$, we define a Feller semigroup $\Pi^{(k)}_t$ acting
on the continuous functions on ${\cal P}(M)^k$ (see \cite{Ma}
for a similar construction when $k=1$).

Let $\cA_k$ denote the algebra of functions
$g:\cP(M)^k\to\RR$ such that \footnote{Here and in the following,
for all measure $\mu$ and $f\in L^1(\mu)$, we denote $\int f~d\mu$ by
$\langle f,\mu\rangle$, $\langle \mu,f\rangle$ or $\mu f$.}

\beq\label{defg} g(\mu_1,\dots,\mu_k) = 
\langle f,\mu_1^{\otimes n_1}\otimes\cdots\otimes\mu_k^{\otimes
n_k}\rangle \eeq
for $f\in C(M^n)$ and $n_1,\dots,n_k$ integers such that
$n=n_1+\cdots+n_k$ (${\cal A}_k$ is the union of an increasing family
of algebras ${\cal A}_{n_1,\dots,n_k}$). For all $g\in{\cal A}_k$,
given by equation (\ref{defg}), let
\beq\label{defpik} \Pi_t^{(k)}g(\mu)=\langle \P^{(n)}_tf, \mu_1^{\otimes
n_1}\otimes\cdots\otimes\mu_k^{\otimes n_k}\rangle.\eeq
with $\mu=(\mu_1,\dots,\mu_k)\in {\cal P}(M)^k$. Since the family of
semigroups $(\P^{(n)}_t,~n\geq 1)$ is compatible, (\ref{defpik}) is
independent of the expression of $g$ in (\ref{defg}).

Let us notice that $\Pi^{(k)}_t$ acts on ${\cal A}_k$ and that, by the
theorem of Stone-Weierstrass, the algebra ${\cal A}_k$ is dense in
$C({\cal P}(M)^k)$.
\blem $\Pi^{(k)}_t$ is a Markovian operator acting on ${\cal A}_k$. \elem 
\prf The only thing to be proved is the positivity property (it is
obvious that $\Pi^{(k)}_t1=1$).

For all integer $N$, let $(X^{j,i},~1\leq i\leq k,~1\leq j\leq N)$ be
a Markov process associated with the Markovian semigroup $\P^{(Nk)}_t$
such that the random variables $(X^{j,i}_0,~1\leq i\leq k,~1\leq j\leq
N)$ are independent and the law of $X^{j,i}_0$ is $\mu_i$, where
$(\mu_1,\dots,\mu_k)\in\mathcal{P}(M)^k$. Let us introduce
the following Markov process on ${\cal P}(M)^k$,
$\mu^{N}_t=(\mu^{N,1}_t,\dots,\mu^{N,k}_t)$ where
\beq \mu^{N,i}_t = \frac{1}{N}\sum_{j=1}^N\delta_{X^{j,i}_t}, \quad
\hbox{ for } 1\leq i\leq k. \eeq

For $g(\mu_1,\dots,\mu_k) = 
\langle f,\mu_1^{\otimes n_1}\otimes\cdots\otimes\mu_k^{n_k}\rangle$,
we have
\beqarr \E[g(\mu_t^{N})] &=&
\E[\la f,(\mu^{N,1}_t)^{\otimes n_1} \otimes \cdots \otimes
(\mu^{N,k}_t)^{\otimes n_k}\ra]\\
&=& \frac{1}{N^n} \sum_{i=1}^k\sum_{l=1}^{n_k}\sum_{j_i^l=1}^N
\E[f(X^{j_1^1,1}_t , X^{j^2_1,1}_t, \dots, X^{j^{n_1}_1,1}_t,
X^{j_2^1,2}_t,\dots, X^{j_k^{n_k},n_k}_t)]\\
&=& \langle \P^{(n)}_t f, \mu_1^{\otimes n_1} \otimes\cdots\otimes
\mu_k^{n_k}\rangle + R_N. \eeqarr
The remainder term $R_N$ comes from terms in which $j_i^a=j_i^b$ for
some $a\neq b$ and some $i$ and is therefore dominated by 
$2\|f\|_\infty(1-\prod_{i=1}^k
(N(N-1)\cdots(n-n_i+1)/N^{n_i}))$. Thus
\beqar
\lim_{N\to\infty} E[g(\mu^N_t)]
&=& \langle \P^{(n)}_t f, \mu_1^{\otimes n_1} \otimes\cdots\otimes
\mu_k^{n_k}\rangle\\
&=& \Pi^{(k)}_t g(\mu_1,\dots,\mu_k). \eeqar
This shows that $\Pi^{(k)}_t$ is positive. \qed 

\medskip
Using this lemma, it is easy to define $\Pi^{(k)}_tg$ for all
continuous function $g$ and to show that $\Pi^{(k)}_t$ is a Markovian
semigroup acting on $C(M^n)$.

\blem $(\Pi^{(k)}_t,~k\geq 1)$ is a compatible family of Feller
semigroups on $\cP(M)$ satisfying (\ref{mapcond}). \elem
\prf Since the semigroups $\P^{(n)}_t$ are Feller, the
semigroups $\Pi^{(k)}_t$ are also Feller~: for all $g$ in $\cA_k$,
then $\Pi^{(k)}_tg$ is continuous and $\lim_{t\to 0}\Pi^{(k)}_t g=g$
and these properties extend to every continuous functions.

It is clear that the family of semigroups $(\Pi^{(k)}_t,~k\geq 1)$ is
compatible (in the sense given in section \ref{hypsg}). Thus
$(\Pi^{(k)}_t,~k\geq 1)$ is a compatible family of Feller semigroups
on $\cP(M)$. We denote $\Pi^{(2)}_{(\mu,\nu)}$ the law of the Markov
process associated with  $\Pi^{(2)}_t$ starting from $(\mu,\nu)$ and
we denote this process by $(\mu_t,\nu_t)$.

For $g\in\cA_1$ in the form (\ref{defg}), $t\geq 0$ and $\mu\in\cP(M)$, we have
$$ \Pi^{(2)}_t g^{\otimes 2}(\mu,\mu) 
= \langle\P^{(2n)}_t f^{\otimes 2},\mu^{\otimes 2n}\rangle
= \Pi^{(1)}_t g^2 (\mu). $$
Thus (\ref{mapcond}) is satisfied for $g\in\cA_1$ and this extends to
$C(\cP(M))$. \qed

\subsubsection{Proof of the first part of theorem \ref{mainthmk}.}
Using theorem \ref{mainthm} we construct $(\Q_t)_{t\geq 0}$ a Feller
convolution semigroup on $(\Gamma,\cG)$. Let $\cJ_t$ be a measurable
presentation of $\Q_t$. Set $\nu_t=(\delta\circ\cJ_t)^*\Q_t$.

\blem\label{lemlin0} For all $\mu\in\cP(M)$ and all $t\geq 0$,
\beq\Q_t(d\Phi)-\hbox{a.s.},\qquad \Phi(\mu)=\mu(\delta\circ\cJ_t(\Phi)).\eeq
And for all $t\geq 0$, $\cI^*(\nu_t)=\Q_t$.
\elem
\prf For all $f\in C(M)$, set $g(\mu)=\mu f$, then
\beqarr \E[(\mu(\delta\circ\cJ_t(\Phi))f-\Phi(\mu)f)^2]
&=& \E\left[\left(\int g(\Phi(\delta_x)) \mu(dx)-g(\Phi(\mu))\right)^2\right]\\
&=& \int \Pi^{(2)}_tg^{\otimes 2}(\delta_x,\delta_y)
     \mu(dx)\mu(dy)  + \Pi^{(2)}_tg^{\otimes 2}(\mu,\mu)\\
&& - \quad 2 \int \Pi^{(2)}_tg^{\otimes 2}(\delta_x,\mu)\mu(dx).
\eeqarr
Since for all $\mu$ and $\nu$ in $\cP(M)$,
$$\Pi^{(2)}_tg^{\otimes 2}(\mu,\nu) =
\int \P^{(2)}_t f^{\otimes 2}(x,y)~\mu(dx)\nu(dy),$$
we get $\E[(\mu(\delta\circ\cJ_t(\Phi))f-\Phi(\mu)f)^2]=0$. This proves the
lemma. \qed

\medskip Lemma \ref{lemlin0} implies that
$\nu=(\nu_t)_{t\geq 0}$ is a Feller convolution semigroup on $(E,\cE)$
(we apply proposition \ref{qnu}) and (\ref{eqrelconv}) holds. This
proves the first part of theorem \ref{mainthmk}. \qed

\subsubsection{Proof of the second part of theorem \ref{mainthmk}.}
Suppose now we are given $\nu=(\nu_t)_{t\geq 0}$ a Feller
convolution semigroup on $(E,\cE)$. For $t\geq 0$, set
$\Q_t=\cI^*(\nu_t)$. Then $\Q=(\Q_t)_{t\geq 0}$ is a Feller
convolution semigroup on $(\Gamma,\cG)$. Using theorem \ref{mainthm},
we construct $\P_\Q$ the law of a stochastic flow of mappings on
$\cP(M)$ associated with $\Q$. 
Let $(\Phi_{s,t},~s\leq t)$ be a stochastic flow of mappings of law $\P_\Q$. 
For $t\ge 0$, there is a measurable presentation $\cJ_t$ of $\Q_t$ and set $\mathfrak{p}_t=\delta\circ\cJ_t\circ\cI$. Then $\mathfrak{p}_t$ is a measurable presentation of $\nu_t$. 
 For $s\leq t$,
set $K_{s,t}=\delta\circ\cJ_{t-s} (\Phi_{s,t})$.

We now show that $K=(K_{s,t},~s\leq t)$ is a stochastic flow of
kernels. Note that the law of $K_{s,t}$ is $\nu_{t-s}$. Thus it is
easy to check that $K$ satisfies {\bf (b)}, {\bf (c)}, {\bf (d)} and
{\bf (e)}. In order to show {\bf (a)}, we use

\blem\label{lemlin} For all $\mu\in\cP(M)$ and all $s\leq t$,
\beq\P-\hbox{a.s.},\qquad \mu K_{s,t}=\Phi_{s,t}(\mu).\eeq
\elem
\prf For all $f\in M$, set $g(\mu)=\mu f$, then like in the proof of
lemma \ref{lemlin0},
\beqarr \E[(\mu K_{s,t}f-\Phi_{s,t}(\mu)f)^2]
&=& \E\left[\left(\int g(\Phi_{s,t}(\delta_x)) \mu(dx)-g(\Phi_{s,t}(\mu))\right)^2\right]\\
&=& \int \Pi^{(2)}_{t-s}g^{\otimes 2}(\delta_x,\delta_y)
     \mu(dx)\mu(dy)  + \Pi^{(2)}_{t-s}g^{\otimes 2}(\mu,\mu),\\
&& - \quad 2 \int \Pi^{(2)}_{t-s}g^{\otimes 2}(\delta_x,\mu)\mu(dx)\\
&=& 0. 
\eeqarr
This proves the lemma. \qed

\brem{\em  If $\Lambda$ is a $\cP(M)$-valued random variable independent of $\Phi_{s,t}$, then (using that $(\mu,K)\mapsto \mu(\mathfrak{p}_{t-s}( K_{s,t}))$ and $(\mu,\Phi)\mapsto \cJ_{t-s}( \Phi_{s,t})(\mu)$ are measurable)
$$\P-\hbox{a.s.},\qquad \Lambda (\mathfrak{p}_{t-s}( K_{s,t}))=\cJ_{t-s}(\Phi_{s,t})(\Lambda).$$}\erem 

Let $s\leq u\leq t$ and $\mu\in\cP(M)$. Lemma \ref{lemlin}
and the cocycle property for $\Phi$ imply that a.s.,
$$\mu K_{s,t} = \Phi_{s,t}(\mu) = \cJ_{t-u}(\Phi_{u,t})\circ\Phi_{s,u}(\mu).$$
Lemma \ref{lemlin} implies that a.s.,
$\cJ_{t-u}(\Phi_{u,t})\circ\Phi_{s,u}(\mu)=\cJ_{t-u}(\Phi_{u,t})(\mu K_{s,u})$.
Fubini's theorem, lemma \ref{lemlin} and the fact that $\mu K_{s,u}$
and $\Phi_{u,t}$ are independent imply that a.s.,
$$\cJ_{t-u}(\Phi_{u,t})(\mu K_{s,u})=\mu K_{s,u}(\mathfrak{p}_{t-u}( K_{u,t})).$$
This proves {\bf (a)}, i.e., a.s. $\mu K_{s,t}=\mu K_{s,u}(\mathfrak{p}_{t-u}( K_{u,t}))$.
We let $\P_\nu$ be the law of $K$. Then $T_h^*(\P_\nu)=\P_\nu$. The
rest of the proof is similar to the end of the proof of theorem
\ref{mainthm}. \qed

\subsection{Sampling the flow.}\label{cadlag}
Let $(K_{s,t},~s\leq t)$ be a measurable stochastic flow of kernels defined on a
probability space $(\O,\cA,\P)$ and $(T_h)_{h\in\RR}$ a one-parametric
group of transformations of $\O$ preserving $\P$ and such that
$K_{s,t}\circ T_h=K_{s+h,t+h}$. In this section, we construct on an
extension of $(\O,\cA,\P)$ a random path $X_t$ starting at $x$ such
that for all positive $t$, 
\beq K_{0,t}f(x)=\E[f(X_t)|\cA]. \eeq 

For $x\in M$ and $\omega\in\Omega$, by Kolmogorov's theorem, we define
on $M^{\RR^+}$, a probability $\P^0_{x,\omega}$ such that
\beq \E^0_{x,\omega}\left[\prod_{i=1}^nf_i(X^0_{t_i})\right] =
K_{0,t_1}(f_1(K_{t_1,t_2}f_2(\cdots(f_{n-1}K_{t_{n-1},t_n}f_n))))(x),
\eeq
for all $f_1,\dots,f_n$ in $C(M)$, $0<t_1<t_2<\cdots<t_n$.

With $\P$ and $\P^0_{x,\omega}$, we construct a probability
$\P^0_x(d\omega,d\omega') = 
\P(d\omega)\otimes \P^0_{x,\omega}(d\omega')$ on $\Omega\times
M^{\RR^+}$. Then, on the probability space 
$(\Omega\times M^{\RR^+},\cA\otimes{\cal B}(M)^{\otimes\RR^+},\P^0_x)$,
the random process $(X^0_t,~t\geq 0)$, defined by
$X^0_t(\omega,\omega')=\omega'(t)$, is a Markov process starting at
$x$ with semigroup $\P^{(1)}_t$ since
\beq \E^0_{x}\left[\prod_{i=1}^nf_i(X^0_{t_i})\right] =
\P^{(1)}_{t_1} (f_1 (\P^{(1)}_{t_2-t_1}f_2 (\cdots
(f_{n-1}\P^{(1)}_{t_{n}-t_{n-1}}f_n))))(x), \eeq 
for all $f_1,\dots,f_n$ in $C(M)$, $0<t_1<t_2<\cdots<t_n$.

Therefore, there is a c\`adl\`ag (or continuous when $\P^{(1)}_t$ is the
semigroup of a continuous Markov process) modification $X=(X_t,~t\geq
0)$ of $(X^0_t,~t\geq 0)$. Let now $\P_{x,\omega}$ be the law of $X$
knowing $\cA$. It is a law on $D(\RR^+,M)$, the space of c\`adl\`ag
functions (or $C(\RR^+,M)$ when $\P^{(1)}_t$ is the semigroup of a
continuous Markov process). Equipped with the Skorohod topology (see
\cite{maison} or \cite{biling}), $D(\RR^+,M)$ becomes a Polish space
(respectively $C(\RR^+,M)$ is equipped with the topology of uniform
convergence on every compact on $\RR^+$). 

On the probability space $(\Omega\times D(\RR^+,M) , \cA\otimes
{\cal B}(D(\RR^+,M)) , \P_x)$  (respectively on $(\Omega\times
C(\RR^+,M) ,  \cA\otimes {\cal B}(C(\RR^+,M)) , \P_x)$), where
$\P_x(d\omega,d\omega') = \P(d\omega)\otimes \P_{x,\omega}(d\omega')$,
let $X$ be the random process $X(\omega,\omega')=\omega'$. Then $X$ is
a c\`adl\`ag (respectively continuous) process and
\beqar \E_x\left[\left.\prod_{i=1}^nf_i(X_{t_i})\right|\cA\right] &=&
\E_{x,\omega}\left[\prod_{i=1}^nf_i(X_{t_i})\right]\\
&=& K_{0,t_1}(f_1(K_{t_1,t_2}f_2(\cdots(f_{n-1}K_{t_{n-1},t_n}f_n))))(x),\non
\eeqar
where $\E_x$ denotes the expectation with respect to $\P_x$.

\medskip
Let $(K'_{s,t},~s\leq t)$ be the stochastic flow of kernels defined on
$(\Omega,\cA,\P)$ by 
\beq K'_{s,t}f(x,\omega)=K'_{0,t-s}f(x,T_s\omega)\eeq
where
\beq K'_{0,t}f(x)=\E_{x}[f(X_t)|\cA] = 
\int f(X_t(\omega,\omega'))\P_{x,\omega}(d\omega') \eeq
for $f\in C(M)$, $x\in M$.
Then $(K'_{s,t},~s\leq t)$ is a c\`adl\`ag in $t$ (respectively continuous
in $t$) modification of $(K_{s,t},~s\leq t)$.

\brem  The concept of sampling will be used in section \ref{solveSDE}.\erem

Replacing $K_{0,t}$ by $K_{0,t}^{\otimes n}$ and $\P^{(1)}_t$ by
$\P^{(n)}_t$ in the above, we obtain a random process $X^{(n)}$ in
$M^n$ which represents an $n$-sampling of the flow. The coordinates of
$X^{(n)}$ are independent given the flow $K$.

Let $(x_i)_{i\geq 1}$ be a sequence in $M$. For $\o\in\O$, let
$\P_{x_1,\dots,x_n,\o}=\otimes_{i=1}^n \P_{x_i,\o}$,
$\P_{(x_i)_{i\geq 1},\o}=\otimes_{i\geq 1}^n \P_{x_i,\o}$,
$\P_{x_1,\dots,x_n}(d\o,d\o'_1,\dots,d\o'_n)=
\P(d\o)\otimes\P_{x_1,\dots,x_n,\o}(d\o'_1,\dots,d\o'_n)$ and
$\P_{(x_i)_{i\geq 1}}(d\o,d\o')=
\P(d\o)\otimes\P_{(x_i)_{i\geq 1},\o}(d\o')$. Then the process
$X^{(n)}(\o,\o')=(\o'_1,\dots,\o'_n)$ defines an $n$-sampling of the
flow (under $P_{x_1,\dots,x_n}$ or $\P_{(x_i)_{i\geq 1}}$). Let
$X^i(\o,\o')=\o'_i$. Then, under $\P_{(x_i)_{i\geq 1}}$, the sequence
$(X^i)_{i\geq 1}$ are independent conditionaly to $\cA$. Moreover, if
for all $i\geq 1$, $x_i=x$, this sequence is identically distributed
and the law of large numbers implies that for all $f\in C_0(M)$,
$\frac{1}{n}\sum_{i=1}^n f(X^i_t)$ converges a.s. towards
$\E_x[f(X^1_t)|\cA]=K_{0,t}f(x)$.

Since, under $\P_{(x_i)_{i\geq 1}}$, $X^{(n)}$ is equal in law to the
$n$-point motion of $K$ starting from $(x_1,\dots,x_n)$, if for all
$n\geq 1$, we let $X^{(n)}$ denote the $n$-point motion starting from
$(x,\dots,x)$, we have that $\frac{1}{n}\sum_{i=1}^n f(X^i_t)$
converges in law towards $K_{0,t}f(x)$ for all $f\in C^0(M^n)$. This
gives an intuitive way to recover $K_{0,t}(x)$ out of the $n$-point
motions. 

\setcounter{equation}{0}
\section{Noise and stochastic flows.} \label{secnoise}
\subsection{The noise generated by a stochastic flow of kernels.}\label{noisedef}
The definition of a noise we give here is very close to the one
given by Tsirelson in \cite{Tsirelson}.
\begin{definition}
A noise consists of a separable probability space 
$(\Omega,\cA,\P)$, a one-parametric group $(T_h)_{h\in\RR}$ of
$\P$-preserving $L^2$-continuous transformations of $\Omega$ and a
family $\{{\cal F}_{s,t},~-\infty\leq s\leq t\leq\infty\}$ of 
sub-$\sigma$-fields of $\cA$ such that
\bdes
\ita $T_h$ sends ${\cal F}_{s,t}$ onto ${\cal F}_{s+h,t+h}$
for all $h\in\RR$ and all $s\leq t$,
\itb ${\cal F}_{s,t}$ and ${\cal F}_{t,u}$ are independent 
for all $s\leq t\leq u$,
\itc ${\cal F}_{s,t}\vee {\cal F}_{t,u}={\cal F}_{s,u}$ 
for all $s\leq t\leq u$. \edes
Moreover, we will assume that, for all $s\leq t$, $\mathcal{F}_{s,t}$
contains all $\P$-negligible sets of $\mathcal{F}_{-\infty,\infty}$,
denoted $\mathcal{F}$. \end{definition}

In the following, $(\O^0,\cA^0,\P_\nu)$ denotes the canonical
probability space of a stochastic flow of kernels on $M$, a locally compact separable metric space, associated with a
Feller convolution semigroup $\nu$. And $K^0=(K^0_{s,t},~s\leq t)$ denotes
this canonical flow. When this stochastic flow is induced by a
flow of maps, one can take for $(\O^0,\cA^0,\P_\nu)$, the canonical
probability space associated to this stochastic flow of mappings.

\medskip
For all $-\infty\leq s\leq t\leq\infty$, let $\cF^{\nu}_{s,t}$ be the
sub-$\sigma$-field of $\cA^0$ generated by the random variables
$K^0_{u,v}$ for all $s\leq u\leq v\leq t$ completed by all
$\P_\nu$-negligible sets of $\cA^0$.
Then the cocycle property of $K^0$ implies that 
$N_\nu:=(\O^0,\cA^0,(\cF^\nu_{s,t})_{s\leq t},\P_\nu,(T_h)_{h\in\RR})$
is a noise ($T_h$ is $L^2$-continuous because of the Feller
property). We call it the noise generated by the canonical flow $K^0$.

\begin{definition} Let $\nu$ be a Feller convolution semigroup,
$N=(\Omega,\cA,({\cal F}_{s,t})_{s\leq t} , \P , (T_h)_{h\in\RR})$ be
a noise and $K$ be a measurable stochastic flow of kernels of law $\P_\nu$
defined on $(\Omega,\cA,\P)$ such that for all $s<t$, $K_{s,t}$ is
${\cal F}_{s,t}$-measurable and for all $h\in\RR$,
\beq K_{s+h,t+h} = K_{s,t}\circ T_h, \qquad \hbox{a.s.}\eeq
We will call $(N,K)$ an extension of the noise $N_\nu$. \end{definition}

Let $(N_1,K_1)$ and $(N_2,K_2)$ be two extensions of the
noise $N_\nu$. Let $\Omega = \Omega_1 \times \Omega_2$, $\mathcal{A} =
\mathcal{A}_1\otimes\mathcal{A}_2$ and $\P$ be the probability
measure on $(\Omega,\cA)$ defined by 
\beq \E[Z] = \int \E_1[Z_1|K_1=K]\E_2[Z_2|K_2=K] ~\P_\nu(dK),\eeq
for any bounded random variable
$Z(\omega_1,\omega_2)=Z_1(\omega_1)Z_2(\omega_2)$. Let
$(T_h)_{h\in\RR}$ be the one-parametric group of $\P$-preserving
transformations of $\Omega$ defined by
$T_h(\omega_1,\omega_2)=(T^1_h(\omega_1),T^2_h(\omega_2))$. For all
$s<t$, let $\mathcal{F}_{s,t} =
\mathcal{F}^1_{s,t}\otimes\mathcal{F}^2_{s,t}$. Then
$N:=(\Omega,\cA,({\cal F}_{s,t})_{s\leq t} , \P ,
(T_h)_{h\in\RR})$ is a noise. And if $K$ denotes the random variable
$K(\omega_1,\omega_2)=K_1(\omega_1)(=K_2(\omega_2)$ $\P$-a.s.),
then $(N,K)$ is an extension of $N_\nu$. We will call $(N,K)$ the
product of the extensions $(N_1,K_1)$ and $(N_2,K_2)$. Note that $N_1$
and $N_2$ are isomorphic to sub-noises of $N$.

\subsection{Filtering kernels}
Let $K$ be a random kernel defined on a probability space $(\Omega,\mathcal{A},\P)$ and let $\bar{\mathcal{A}}$ be a sub $\sigma$-field. Denote by $\nu$ the law of $K$ and set $\Q=\mathcal{I}^*(\nu)$. Then $\Q$ is a law on $(\Gamma,\mathcal{G})$.

\blem\label{lem:321}
Suppose that 
\begin{equation}\label{eq:regcond}
\lim_{\rho(\mu_1,\mu_2)\to 0} \E[\rho(\mu_1 K, \mu_2 K)^2]=0.
\end{equation}
Then $\nu$ is regular and there is $\bar{K}$ an  $\bar{\mathcal{A}}$-measurable random kernel, with law denoted by $\bar{\nu}$, such  that
\begin{enumerate}[label=(\roman*)]
\item For all $\mu \in \mathcal{P}(M)$, $\mu \bar{K}=\E[\mu K|\bar{\mathcal{A}}]$;
\item $\bar{\nu}$ is regular;
\item Let $\mathfrak{p}$ and $\bar{\mathfrak{p}}$ be  measurable presentations respectively of $\nu$ and of $\bar{\nu}$. Let 
$\Lambda$ be a $\mathcal{P}(M)$-valued random variable and
$\bar{\mathcal{A}}'$ be another sub $\sigma$-field of $\mathcal{A}$.
Suppose that  $\sigma(\Lambda)\vee\bar{\mathcal{A}}'$ and $\sigma(K)\vee \bar{\mathcal{A}}$ are independent. Then, if $\bar{\Lambda}=\E[\Lambda|\bar{\mathcal{A}}']$,
\begin{equation}\label{eq:conditioning}
\bar{\Lambda} \bar{\mathfrak{p}}(\bar{K})=\E[\Lambda \mathfrak{p}(K) | \bar{\mathcal{A}}'\vee \bar{\mathcal{A}})].
\end{equation}
\end{enumerate}
\elem
\begin{proof} In the proof of this lemma, we will use regular probability measures:
Let $Y$ be a random variable taking its values in a Borel space $S$ equipped with its Borel $\sigma$-field $\mathcal{S}$. 
Applying theorems 6.3 and 6.4 in \cite{kal}, 
there is a regular probability measure $\nu:\Omega\times \mathcal{S} \to [0,1]$ such that
\begin{itemize}
\item[(i)] $\nu(\omega,\cdot)$ is a probability measure on $\mathcal{S}$;
\item[(ii)] $\nu(\cdot,A)$ is $\bar{\mathcal{A}}$-measurable for all $A\in \mathcal{S}$;
\end{itemize}
and such that for all bounded measurable function $g$, and $\Lambda$ an $\bar{\mathcal{A}}$-measurable random variable,
$$\E[g(\Lambda,Y)|\bar{\mathcal{A}}]=\int_{\mathcal{S}} g(\Lambda(\omega),y) \nu(\omega,dy).$$
In other words, $\nu$ is the conditional law of $Y$ given $\mathcal{G}$. (Regular conditional distributions  were actually implicitly used in section 2.6).

In general, if $Y$ is a $\mathcal{P}(M)$-valued random variable of law $\P_Y$, $f$ an element of $C_c(M)$, $\varphi(f)=\E[Y f]$ is a positive linear form on $C_c(M)$. By the Riesz-Markov-Kakutani representation theorem there is $\mu\in \mathcal{P}(M)$ such that $\E[Y f] = \mu f.$
Then $\mu$ can be denoted$\int_{\mathcal{P}(M)} y \P_Y(dy)$.
Using regular probability measures as above, one can define $\E[Y|\bar{\mathcal{A}}]=\int_{\mathcal{P}(M)} y \nu(\cdot,dy)$, with $\nu$ the conditional law of $Y$ given $\bar{\mathcal{A}}$.

\smallskip
We now prove the lemma.
Let us first show that $\nu$ is regular, i.e. that $\Q$ is regular.
Let $(\mu_k,\,k\in\mathbb{N})$ be a dense sequence in $\mathcal{P}(M)$ and let $j:\Gamma\to \mathcal{P}(M)^{\otimes \mathbb{N}}$ be the measurable function defined by $j(\Phi)=(\Phi(\mu_k), k\in \mathbb{N})$.
Following section 1, \eqref{eq:regcond} allows to construct $i:\mathcal{P}(M)^{\otimes \mathbb{N}}\to \Gamma$ a measurable function such that $\mathcal{J}= i\circ j$ is a measurable presentation of $\Q$. This shows that $\Q$ and therefore $\nu$ are regular. 
With this construction, as is shown in section 2, $\delta\circ \mathcal{J}$ is measurable, $(\delta\circ \mathcal{J})^*(\Q)=\nu$ and $\mathfrak{p}=\delta\circ\mathcal{J}\circ \mathcal{I}$ is a measurable presentation of $\nu$. 

Set $\P^{(\infty)}=j^*(\Q)$.
Then $i^*(\P^{(\infty)})=\Q$ and  $(\delta\circ i)^*(\P^{(\infty)})=\nu$ (indeed $\delta\circ i = (\delta\circ \mathcal{J})\circ i$ is measurable and $(\delta\circ i)^*(\P^{(\infty)})=(\delta\circ \mathcal{J}\circ i)^*(\P^{(\infty)})=(\delta\circ \mathcal{J})^*(\Q)=\nu$).

Set $Y=j\circ\mathcal{I}(K)$. Then $Y$ is a random variable of law $\P^{(\infty)}$.
Since $M^\mathbb{N}$ is a Polish space, one can define $\bar{\P}^{(\infty)}(\omega,dy)$ the conditional law of $Y$ given $\bar{\mathcal{A}}$.
In particular for all bounded measurable function $g$ and $\Lambda$ an $\bar{\mathcal{A}}$-measurable random variable,
$$\E[g(\Lambda,Y)|\bar{\mathcal{A}}]=\int_{\mathcal{P(M)}^{\mathbb{N}}} g(\Lambda(\omega),y) \bar{\P}^{(\infty)}(\omega,dy).$$

Therefore, if $K'=\mathfrak{p}(K)=
\delta\circ i(Y)$, we have for all bounded measurable function $g$ and $\Lambda$ an $\bar{\mathcal{A}}$-measurable random variable,
$$\E[g(\Lambda,\mathfrak{p}(K))|\bar{\mathcal{A}}]=\int_{\mathcal{P}(M)^{\mathbb{N}}} g(\Lambda(\omega),(\delta\circ i)(y)) \bar{\P}^{(\infty)}(\omega,dy).$$

For all $\omega\in \Omega$, define $\bar{\nu}(\omega,\cdot)=(\delta\circ i)^*(\bar{\P}^{(\infty)}(\omega,\cdot))$. Then $\bar{\nu}(\omega,d \mathsf{k})$ is a regular probability measure and we have for all bounded measurable function $g$ and $\Lambda$ an $\bar{\mathcal{A}}$-measurable random variable,
$$\E[g(\Lambda,\mathfrak{p}(K))|\bar{\mathcal{A}}]=\int_{E} g(\Lambda(\omega),\mathsf{k}) \bar{\nu}(\omega,d\mathsf{k}).$$

For $x\in M$, set 
$$\bar{K}(x)=\int_{E} \mathfrak{p}(\mathsf{k})(x) \bar{\nu}(\omega,d\mathsf{k}).$$
Then $\bar{K}$ is a random kernel is $\bar{\mathcal{A}}$-measurable and a.s. $\bar{K}(x)=\int_{E} \mathsf{k}(x) \bar{\nu}(\omega,d\mathsf{k})$.
Denote by $\bar{\nu}$ the law of $\bar{K}$ and set $\bar{Q}=\mathcal{I}^*(\bar{\nu})$.
Let us now show (i), (ii) and (iii) are satisfied. For $\mu \in \mathcal{P}(M)$, a.s. $\mu \mathfrak{p}(K)=\mu K$ and therefore, a.s.
$$\E[\mu K |\bar{\mathcal{A}}]
=\E[\mu\mathfrak{p}(K)|\bar{\mathcal{A}}]=\int_{E} \mu\mathsf{k}\,\bar{\nu}(\omega,d\mathsf{k})=\mu \bar{K}$$
where in the last equality we have used Fubini's theorem. This proves (i).

To prove (ii), observe that for $f\in C(M)$ and $(\mu_1,\mu_2)\in \mathcal{P}(M)^2$,
\begin{align*}
\E[(\mu_2 \bar{K}f-\mu_1 \bar{K}f)^2]
&= \E\left[\E[\mu_2 K f - \mu_1 Kf |\bar{\mathcal{A}}]^2\right]\\
&\le \E\left[\left(\mu_2 K f - \mu_1 Kf\right)^2\right].
\end{align*}
This implies that $\bar{\nu}$ satisfies \eqref{eq:regcond} and as a consequence, by construction, that $\mathcal{J}$ and $\mathfrak{p}$ are also measurable presentations respectively of $\bar{\Q}$ and of $\bar{\nu}$.

Let us finally prove (iii). Note first that item (i) implies that   \eqref{eq:conditioning} holds if $\Lambda$ is not random (indeed, for all $\mu\in \mathcal{P}(M)$, a.s. $\mu K=\mu \mathfrak{p}(K)$ and $\mu \bar{K}=\mu \bar{\mathfrak{p}}(\bar{K})$). 
Let now $\bar{Z}'$ and $\bar{Z}$ be bounded  random variables, respectively $\bar{\mathcal{A}}'$-measurable and $\bar{\mathcal{A}}$-measurable. 
Then (using that $\sigma(\Lambda)\vee \bar{\mathcal{A}}'$ and $\sigma(K)\vee \bar{\mathcal{A}}$ are independent  and using the notation $\P_Z$ for the law of a random variable $Z$),
\begin{align*}
\E[\Lambda\mathfrak{p}(K) \,  \bar{Z}' \bar{Z}]
&= \int_{\mathcal{P}(M)\times \mathbb{R}} \E[\mu \mathfrak{p}(K) \, \bar{Z}]\, \bar{z}' \;\P_{(\Lambda,\bar{Z}')}(d\mu,d\bar{z}')\\
&= \int_{\mathcal{P}(M)\times\mathbb{R}} \E[\mu \bar{\mathfrak{p}}(\bar{K})\, \bar{Z}] \, \bar{z}' \;\P_{(\Lambda,\bar{Z}')}(d\mu,d\bar{z}')\\
&=\E[\Lambda\bar{\mathfrak{p}}(\bar{K})\, \bar{Z}' \bar{Z}]\\
&= \int_{E\times\mathbb{R}}
\E[\Lambda \bar{\mathfrak{p}}(\bar{\mathsf{k}})\,\bar{Z}']\,\bar{z} \;\P_{(\bar{K},\bar{Z})}(d\bar{\mathsf{k}},d\bar{z})\\
&= \int_{E\times\mathbb{R}}
\E[\bar{\Lambda} \bar{\mathfrak{p}}(\bar{\mathsf{k}})\, \bar{Z}']\,\bar{z} \;\P_{(\bar{K},\bar{Z})}(d\bar{\mathsf{k}},d\bar{z})\\
&= \E[\bar{\Lambda }\bar{\mathfrak{p}}(\bar{K})\, \bar{Z}' \bar{Z}]
\end{align*}
This implies (iii).
\end{proof}

\subsection{Filtering by a sub-noise.}\label{filtering}
Let $\bar{N}$ be a sub-noise of an extension $(N,K)$ of $N_\nu$,
i.e. $\bar{N}$ is a noise $(\Omega,\cA,({\bar{\cal F}_{s,t}})_{s\leq
t},\P,(T_h)_{h\in\RR})$) such that $\bar{{\cal F}}_{s,t}\subset
{\cal F}_{s,t}$ for all $s\leq t$.

\brem  A sub-noise is characterized by
$\bar{\mathcal{F}}_{-\infty,\infty}$, denoted
$\bar{\mathcal{F}}$. This $\sigma$-field has to be stable under $T_h$,
to contain all $\P$-negligible sets of $\cF$, and be such that $\bar{\cF} =
(\bar{\cF}\cap \cF_{-\infty,0})\vee (\bar{\mathcal{F}}\cap
\mathcal{F}_{0,\infty})$. \erem

For all $n\geq 1$, let $\bar{\P}^{(n)}_t$ be the operator acting on
$C(M^n)$ defined by 
\beq\label{defpbar}
\bar{\P}^{(n)}_t(f_1 \otimes\cdots\otimes f_n)(x_1,\dots,x_n) = 
\E\left[\prod_{i=1}^n \E[K_{0,t}f_i(x_i)|\bar{\cF}_{0,t}]\right],\eeq
for all $x_1,\dots,x_n$ in $M$ and all $f_1,\dots,f_n$ in $C(M)$.
\blem The family $(\bar{\P}^{(n)}_t,~n\geq 1)$ is a compatible family
of Feller semigroups. \elem

\prf The semigroup property of $\bar{\P}^{(n)}_t$ follows
directly from the independence of the increments of the flow. The
Markovian property and in particular the positivity property holds
since for all $h\in C(M^n)$, 
\beq \bar{\P}^{(n)}_th(x_1,\dots,x_n) =
\E[\langle h,\otimes_{i=1}^n\E[K_{0,t}(x_i)|\bar{{\cal F}}_{0,t}]\rangle].\eeq
From this, it is clear that $(\bar{\P}^{(n)}_t,~n\geq 1)$ is a compatible family
of Markovian semigroups respectively acting on $C(M^n)$.

It remains to prove the Feller property. For all continuous
functions $f_1,\dots,f_n$, $h=f_1\otimes\cdots\otimes f_n$,
$x=(x_1,\dots,x_n)$ and $y=(y_1,\dots,y_n)$ in $M^n$, for $M$ large
enough,
\beqar |\bar{\P}^{(n)}_t h(x)-\bar{\P}^{(n)}_t h(y)| &\leq&
M\sum_{i=1}^n \E[(\E[K_{0,t}f_i(x_i)-K_{0,t}f_i(y_i)|
\bar{{\cal F}}_{0,t}])^2]^{\frac{1}{2}}\non\\ 
&\leq& M\sum_{i=1}^n
\E[(K_{0,t}f_i(x_i)-K_{0,t}f_i(y_i))^2]^{\frac{1}{2}} \eeqar
which converges towards 0 as $y$ tends to $x$ since {\bf (e)} in
definition \ref{defsfk} is satisfied.

We also have, for all $h=f_1\otimes\cdots\otimes f_n$ and
$x=(x_1,\dots,x_n)$ in $M^n$, for $M$ large enough,
\beqar |\bar{\P}^{(n)}_t h(x)-h(x)| &\leq& M\sum_{i=1}^n
\E[(\E[(K_{0,t}f_i(x_i)-f_i(x_i)|
\bar{{\cal F}}_{0,t}])^2]^{\frac{1}{2}}\non\\ 
&\leq& M\sum_{i=1}^n \E[(K_{0,t}f_i(x_i)-f_i(x_i))^2]^{\frac{1}{2}}
\eeqar 
which converges towards 0 as $t$ tends to 0 since {\bf (d)} in
definition \ref{defsfk} is satisfied. Hence, for all function $h\in
C(M^n)$ such that $h$ is a linear combinaison of functions of the type
$f_1\otimes\cdots\otimes f_n$, we have $\bar{\P}^{(n)}_th$ is continuous
and $\lim_{t\to 0} \bar{\P}^{(n)}_th(x)=h(x)$ for all $x\in M^n$. This
extends to all functions $h\in C(M^n)$. \qed

\medskip
Let us denote by $\bar{\nu}=(\bar{\nu}_t)_{t\geq 0}$ the Feller
convolution semigroup on $(E,{\cal E})$ associated with
$(\bar{\P}^{(n)}_t,~n\geq 1)$. Note that the one-point motion of $\nu$
and $\bar{\nu}$ are the same, i.e. $\bar{\P}^{(1)}_t=\P^{(1)}_t$.

\blem \label{lemfiltr}
Let $(N,K)$ be an extension of $N_\nu$ and $\bar{N}$ be a
sub-noise of $N$. Then there exists $\bar{K} = (\bar{K}_{s,t},~s\leq
t)$ a stochastic flow of kernels of law $\P_{\bar{\nu}}$ such that
$(\bar{N},\bar{K})$ is an extension of $N_{\bar{\nu}}$ and
\beq \bar{K}_{s,t}f(x)=\E[K_{s,t}f(x)|\bar{\cal F}_{s,t}] =
\E[K_{s,t}f(x)|\bar{\mathcal{F}}] \qquad \P-\hbox{a.s.} \eeq
for all $s\leq t$, $x\in M$ and $f\in C_0(M)$.
We say $\bar{K}$ is obtained by filtering $K$ with respect to
$\bar{N}$. 
\elem
\begin{proof}
For $s\le t$, $K_{s,t}$ satisfies \eqref{eq:regcond} and there is  $\bar{K}_{s,t}$, the random kernel obtained by filtering $K_{s,t}$ with respect to $\bar{\mathcal{F}}_{s,t}$
(Remark: one can first define $\bar{K}_ {0,t}$ for $t\ge 0$ and set $\bar{K}_{s,t}=\bar{K}_{0,t-s}\circ T_s$ in order to ensure that $\bar{K}_{s,t}\circ T_h=\bar{K}_{s+h,t+h}$).
For all $s\le t$, the law of $\bar{K}_{s,t}$ is $\bar{\nu}_{t-s}$ and $\bar{K}$ is a stochastic flow of kernels. 
Indeed, except for (a), all the conditions given in definition 2.3  are easy to check. 
To prove the cocycle property (a), let $\mu \in  \mathcal{P}(M)$ and $s\le u\le t$. Then
\begin{align*}
\mu  \bar{K}_{s,t}
&= \E[\mu K_{s,t} | \bar{\mathcal{F}}_{s,t}]\\
&= \E[\mu K_{s,u}\mathfrak{p}_{t-u}( K_{u,t}) | \bar{\mathcal{F}}_{s,u}
\vee \bar{\mathcal{F}}_{u,t}].
\end{align*}
Since $\sigma(\mu K_{s,u})\vee \bar{\mathcal{F}}_{s,u}(\subset \mathcal{F}_{s,u})$ and $\sigma(\mu K_{u,t})\vee \bar{\mathcal{F}}_{u,t}(\subset \mathcal{F}_{u,t})$ are independent, one can apply (iii) of lemma \ref{lem:321} and easily obtain (a).
\end{proof}

\begin{definition}\label{defdomination} Given two Feller convolution
semigroups on $(E,{\cal E})$, $\nu^1$ and $\nu^2$, we say that $\nu^1$
dominates (respectively weakly dominates) $\nu^2$, denoted $\nu^1\succeq \nu^2$
(respectively $\nu^1\succeq^w \nu^2$), if there exists a sub-noise
of $N_{\nu^1}$ (respectively of an extension $(N^1,K^1)$ of
$N_{\nu^1}$) such that $\P_{\nu^2}$ is the law of the flow obtained by
filtering the canonical flow of law $\P_{\nu^1}$ (respectively by
filtering $K^1$) with respect to this sub-noise.  \end{definition}

Notice that in lemma \ref{lemfiltr}, $\nu$ weakly dominates
$\bar{\nu}$ and $\nu$ dominates $\bar{\nu}$ if $\bar{N}$ is a
sub-noise of $N_\nu$. Note that the domination relation is in fact an
extension of the notion of barycenter.

\blem\label{propfiltr} Let $\nu$ and $\bar{\nu}$ be two Feller
convolution semigroups such that $\nu$ dominates $\bar{\nu}$. Let
$(N,K)$ be an extension of $N_\nu$. Let $\tilde{N}_\nu$ be the
sub-noise (isomorphic to $N_\nu$) of $N$ generated by $K$. Then there
exists a sub-noise $\bar{N}$ of $\tilde{N}_\nu$ such that
$\P_{\bar{\nu}}$ is the law of the flow obtained by filtering $K$ with
respect to $\bar{N}$. \elem

\prf Let $N_\nu:=(\Omega^0 , \cA^0 , ({\cal F}^\nu_{s,t})_{s\leq t} ,
\P_\nu , (T_h)_{h\in\RR})$ be the noise generated by the canonical
flow associated with $\nu$.
Notice that $\nu\succeq \bar{\nu}$ means the existence of
$\bar{N}^0$ a sub-noise of $N_\nu$ such that $\P_{\bar{\nu}}$ is the
law of $\bar{K}^0$, the flow obtained by filtering the canonical flow
of law $\P_\nu$ with respect to $\bar{N}^0$.

Note that the mapping $K~:~(\O,\cA) \to (\O^0,\cA^0)$ is
measurable. Let $\bar{\cF}$ be the completion  of
$K^{-1}(\bar{\cF}^0)$ by all $\P$-negligible sets of $\cA$ and,
for all $s\leq t$, set $\bar{\cF}_{s,t} = \bar{\cF}\cap
{\cF}_{s,t}$. Then $\bar{N} = (\O,\cA,(\bar{\cF}_{s,t})_{s\leq
  t},\P,(T_h)_{h\in\RR})$  is a sub-noise of 
$N$. Lemma \ref{lemfiltr} allows us to define $\bar{K}$ the flow
obtained by filtering $K$ with respect to $\bar{N}$. One can check
that $\bar{K}=\bar{K}^0(K)$. This implies that the law of $\bar{K}$ is
$\P_{\bar{\nu}}$. Thus the proposition is proved. \qed

\bprop The domination relation and the weak domination relation are
partial orders on the class of Feller convolution semigroups. \eprop

\prf {\bf 1)} The transitivity of the domination relation follows from
lemma \ref{propfiltr} by the chain rule for conditional expectations.

Let us observe that if $\nu^1\preceq \nu^2$ and $\nu^2\preceq
\nu^1$ then $\nu^1=\nu^2$. Indeed, if $\nu^1\succeq \nu^2$, Jensen's
inequality shows that for all $f_1,\dots,f_n$ in $C(M)$,
$x_1,\dots,x_n$ in $M$  and $t\geq 0$,
\beq \E_{\nu^1}\left[\exp\left(\sum_{i=1}^n K_{0,t}f_i(x_i)\right)\right] \geq
\E_{\nu^2}\left[\exp\left(\sum_{i=1}^n K_{0,t}f_i(x_i)\right)\right]. \eeq
Therefore, if moreover $\nu^1\preceq\nu^2$, the preceding inequality
becomes an equality. This proves $\nu^1=\nu^2$. 

\medskip {\bf 2)} For the weak domination relation, the proof is
similar. We prove the transitivity using the product of
extensions. Indeed, if $\bar{\nu}\preceq^w\nu$, given any extension
$(N^1,K^1)$ of $N_\nu$, there exist a larger extension $(N,K)$ and a subnoise
$\bar{N}$ of $N$ such that $\bar{K}$ has law $\P_{\bar{\nu}}$: let
$\bar{N}^2$ be a sub-noise of an extension $(N^2,K^2)$ of $N_\nu$ such
that $\bar{K}^2$ has law $\P_{\bar{\nu}}$, then $(N,K)$ is taken as
the product of the extensions $(N^1,K^1)$  and $(N^2,K^2)$, and
$\bar{N}$ is induced by $\bar{N^2}$.   \qed

\brem  The concept of filtering will be used in sections \ref{coalsemi}, \ref{34}, \ref{classi}
and an example is given in the following section. \erem

\subsection{An example of filtering.}\label{exemple}
Let $M=\{0,1\}$. Then $F$, the set of maps from $\{0,1\}$ on
$\{0,1\}$ is constituted of the maps $\sigma$, $I$, $f_0$ and $f_1$,
with $I$ the identity, $\sigma(0)=1$, $\sigma(1)=0$, $f_0=0$ and
$f_1=1$. Let $(N_t)$ be a Poisson process on $\RR$ and
$(\p_n)_{n\in\ZZ}$ be a sequence, independent of the Poisson process,
of independent random variables taking their values in $F$ with law
$$\frac{1}{4}(\delta_{f_0}+\delta_{f_1}+\delta_{I}+\delta_{\sigma}).$$

We then define a stochastic flow of mappings on $\{0,1\}$ by
$$\left\{\begin{array}{lll} \p_{s,t}=I &\hbox{if}& N_t-N_s=0\\
\p_{s,t}=\p_{N_t-1}\circ\cdots\circ\p_{N_s} &\hbox{if}& N_t-N_s>0
\end{array}\right.$$
for all $s\leq t$. Note that $\p$ is a coalescing flow since for all
$s$, there is a.s. a finite time $T$ such that for all $t\geq T$,
$\p_{s,t}(0)=\p_{s,t}(1)$. The one-point motion of this flow is given
by the symmetric random walk with generator $A^{(1)}$ given by
$$A^{(1)} = \left(\begin{array}{ll} 1/2 & 1/2 \\ 1/2 & 1/2
\end{array}\right).$$
Note also that, since the $\{0,1\}$ has only two points, the $n$-point
motions associated with this stochastic flow of mappings are determined by
the two-point motion. The generator $A^{(2)}$ of the
two-point motion is (the state space is $\{(0,0),(1,1),(0,1),(1,0)\}$)
$$A^{(2)}=\left(\begin{array}{llll}
1/2 & 1/2 & 0 & 0\\
1/2 & 1/2 & 0 & 0\\
1/4 & 1/4 & 1/4 & 1/4\\
1/4 & 1/4 & 1/4 & 1/4 \end{array}\right)$$

With the stochastic flow $\p$ and an independent sequence of random
variables $(Y_n)_{n\in\ZZ}$ with $\P[Y_n=1]=p=1-\P[Y_n=0]$, we define
a stochastic flow of kernels $K$, by
$$\left\{\begin{array}{lll} K_{s,t}(i)=\delta_i &\hbox{if}& N_t-N_s=0\\
K_{s,t}=K_{N_s}\cdots K_{N_t-1} &\hbox{if}& N_t-N_s>0
\end{array}\right.$$
where $K_n=Y_n\delta_{\p_n}+(1-Y_n)\frac{1}{2}(\delta_0+\delta_1)$.

Denote by $N^c$ the noise of $\p$, by $N$ the noise of $K$ and by
$\hat{N}$ the noise of $(\p,Y)$. Then $N^c$ is the noise of
$(N_t,\p_{N_t})$, $\hat{N}$ is the noise of  $(N_t,\p_{N_t},Y_{N_t})$
and $N$ is the noise of $(N_t,K_{N_t})$.
The noises $N^c$ and $N$ are subnoises of $\hat{N}$. And $N$ cannot be
isomorphic to a subnoise of $N^c$. Indeed for $\eps$ small,
$\cF_{0,\eps}^{N^c}$ has one atom of probability $e^{-\eps}$ and four
atoms of probability $\frac{1}{4}\eps e^{-\eps}$, and
$\cF_{0,\eps}^{N}$ has one atom of probability $e^{-\eps}$ as well but
one atom of probability $(1-p)\eps e^{-\eps}$ and four atoms of
probability $\frac{p}{4}\eps e^{-\eps}$. 

The flow $K$ coincides with the flow obtained by filtering $\p$
with respect to $N$. Thus the law of $K$ is weakly dominated by the
law of $\p$ but is not dominated. 

\subsection{Continuous martingales.}
Let $(K_{s,t},~s\leq t)$ be a stochastic flow of kernels. For all
$s\leq t$ set ${\cal F}_{s,t}=\sigma(K_{u,v},~s\leq u\leq v\leq t)$.
Let ${\cal F}$ be the filtration $({\cal F}_{0,t})_{t\geq 0}$.
Let ${\cal M}({\cal F})$ be the space of locally square integrable
${\cal F}$-martingales.
\bprop\label{martcont} Suppose that $\P^{(1)}_t$ is the semigroup of a
Markov process with continuous paths, then all martingales of 
${\cal M}({\cal F})$ are continuous. \eprop
\prf Let $M\in {\cal M}({\cal F})$ be a martingale in the form
$\E[F|{\cal F}_{0,t}]$ where $F=\prod_{i=1}^n K_{s_i,t_i}f_i(x_i)$,
with $f_1,\dots,f_n$ in $C(M)$, $x_1,\dots,x_n$ in $M$ and
$0\leq s_i<t_i$ (we take here the continuous modification in $t$ of the
stochastic flow of kernels). By definition of the filtration,
functions in this form are dense in $L^2(\cF_{0,\infty})$. This
implies that martingales of this form are dense in ${\cal M}({\cal
F})$. Since the space of continuous martingales is closed in ${\cal M}({\cal
F})$, it is enough to prove the continuity of these
martingales.

\smallskip
For all $t$, let $\tilde{K}_t$ be the kernel defined on $\RR^+\times
M$ by
\beq \tilde{K}_t(s,x) = \left\{
\begin{array}{lll} \delta_{s-t}\otimes\delta_x &\hbox{for} & s\geq t,\\
\delta_{0}\otimes K_{s,t}(x) & \hbox{for} & s\leq t. 
\end{array}\right. \eeq
Then we can rewrite $F$ in the form $\prod_{i=1}^n
\tilde{K}_{t_i}\tilde{f}_i(s_i,x_i)$, where $\tilde{f}_i(s,x)=f_i(x)$.

Note that $(\tilde{K}_{t_i}(s_i,x_i),~1\leq i\leq n)$ is a Markov
process on $({\cal B}(\RR^+)\otimes{\cal P}(M))^n$. This Markov process
is continuous and Feller (the Feller property follows from the Feller
property of the semigroups $(\Pi^{(k)}_t,~k\geq 1)$). It is well known
that the martingales relative to the filtration denoted here
$({\cal F}^{\{s_i,x_i\}_{1\leq i\leq n}}_t,~t\geq 0)$ generated by
such a process are continuous (see \cite{R-W} tome II).

This proves that $\E[F|{\cal F}^{\{s_i,x_i\}_{1\leq i\leq n}}_t]$ is a
continuous martingale. We conclude after remarking that 
$M_t=\E[F|{\cal F}^{\{s_i,x_i\}_{1\leq i\leq n}}_t]$, which holds
since the $\sigma$-field ${\cal F}^{\{s_i,x_i\}_{1\leq i\leq n}}_t$
is a sub-$\sigma$-field of ${\cal F}_t$ and $M_t$ is easily seen to be
${\cal F}^{\{s_i,x_i\}_{1\leq i\leq n}}_t$-measurable. \qed

\setcounter{equation}{0}
\section{Stochastic coalescing flows.}\label{coal}
In this section we study stochastic coalescing flows, we denote by
$(\p_{s,t},~s\leq t)$. In section \ref{ptmeas}, it is shown that for
all $s<t$, $\p_{s,t}^*(\lambda)$ is atomic (where $\lambda$ denotes
any positive Radon measure on $M$). We study this point measure valued
process which provides a description of the coalescing flow.

In section \ref{coalsemi}, starting from a compatible family of Feller
semigroups, under the hypothesis that starting close to the
diagonal the two-point motion hits the diagonal with a probability
close to 1, we construct another compatible family of Feller
semigroups to which is associated a stochastic coalescing
flow. We then show that the stochastic flow of kernels associated with
the first family of semigroups can be defined by filtering the
stochastic coalescing flow with respect to a sub-noise of an extension
of its canonical noise.

Finally, we give three examples. The first one, due to Arratia
\cite{Arratia}, describes the flow of independent Brownian motions
sticking together when they meet. The second one is due to Propp and
Wilson \cite{P-W} in the context of perfect simulation of the
invariant distribution of a finite state irreducible Markov chain,
their stochastic flows being indexed by the integers. The third one is
the construction of a stochastic coalescing flow solution of Tanaka's SDE
\beq dX_t=\hbox{sgn}(X_t)dW_t, \eeq
where $W$ is a real white noise. This coalescing flow was
constructed by Watanabe in \cite{Watanabe} and Warren in
\cite{warren}. In \cite{ljr}, a stochastic flow of kernels solution
of this SDE was constructed as the only strong solution of this SDE.
\subsection{Definition.}
Let $M$ be a locally compact separable metric space.
\begin{definition} A stochastic flow of mappings on $M$,
$(\p_{s,t},~s\leq t)$, is called a stochastic coalescing flow if
for some $(x,y)\in M^2$, $T_{x,y}=\inf\{t\geq 0,~\p_{0,t}(x)=\p_{0,t}(y)\}$
is finite with a positive probability and for all $t\geq T_{x,y}$,
$\p_{0,t}(x)=\p_{0,t}(y)$.  In other words, a pair of points stick
together after a finite time with a positive probability. 
\end{definition}

\brem {\em This definition depends only on the two-point motion.}\erem

Let $(\P^{(n)}_t,~n\geq 1)$ be a compatible family of Feller
semigroups. We denote by $\P^{(2)}_{(x,y)}$ the law of the Markov
process associated with $\P^{(2)}_t$ starting from $(x,y)$ and we denote
this process $(X_t,Y_t)$ or $X^{(2)}_t$. Let $T_\Delta=\inf\{t\geq
0,~X_t=Y_t\}$.

\brem  A compatible family $(\P^{(n)}_t,~n\geq 1)$ of Feller
semigroups  defines a stochastic coalescing flow if and only if
for all $(x,y)\in M^2$, for all $t\geq T_\Delta$, 
$X_t=Y_t$, $\P^{(2)}_{(x,y)}$-almost surely, and for some $(x,y)\in M^2$, 
$\P^{(2)}_{(x,y)}[T_{\Delta}<\infty]>0$. \erem

\subsection{A point measure valued process associated with a stochastic
coalescing flow.}\label{ptmeas}

In this subsection, we suppose we are given a compatible family of Feller
semigroups $(\P^{(n)}_t,~n\geq 1)$ such that
\beq \label{condcoal} 
\left\{\begin{array}{l}\forall x\in M,~\forall t>0,\quad \lim_{y\to x}
\P^{(2)}_{(x,y)}[X_t\neq Y_t] = 0,\\
\forall (x,y)\in M^2,\qquad \P^{(2)}_{(x,y)}[T_{\Delta}<\infty] > 0.\end{array}\right. \eeq
\brem  The assumption (\ref{condcoal}) implies that the associated
stochastic flow is a stochastic coalescing flow and is verified in all
the examples of coalescing flows we will study except for the example
presented in section \ref{sgn}, where $\P^{(2)}_{(x,y)}[X_t\neq Y_t]$
does not converge towards 0 as $y$ tends to $x$ when $x\neq 0$.\erem

Let $\p=(\p_{s,t},~s\leq t)$ be a measurable stochastic coalescing
flow associated with $(\P^{(n)}_t,~n\geq 1)$. For all $s<t\in \RR$,
let $\mu_{s,t}=\p^*_{s,t}(\lambda)$, where $\lambda$ is any positive
Radon measure on $M$.

\bprop \bdes \ita For all $s<t\in\RR$, almost surely, $\mu_{s,t}$ is
atomic.
\itb For all $s<u<t\in\RR$, almost surely, $\mu_{s,t}$ is absolutely
continuous with respect to $\mu_{u,t}$. \edes
\eprop
\prf
Fix $s<t\in\RR$. For all positive $\eps$ and all $x\in M$, let
$m_\eps^x=\int_{B(x,\eps)}1_{\p_{s,t}(x)=\p_{s,t}(y)}~\lambda(dy)$ ($m_\eps^x$ is
well defined since $(x,\o)\mapsto\p_{s,t}(x,\o)$ is measurable). For all
$\a\in ]0,1[$ and $x\in M$, let
\beq A_n^{\a,x}=\{m_{\eps^x_n}^x<(1-\a)\lambda(B(x,\eps^x_n))\},\eeq
where $\eps^x_n$ is a positive sequence such that $d(x,y)\leq\eps_n^x$ implies
$$\P^{(2)}_{(x,y)}[X_{t-s}\neq Y_{t-s}]\leq 2^{-n}.$$

\blem For all positive $\a$, $x\in M$ and $n\in \NN$,
$\P(A_n^{\a,x})\leq\frac{1}{\a 2^n}$.\elem 
\prf For all integer $n$, we have
\beqarr
\E[m_{\eps^x_n}^x]&=&\int_{B(x,\eps^x_n)}\P^{(2)}_{(x,y)}[X_{t-s}=Y_{t-s}]~\lambda(dy)\\
&=&\int_{B(x,\eps^x_n)}(1-\P^{(2)}_{(x,y)}[X_{t-s}\neq Y_{t-s}])~\lambda(dy)\\
&\geq&(1-2^{-n})\lambda(B(x,\eps^x_n)).
\eeqarr
And we conclude since
$$\E[m_{\eps^x_n}^x]\leq
\P(A_n^{\a,x})(1-\a)\lambda(B(x,\eps^x_n)) +
(1-\P(A_n^{\a,x}))\lambda(B(x,\eps^x_n))\qquad$$ 
(we use the fact that $m_{\eps^x_n}^x\leq\lambda(B(x,\eps^x_n))$). \qed

\blem \label{equiv}For all $x\in M$, almost surely,
$m_{\eps^x_n}^x\sim\lambda(B(x,\eps^x_n))$ as $n\to\infty$.\elem
\prf Using Borel-Cantelli's lemma, for all $\a\in ]0,1[$
$$1 - \a\leq\liminf_{n\to\infty}
\frac{m_{\eps^x_n}^x}{\lambda(B(x,\eps^x_n))} \leq \limsup_{n\to\infty}
\frac{m_{\eps^x_n}^x}{\lambda(B(x,\eps^x_n))} \leq 1$$
almost surely. This implies
$\lim_{n\to\infty}\frac{m_{\eps^x_n}^x}{\lambda(B(x,\eps^x_n))}=1$ a.s. \qed

\medskip
Since for all $(x,\o)\in M\times \Omega$,
\beqarr\mu_{s,t}(\{\p_{s,t}(x)\})&=&\lambda(\{y,~\p_{s,t}(y)=\p_{s,t}(x)\})\\
&\geq&\lambda(\{y\in B(x,\eps_n),~\p_{s,t}(y)=\p_{s,t}(x)\}),\eeqarr
lemma \ref{equiv} implies that for all $x\in M$, 
\beq\label{as}\mu_{s,t}(\{\p_{s,t}(x)\})>0\eeq
almost surely. Since $(x,\o)\mapsto
\mu_{s,t}(\{\p_{s,t}(x)\})$ is measurable,
\beq \lambda(dx)\otimes\P(d\o)\hbox{-a.e.},\qquad \mu_{s,t}(\{\p_{s,t}(x)\})>0.\eeq
This equation implies (since $\mu_{s,t}=\p_{s,t}^*(\lambda)$)
\beq \mu_{s,t}(dy)\hbox{-a.e.},\qquad \mu_{s,t}(\{y\})>0\eeq
almost surely. This last equation is one characterization of the atomic
nature of $\mu_{s,t}$ and {\bf (a)} is proved.

\medskip
To prove {\bf (b)}, note first that $\lambda(dx)\otimes \P(d\o)$-a.e.,
$\p_{u,t}^*(\delta_x)=\delta_{\p_{u,t}(x)}$ is absolutely continuous
with respect to $\p_{u,t}^*(\lambda)=\mu_{u,t}$ since (\ref{as})
holds. Note also that $\lambda(dx)\otimes\P(d\o)$-a.e.,
$\p_{s,t}(x)=\p_{u,t}\circ \p_{s,u}(x)$. This implies
\beq \mu_{s,t}=\p_{u,t}^*(\mu_{s,u})\qquad \hbox{a.s.}\eeq
Since $\mu_{s,u}$ is atomic, independent of $\p_{u,t}$ and
$\E[\mu_{s,u}]=\lambda$, it follows that
$\mu_{s,t}$ is absolutely continuous with respect to $\mu_{u,t}$. This
proves {\bf (b)}. \qed

\brem  {\em $\bullet$ $(\mu_{s,t},~s\leq t)$ is Markovian in $t$.

\medskip
$\bullet$ Since $\mu_{s,t}$ is atomic for $t>s$, there exists a point process
$\xi_{s,t}=\{\xi_{s,t}^i\}$ and weights $\{\a^i_{s,t}\}\in\RR^\NN$
such that  $\mu_{s,t}=\sum_i\a^i_{s,t}\delta_{\xi^i_{s,t}}$. The point
process $\xi_{s,t}$ and the marked point process
$(\xi_{s,t},\alpha_{s,t})$ are Markovian in $t$ since for all $s<u<t$,
$\xi_{s,t}=\p_{u,t}(\xi_{s,u})$ and
$\a^i_{s,t}=\sum_{\{j,~\xi^i_{s,t}=\p_{u,t}(\xi^j_{s,u})\}}\a^j_{s,u}$.

\medskip
$\bullet$ Let $A^i_{s,t}=\p_{s,t}^{-1}(\xi^i_{s,t})$ and $\Pi_{s,t}$ be the
collection of the sets $A^i_{s,t}$. Note that $\cup_i A^i_{s,t}=M$
$\lambda$-a.e, the union being disjoint. Note also that 
$\xi_{s,t}$ and $\Pi_{s,t}$ determines $\p_{s,t}$ $\lambda$-a.e. Note
finally that $\Pi_{s,t}$ is Markovian in $s$ when $s$ decreases, since
for all $s<u<t$, $\Pi_{s,t}=\{\p_{s,u}^{-1}(A^{i}_{u,t})\}$. This
Markov process has also a coalescence property~: one can have for
$i\neq j$, $\p_{s,u}^{-1}(A^{i}_{u,t})=\p_{s,u}^{-1}(A^{j}_{u,t})$. 
When $s$ decreases, the partition $\Pi_{s,t}$ becomes coarser.} 
\erem

\subsection{Construction of a family of coalescent semigroups.}\label{coalsemi}
Let $(\P^{(n)}_t,~n\geq 1)$ be a compatible family of Feller 
semigroups on a separable locally compact metric space $M$ and
$\nu=(\nu_t)_{t\in\RR}$ the associated Feller convolution semigroup on
$(E,{\cal E})$. Let $\Delta_n=\{x\in M^n,~\exists i\neq j, ~x_i=x_j\}$
and $T_{\Delta_n}=\inf\{t\geq 0,~X^{(n)}_t\in\Delta_n\}$, where
$X^{(n)}_t$ denotes the $n$-point motion, i.e. the Markov process on
$M^n$ associated with the semigroup $\P^{(n)}_t$. We will denote
$\Delta_2$ by $\Delta$. 

\bthm\label{sgcoal} There exists a unique compatible family
$(\P^{(n),c}_t,~n\geq 1)$ of Markovian semigroups on $M$ such that if
$X^{(n),c}$ is the associated $n$-point motion and
$T^c_{\Delta_n}=\inf\{t\geq 0,~X^{(n),c}_t\in\Delta_n\}$, then 
\bit
\item $(X^{(n),c}_t,~t\leq T^c_{\Delta_n})$ is equal in law to
$(X^{(n)}_t,~t\leq T_{\Delta_n})$,
\item for $t\geq T^c_{\Delta_n}$, $X^{(n),c}_t\in \Delta_n$.\eit
Moreover, this family is constituted of Feller semigroups if condition
{\bf (C)} below is satisfied, 
\bdes\item[(C)] For all $t>0$, $\eps>0$ and $x\in M$,
$$\lim_{y\to x}
\P^{(2)}_{(x,y)}[\{T_\Delta>t\}\cap\{d(X_t,Y_t)>\eps\}]=0$$
where $(X_t,Y_t)=X^{(2)}_t$. And for some $x$ and $y$ in $M$,
$\P^{(2)}_{(x,y)}[T_\Delta<\infty]>0$. \edes
In this case, $(\P^{(n),c}_t,~n\geq 1)$ satisfies (\ref{mapcond}) and
is associated with a coalescing flow. \ethm 

\prf For all $n\geq 1$, let $\cP_n$ be the set of all partitions of
$\{1,\dots,n\}$. The number of elements of $\pi\in\cP_n$ is denoted
$|\pi|$. For all $\pi\in\cP_n$, we define the equivalent relation
$i\sim_\pi j$ if $i$ and $j$ belong to the same element of $\pi$. We
define a partial order on $\cP_n$ by $\pi'\leq\pi$ if $i\sim_{\pi'}j$
implies $i\sim_\pi j$ ($\pi$ is finer that $\pi'$).

For all $\pi\in\cP_n$, we let $E_\pi$ be the set of elements $x\in
M^n$ such that $x_i=x_j$ if $i\sim_\pi j$ and
$\partial E_\pi=\cup_{\pi'<\pi} E_{\pi'}$, the set of elements
$x\in E_\pi$ such that there exists $i$ and $j$ with $i\not\sim_\pi j$
and $x_i=x_j$. Let $j_\pi$ be an isometry between $M^{|\pi|}$ and
$E_\pi$.

By induction on $k=|\pi|$, we define a Markov process $X^{\pi}$ on $E_\pi$.
For $k=1$, we let $X^\pi=j_\pi(X^{(1)})$. Assume now we have defined a
Markov process on $E_\pi$ for all $\pi$ such that $|\pi|\leq k$. Let
$\pi\in\cP_n$ with $|\pi|=k+1$, we define $X^\pi$ concatenating the
process $j_\pi(X^{(k+1)})$ stopped at the entrance time $T$
in $\partial E_\pi$ with the process $X^{\pi'}$ starting from the
corresponding point and where $\pi'$ is the finest partition such
that $X_T^{\pi'}\in E_{\pi'}$. This way, we construct a Markov process
on $M^n$, $X^{(n),c}=X^{\pi}$ for $\pi=\{\{1\},\dots,\{n\}\}$.

\medskip
For all integer $n$, let $\P^{(n),c}_t$ be the Markovian semigroup
associated with the Markov process $X^{(n),c}$. From the above
construction, it is clear that the family $(\P^{(n),c}_t,~n\geq 1)$
of Markovian semigroups is compatible. 

\medskip It remains to prove that when {\bf (C)} is satisfied, this
family of Markovian semigroups is constituted of Feller
semigroups. This holds since {\bf (C)} implies {\bf (F)} in lemma
\ref{Fel2n}~: for all positive $\eps$,
$\P^{(2),c}_{(x,y)}[d(X_t,Y_t)>\eps] \leq
\P^{(2)}_{(x,y)}[\{T_\Delta>t\}\cap \{d(X_t,Y_t)>\eps\}]$ which
converges towards 0 as $y\to x$. Note that when {\bf (C)} holds, it is
easy to see that the canonical flow is a coalescing flow. \qed

\medskip
We now suppose that $(\P^{(n),c}_t,~n\geq 1)$ is constituted of Feller
semigroups (which is true when {\bf (C)} holds). We denote by $\nu^c$
the associated Feller convolution semigroup.

\bthm \label{filtrage} The convolution semigroup $\nu^c$ weakly dominates
$\nu$. \ethm

\prf The idea of the proof is to define a coupling between the flows
of kernels $K$ and $K^c$ respectively of law $\P_\nu$ and
$\P_{\nu^c}$. (Since we did not assume {\bf (C)} holds, it is not
clear that $K^c$ is a flow of mappings.)

In a way similar to the construction of the Markov process
$X^{(n),c}$ in the proof of theorem \ref{sgcoal}, for all integer
$n\geq 1$, we construct a Markov process $\hat{X}^{(n)}$ on $(M\times
M)^n$ identified with $M^n\times M^n$ such that:
\begin{itemize}
\item $(\hat{X}^{(n)}_1,\dots,\hat{X}^{(n)}_n)$ is the
$n$-point motion of $\nu^c$,
\item $(\hat{X}^{(n)}_{n+1},\dots,\hat{X}^{(n)}_{2n})$
is the $n$-point motion of $\nu$,
\item  between the coalescing times, $\hat{X}^{(n)}$ is described by the
$(k+n)$-point motion of $\nu$ (when
$(\hat{X}^{(n)}_1,\dots,\hat{X}^{(n)}_n)$ belongs to $E_\pi$, with
$|\pi|=k$).
\end{itemize}

Let $\hat{\P}^{(n)}_t$ be the Markovian semigroup associated with
$\hat{X}^{(n)}$. One easily checks that this semigroup is Feller using
the fact that $\P^{(n)}_t$ and $\P^{(n),c}_t$ are Feller. Then
$(\hat{\P}^{(n)}_t,~n\geq 1)$ is a compatible family of Feller
semigroups, associated with a Feller convolution semigroup $\hat{\nu}$.

Let $\hat{K}$ be the canonical stochastic flow associated with this
family of semigroups. Straightforward computations show that for all
$s<t$, $(f,g)\in C(M)^2$ and $(x,y)\in M^2$,
\beqarr
\E[(\hat{K}_{s,t}(f\otimes g)(x,y))^2] &=& \P^{(3)}_{t-s}f^2\otimes
g\otimes g(x,y,y),\\
\E[(\hat{K}_{s,t}(f\otimes 1)\hat{K}_{s,t}(1\otimes g))^2(x,y)] &=&
\P^{(3)}_{t-s}f^2\otimes g\otimes g(x,y,y),\\
\E[(\hat{K}_{s,t}(f\otimes g)\hat{K}_{s,t}(f\otimes
1)\hat{K}_{s,t}(1\otimes g))(x,y)] &=& \P^{(3)}_{t-s}f^2\otimes
g\otimes g(x,y,y). \eeqarr
This implies that
\beq \E[(\hat{K}_{s,t}(f\otimes g)-\hat{K}_{s,t}(f\otimes
1)\hat{K}_{s,t}(1\otimes g))^2(x,y)] = 0. \eeq
Thus we have $\hat{K}_{s,t}(x,y)=K^c_{s,t}(x)\otimes K_{s,t}(y)$ and
it is easy to check that the laws of $K^c$ and of $K$ are respectively
$\P_{\nu^c}$ and $\P_\nu$. Thus $(N_{\hat{\nu}},K^c)$ is an extension
of $N_{\nu^c}$.  Let $\tilde{N}_\nu$ be the sub-noise of
$N_{\hat{\nu}}$ generated by $K$. 

Let us notice now that for all $g,f_1,\dots,f_n$ in $C_0(M)$,
all $y,x_1,\dots,x_n$ in $M$ and all $s<t$, we have (setting
$y_i=x_{n+1}=y$ and for $i\leq n$, $h_i=f_i\otimes 1$ and
$h_{n+1}=1\otimes g$)
\beqarr
\E\left[K^c_{s,t}g(y)\prod_{i=1}^nK_{s,t}f_i(x_i)\right] &=&
\E\left[\prod_{i=1}^{n+1}\hat{K}_{s,t}h_i(x_i,y_i)\right]\\
&=& \P^{(n+1)}_{t-s} f_1\otimes\cdots\otimes f_n\otimes
g(x_1,\cdots,x_n,y). \eeqarr
More generally one can prove in a similar way that for all
$g,f_1,\dots,f_n$ in $C_0(M)$, all $y,x_1,\dots,x_n$ in $M$, all $s<t$
and all $(s_i,t_i)_{1\leq i\leq n}$ with $s_i\leq t_i$ that
\beq \E\left[K^c_{s,t}g(y)\prod_{i=1}^nK_{s_i,t_i}f_i(x_i)\right] =
\E\left[K_{s,t}g(y)\prod_{i=1}^nK_{s_i,t_i}f_i(x_i)\right].\eeq
This implies that $K_{s,t}g(y)=\E[K^c_{s,t}g(y)|\sigma(K)]$ and
therefore that $\nu^c\succeq^w\nu$. \qed

\brem  \label{remfel}
{\em Let $(X^{(n)},~n\geq 1)$ be a family of strong Markov
processes respectively taking their values in $M^n$. We suppose that the
associated family of Markovian semigroups $(\P^{(n)}_t,~n\geq 1)$ is
compatible and that for all $x\in M$,
\beq \lim_{y\to x}\P^{(2)}_{(x,y)}[\{T_\Delta>t\}\cap
\{d(X_t,Y_t)>\eps\}]=0 \eeq
for all $\eps>0$ and $t>0$.  Then $(\P^{(n)}_t,~n\geq 1)$ (and
$(\P^{(n),c}_t,~n\geq 1)$) are Feller semigroups. 

One can prove this with a coupling similar to the coupling given in
the proof of the previous theorem~: the idea is
to construct on the same probability space two Markov processes
$X^{(n)}$ and $Y^{(n)}$ associated to $\P^{(n)}_t$ and such that
$X^{(n)}_i(t)=Y^{(n)}_i(t)$ if $t\geq
\inf\{s,~X^{(n)}_i(s)=Y^{(n)}_i(s)\}$. }\erem

\brem  {\em The example given in section \ref{exemple} gives an
illustration of the two theorems of this section, first with
$\P^{(n)}_t=\P_t^{\otimes n}$ then with $\P^{(n)}_t$ the semigroup of the
$n$-point motion of $K$. This example shows in particular that one can
have $\nu\preceq^w\nu^c$ and $\nu\not\preceq \nu^c$.} \erem

\subsection{Examples.}
\subsubsection{Arratia's coalescing flow of independent Brownian motions.}
The first example of coalescing flows was given by Arratia
\cite{Arratia}. On $\RR$, let $\P_t$ be the semigroup of a Brownian
motion. With this semigroup we define the compatible family
$(\P^{\otimes n}_t,~n\geq 1)$ of Feller semigroups. 
Note that the $n$-point motion of this family of semigroups is given
by $n$ independent Brownian motions. Let us also remark that the
canonical stochastic flow of kernels associated with this family of
semigroups is not random and is given by $(\P_{t-s},~s\leq t)$. 

\smallskip
Let $(\P^{(n)}_t,~n\geq 1)$ be the compatible family
of Markovian coalescent semigroups associated with $(\P^{\otimes
n}_t,~n\geq 1)$ (see section \ref{coalsemi}). Note that the $n$-point
motion of this family of semigroups is given by $n$ independent
Brownian motions who stick together when they meet. 

\bprop The family $(\P^{(n)}_t,~n\geq 1)$ is constituted of Feller
semigroups and is associated with a coalescing flow. \eprop 
\prf It is obvious after remarking that two real independent
Brownian motions meet each other almost surely (condition {\bf (C)} is
verified). \qed

\subsubsection{Propp-Wilson algorithm.} 
Similarly to Arratia's coalescing flow, let $\P_t$ be the semigroup of
an irreducible aperiodic Markov process on a finite set $M$, with
invariant probability measure $m$. Let $(\P^{(n)}_t,~n\geq 1)$ be the
compatible family of Markovian coalescent semigroups associated with 
$(\P^{\otimes n}_t,~n\geq 1)$. The coalescing flow in section
\ref{exemple} is of this type. 

\bprop The family $(\P^{(n)}_t,~n\geq 1)$ is constituted of Feller
semigroups and is associated with a coalescing flow. \eprop 
\prf It is obvious since the two-point motion defined by $\P^{\otimes
2}_t$ hits the diagonal almost surely. \qed

\medskip
Let $\p=(\p_{s,t},~s\leq t)$ denote this coalescing flow. Then almost
surely, for all $x,y$ in $M$,
$\tau_{x,y}=\inf\{t>0,~\p_{0,t}(x)=\p_{0,t}(y)\}$ is
finite. Therefore, after a finite time $\hbox{Card}\{\p_{0,t}(x),~x\in
M\}=1$.

\medskip
In Propp-Wilson \cite{P-W}, an algorithm to exactly simulate a random
variable distributed according to the invariant probability measure of
a Markov chain with finite state space is given. The method consists in
constructing a stochastic coalescing flow. We explain this in our
context.

Let $\tau=\inf\{t>0,~\p_{-t,0}(x)=\p_{-t,0}(y) \hbox{ for all }
(x,y)\in M^2\}$. 
\bprop $\tau$ is almost surely finite and the law of $X_\tau$, the random
variable $\p_{-\tau,0}(x)$ (independent of $x\in M$), is $m$. \eprop
\prf Let us remark that for $t>\tau$ and all $x\in M$, the cocycle property
implies that $\p_{-t,0}(x)=X_\tau$. 

Since for all positive $t$, 
\beq\begin{array}{lll}
\P[\tau\geq t] &=& \P[\exists x,y,~\p_{-t,0}(x)\neq \p_{-t,0}(y)]\\
&\leq& \sum_{(x,y)\in M^2}\P[\tau_{x,y}\geq t]
\end{array} \eeq
which converges towards 0 as $t$ goes to
$\infty$. Thus $\tau<\infty$ a.s.

\medskip
For all function $f$ on $M$ and all $x\in M$,
$\lim_{t\to\infty}\P_tf(x)=\sum_{y\in M} f(y)m(y)$ and
\beq \P_tf(x)=\E[f(\p_{-t,0}(x))]=\E[f(\p_{-t,0}(x))1_{t\leq \tau}] 
+ \E[f(X_\tau)1_{\tau <t}].\eeq
Since $\tau$ is almost surely finite, as $t$ goes to $\infty$,
the first term of the right hand side of the preceding equation
converges towards 0 and the second term converges towards
$\E[f(X_\tau)]$. Therefore we prove that $\E[f(X_\tau)]=\sum_{y\in M}
f(y)m(y)$. \qed 

\subsubsection{Tanaka's SDE.}\label{sgn}
In \cite{ljr}, starting from a real Brownian motion $B$, we
constructed a family of random operators $(S_t,~t\geq 0)$, strong
solution of the SDE
\beq dX_t=\hbox{sgn}(X_t)dB_t,\quad t\geq 0. \label{sgnSDE}\eeq
For $f$ continuous,
\beq S_tf(x) = f(R_t^x) 1_{t<T_x} + \frac{1}{2}(f(R_t^x)+f(-R_t^x))
1_{t\geq T_x},\eeq 
where $R_t^x$ is the Brownian motion $x+B_t$ reflected at 0 and $T_x$
the first time it hits 0.
For all continuous functions $f_1,\dots,f_n$, let 
\beq \P^{(n)}_t(f_1\otimes\cdots\otimes f_n)(x_1,\dots,x_n) =
\E\left[\prod_{i=1}^nS_tf_i(x_i)\right].\eeq
Then it is easy to see that $(\P^{(n)}_t,~n\geq 1)$ is a compatible
family of Feller semigroups. Let $(\P^{(n),c}_t,~n\geq 1)$
be the family of semigroups constructed by theorem \ref{sgcoal}.

\medskip Let us describe the $n$-point motion associated with
$(\P^{(n),c}_t,~n\geq 1)$. Let $(X_t,~t\geq 0)$ be a Brownian motion
starting at 0. Let $B_t=\int_0^t \hbox{sgn}(X_s)~dX_s$, 
$(B_t,~t\geq 0)$ is also a Brownian motion starting at 0. For all
$x\in\RR$, let $\tau_x=\inf\{t\geq 0,~|x|+B_t=0\}$. Note that
$X_{\tau_x}=0$. For all $x\in\RR$, let
\beq X_t^x = \left\{\begin{array}{lll}
x + \hbox{sgn}(x) B_t &\hbox{ if }& t<\tau_x,\\ 
X_t &\hbox{ if }& t\geq \tau_x. \end{array}
\right.\eeq

Then $B_t = \int_0^t \hbox{sgn}(X_t^x)~dX_t^x$ and $X^x$ is a solution
of the SDE
\beq dX_t^x = \hbox{sgn}(X_t^x)~dB_t,\quad t\geq 0,\quad X_0^x=x. \eeq
Thus, for all $x_1,\dots,x_n$ in $M$,
$((X_t^{x_1},\dots, X_t^{x_n}),~t\geq 0)$ is the $n$-point motion of
the family of semigroups $(\P^{(n),c}_t,~n\geq 1)$.

\bprop The family $(\P^{(n)}_t,~n\geq 1)$ is constituted of Feller
semigroups and is associated with a coalescing flow. \eprop
\prf It is easy to see that $(\P^{(n),c}_t,~n\geq 1)$ is constituted of Feller
semigroups since for all $t$ and $x_0$, $x\mapsto X_t^x$ is
a.s. continuous at $x_0$ (it
implies that {\bf (F)} in lemma \ref{Fel2n} is satisfied). This also
implies that (\ref{mapcond}) is satisfied. Thus, the associated
stochastic flow is a flow of mappings. And it is a coalescing flow
since almost surely, every pair of point meets after a finite
time. Note that condition {\bf (C)} is verified. \qed 

\setcounter{equation}{0}
\section{Stochastic flows of kernels and SDEs.} \label{secW} 
\subsection{Hypotheses.}\label{hypo}
In this section, $M$ is a smooth locally compact manifold and we
suppose we are given $(\P^{(n)}_t,~n\geq 1)$, a compatible family of
Feller semigroups, or equivalently a Feller convolution
semigroup $\nu=(\nu_t)_{t\geq 0}$ on $(E,{\cal E})$. For all positive
integer $n$, we will denote by $X^{(n)}_t$ the $n$-point motion,
i.e. the Markov process associated with the semigroup $\P^{(n)}_t$.
We denote by $A^{(n)}$ the infinitesimal generator
of $\P^{(n)}_t$ and by $\cD(A^{(n)})$ its domain \footnote{$f$ is in
  the domain of the infinitesimal generator $A$ of a Feller  semigroup
  $\P_t$ if and only if $\frac{\P_tf-f}{t}$ converges uniformly as $t$
  goes towards 0. Its limit is denoted $Af$.}. We assume that 
\bdes \iti The space $C_K^2(M)\otimes C_K^2(M)$ \footnote{$C_K(M)$
(respectively $C^2_K(M)$) denotes the set of continuous (respectively
$C^2$) functions with compact support.} of functions of the form
$f(x)g(y)$, with $f,g$ in $C_K^2(M)$ and $x,y$ in $M$, is included in
${\cal D}(A^{(2)})$.
\itii The one-point motion $X^{(1)}_t$ has continuous paths.
\edes
In that case, we say that $\nu$ is a {\em diffusion convolution semigroup}
on $(E,{\cal E})$ and that the $\P^{(n)}_t$ are diffusion semigroups.

\subsection{Local characteristics of a diffusion convolution semigroup.}
Let us denote by $A$ the restriction of $A^{(1)}$ to $C^2_K(M)$. Note
that it follows easily from {\bf (i)} and {\bf (ii)} that for all
$f\in C_K^2(M)$,
\beq M^f_t=f(X^{(1)}_t)-f(X^{(1)}_0)-\int_0^tAf(X^{(1)}_s)~ds\eeq
is a martingale. Since $f^2$ also belongs to $C_K^2(M)$, using the
martingale $M^{f^2}$, it is easy to see that
\beq \langle M^f\rangle_t =\int_0^t\Gamma(f)(X^{(1)}_s)~ds \eeq
where
\beq\Gamma(f) = Af^2-2fAf.\eeq
In the following $\Gamma(f,g)$ will denote $A(fg)-fAg-gAf$, for $f$
and $g$ in $C^2_K(M)$.

\blem On a smooth local chart on an open set $U\subset M$, there exist
continuous functions on $U$, $a^{i,j}$ and $b^i$ such that for all
$f\in C_K^2(M)$, 
\beq\label{localchart} Af = \frac{1}{2}a^{i,j}\frac{\partial^2 f}{\partial x^i\partial x^j} +
b^i\frac{\partial f}{\partial x^i}. \eeq \elem
\prf For all $x\in U$, let $\p^i(x)=x^i$ denote the coordinate
functions of the local chart. We can extend $\p^i$ into an element of $C_K^2(M)$. For $f\in
C_K^2(M)$, using It\^o's formula, for $t<T_U$, the exit time of $U$,
$$ f(X^{(1)}_t) - f(X^{(1)}_0) -
\int_0^t\left(\frac{1}{2}a^{i,j}(X^{(1)}_s)\frac{\partial^2 f}{\partial
x^i\partial x^j}(X^{(1)}_s) +  b^i(X^{(1)}_s)\frac{\partial
f}{\partial x^i}(X^{(1)}_s)\right)~ds $$
is a martingale, where $b^i(x)=A\p^i(x)$ and
$a^{i,j}(x)=\Gamma(\p^i,\p^j)(x)$. And we get (\ref{localchart}) since
for all $x\in U$, $Af(x)=\lim_{t\to 0}\frac{\P^{(1)}_tf(x)-f(x)}{t}$. \qed

\medskip Note that the two-point motion $X^{(2)}_t$ has also
continuous trajectories and these results also apply to functions in
$C_K^2(M)\otimes C_K^2(M)$. For all $f,g$ in $C_K^2(M)$, let
\beq \label{defC} C(f,g) = A^{(2)}(f\otimes g) -f\otimes Ag - 
Af\otimes g. \eeq
It is clear that on a local chart on $U\times V\subset M\times M$,
\beq C(f,g)(x,y) = c^{i,j}(x,y)\frac{\partial f}{\partial x^i}(x)
\frac{\partial g}{\partial y^j}(y),\eeq
where $c^{i,j}\in C(U\times V)$. Then we can shortly write
$A^{(2)}=A\otimes I + I\otimes A + C$. On a local chart on $U\times
V$, for all $h\in C_K^2(M)\otimes C_K^2(M)$,
\beqar A^{(2)}h(x,y) &=& \frac{1}{2}a^{i,j}(x)\frac{\partial^2}{\partial
x^i\partial x^j}h(x,y) + b^i(x)\frac{\partial}{\partial x^i}h(x,y)\non\\
&+& \frac{1}{2}a^{i,j}(y)\frac{\partial^2}{\partial y^i\partial y^j}h(x,y) +
b^i(y)\frac{\partial}{\partial y^i}h(x,y)\label{gene2}\\ 
&+& C^{i,j}(x,y)\frac{\partial^2}{\partial x^i\partial y^j}h(x,y).\non
\eeqar

We will call $\Gamma(f,g)(x)-C(f,g)(x,x) = \frac{1}{2}
A^{(2)}(1\otimes f-g\otimes 1)^2(x,x) - 
(1\otimes f-g\otimes 1)(1\otimes Af-Ag\otimes 1)(x,x)$ 
the {\em pure diffusion form}. It can easily be seen that it is
nonnegative and it vanishes if the associated canonical flow is a flow
of maps. Indeed
\beqarr \Gamma(f,f)(x) 
&=& \lim_{t\to 0}\frac{1}{t}
\left(\P^{(1)}_tf^2(x)-\P^{(2)}_tf^{\otimes 2}(x,x)\right)\\
&=& \lim_{t\to 0}\frac{1}{2t}
\left(\P^{(2)}_t(1\otimes f-f\otimes 1)^2(x,x)\right).
\eeqarr
 The converse is
not true (see examples in section \ref{isotrope}). Diffusive flows for
which the pure diffusion form vanishes may be called {\em turbulent}.

\medskip
The two-point motion $X^{(2)}_t=(X_t,Y_t)$ solves the
following martingale problem associated with $A^{(2)}$~:
\beq \label{martprob} M^{f\otimes g}_t:=f(X_t)g(Y_t)-f(X_0)g(Y_0)-\int_0^t
A^{(2)}(f\otimes g) (X_s,Y_s)~ds\eeq
is a martingale for all $f$ and $g$ in $C_K^2(M)$.

\medskip
Note that for all functions $h_1$ and $h_2$ in $C_K^2(M)\otimes
C_K^2(M)$, the martingale bracket 
$\langle h_1(X^{(2)}),h_2(X^{(2)})\rangle_t$ is equal to 
\beq\int_0^t(A^{(2)}(h_1h_2)-h_1A^{(2)}h_2-h_2A^{(2)}h_1)(X^{(2)}_s)~ds\eeq 
and for all functions $f$ and $g$ in $C^2_K(M)$,
\beq\label{crochet} \langle f(X),g(Y)\rangle_t=\int_0^tC(f,g)(X_s,Y_s)~ds.\eeq

\begin{definition} \bdes \ita A covariance function on the space of vector
fields is a symmetric map $C$ from $T^*M^2$ in $\RR$
such that its restriction to $T^*_xM\times T^*_yM$ is bilinear and for
any $n$-uples $(\xi_1,\dots,\xi_n)$ of $T^*M^2$, $\sum_{1\leq i,j\leq 
n}C(\xi_i,\xi_j)\geq 0$, (see \cite{ljr}). For $f$ and $g$ in
$C_K^1(M)$, we denote $C(df(x),dg(y))$ by $C(f,g)(x,y)$.
\itb We say the covariance function is continuous if $C(f,g)$ is
continuous for all $f$ and $g$ in $C_K^1(M)$. \edes \end{definition}

\bprop \bdes \ita $C$ is a continuous covariance function on the space
of vector fields.
\itb For all $f_1,\dots,f_n$ in $C^2_K(M)$, then $g=f_1\otimes\cdots\otimes
f_n\in {\cal D}(A^{(n)})$ and for $x=(x_1,\dots,x_n)\in M^n$,
\beq A^{(n)}g(x) = \sum_i\prod_{j\neq i}f_j(x_i)Af_i(x_i) +
\sum_{i<j} C(f_i,f_j)(x_i,x_j)\prod_{k\neq i,j}f_k(x_k).
\label{genen} \eeq\edes\eprop 
\prf For all $f$ and $g$ in $C^2_K(M)$, $C(f,g)(x,y)$ is a
function of $df(x)$ and of $dg(y)$ we denote $C(df(x),dg(y))$. Hence
$C$ is a symmetric map from $T^*M^2$ in $\RR$ and its restriction to
$T^*_xM\times T^*_yM$ is bilinear. To prove {\bf (a)}, it remains to prove
$\sum_{i,j}C(\xi_i,\xi_j)\geq 0$ for all $\xi_1,\dots,\xi_n$ in
$T^*M^2$. This holds since, for all $f_1,\dots,f_n$ in $C^2(M)$ and
all $x_1,\dots,x_n$ in $M$,
\beq \sum_{i,j} C(f_i,f_j)(x_i,x_j) =
(A^{(n)}g^2-2gA^{(n)}g)(x_1,\dots,x_n) \eeq
where $g(x_1,\dots,x_n)=\sum_{i=1}^n f_i(x_i)\in{\cal D}(A^{(n)})$.
This expression is nonnegative since $A^{(n)}g^2-2gA^{(n)}g=\lim_{t\to
0}\frac{1}{t}(\P^{(n)}_t g^2-(\P^{(n)}_tg)^2+(\P^{(n)}_t g-g)^2)$. 

\smallskip The proof of {\bf (b)} is an application of It\^o's
formula. \qed

\begin{definition} The diffusion generator $A$ and the
covariance function $C$ are called the local characteristics of the
family $(\P^{(n)}_t,~n\geq 1)$ or of the diffusion convolution
semigroup. \end{definition} 

When there is no pure diffusion, to give the local characteristics
$(A,C)$ in a system of local charts is equivalent to give a drift $b$
and $C$ (this corresponds to the usual definition of a local
characteristics of a stochastic flow) since in this case $a^{i,j}(x) =
c^{i,j}(x,x)$.

\brem  When $(\P^{(n)}_t,~n\geq 1)$ satisfies {\bf (C)}, {\bf
(i)} and {\bf (ii)}, then $(\P^{(n),c}_t,~n\geq 1)$ also satisfies
{\bf (i)} if and only if for all $x$ in $M$ and all $f$, $g$ in
$C^2_K(M)$, $C(f,g)(x,x)=\Gamma(f,g)(x)$ (this holds since we have
$C(f,g)(x,x)-\Gamma(f,g)(x)=\lim_{t\to 0}
\frac{1}{t}(\P^{(2),c}_t(f\otimes g)(x,x)-\P^{(1)}_t(fg)(x))$),
i.e. when there is no pure diffusion. So that the results of this
section also apply to $(\P^{(n),c}_t,~n\geq 1)$.

Then in this case $(\P^{(n)}_t,~n\geq 1)$ and $(\P^{(n),c}_t,~n\geq
1)$ have the same local characteristics. \erem

\medskip
Let $K=(K_{s,t},~s\leq t)$ be a measurable stochastic flow of kernels
  associated with $(\P^{(n)}_t,~n\geq 1)$ defined on a probability
  space $(\Omega,{\cal A},\P)$. 

\begin{definition}
Let $C$ be a covariance function on the space of vector fields.
A two parametric family $W=(W_{s,t},~s\leq t)$ of random variables
taking their values in the space of vector fields on $M$ is
called a vector field valued white noise of covariance $C$ if 
\bdes\iti for all
$s_i\leq t_i\leq s_{i+1}$, the random variables $(W_{s_i,t_i},~1\leq
i\leq n)$ are independent, 
\itii for all $s\leq u\leq t$,
$W_{s,t}=W_{s,u}+W_{u,t}$ a.s. and 
\itiii for all $s\leq t$, $\{\langle
W_{s,t},\xi\rangle,~\xi\in T^*M\}$  \footnote{when $\xi=(x,u)$,
$\langle W_{s,t},\xi\rangle=\langle W_{s,t}(x),u\rangle$.} is a
centered Gaussian process of covariance given by
\beq \E[\langle W_{s,t},\xi\rangle\langle W_{s,t},\xi'\rangle] =
(t-s)C(\xi,\xi'),\eeq 
for $\xi$ and $\xi'$ in $T^*M$. \edes
\end{definition}

In this section, we intend to define on $(\Omega,\cA,\P)$ a
vector field valued white noise $W$ of covariance $C$ such that
$K$ solves a SDE driven by $W$. 

\medskip
In section \ref{noise}, under an additional assumption, we will prove
that the linear (or Gaussian) part of the noise generated by
$K$ (in the case it is the canonical flow) is the
noise generated by the vector field valued white noise $W$. 

\subsection{The velocity field $W$.}\label{constrW}
For all $s\leq t$, all $f\in C_K^2(M)$ and all $x\in M$, let
\beq M_{s,t}f(x)=K_{s,t}f(x)-f(x)-\int_{s}^{t}K_{s,u}(Af)(x)~du.\eeq
\blem\label{lemmart} For all $s\in \RR$, $f\in C_K^2(M)$ and $x\in M$,
$M^{f,x}_s=(M_{s,t}f(x),~t\geq s)$ is a martingale with 
respect to the filtration ${\cal F}^s=({\cal F}_{s,t},~t\geq s)$ and 
\beq \frac{d}{dt}\langle M^{f,x}_s,M^{g,y}_s\rangle_t =
K_{s,t}^{\otimes 2} C(f,g)(x,y),\eeq 
for all $f$, $g$ in $C^2_K(M)$ and all $x$, $y$ in $M$. \elem
\prf Since $K$ is a measurable stochastic flow of kernels and
that for all positive $h$ and all $f$ in $C^2_K(M)$, a.s.
\beq\label{Mst}M_{s,t+h}f(x)-M_{s,t}f(x) = K_{s,t}(M_{t,t+h}f)(x),\eeq
$M^{f,x}_s$ is a martingale. Note that equation (\ref{Mst}) also implies
that for all positive $h$, all $f$, $g$ in $C^2_K(M)$ and all $x$, $y$
in $M$,
\beqarr \E[(M_{s,t+h}f(x)-M_{s,t}f(x))(M_{s,t+h}g(y)-M_{s,t}g(y))
|{\cal F}_{s,t}] \hskip-100pt &&\\ 
&=& K_{s,t}^{\otimes 2} (\E[M_{t,t+h}f\otimes
M_{t,t+h}g])(x,y). \eeqarr
The stationarity implies that
$\E[M_{t,t+h}f(x)M_{t,t+h}g(y)]=\E[M_{0,h}f(x)M_{0,h}g(y)]$.
Computation using the fact that
$\P^{(1)}_tf-f=\int_0^t\P^{(1)}_sAf~ds$  and $\P^{(2)}_t(f\otimes
g)-f\otimes g=\int_0^t\P^{(2)}_sA^{(2)}(f\otimes g)~ds$  gives 
\beq \E[M_{0,h}f(x)M_{0,h}g(y)]=\int_0^h\P^{(2)}_s(C(f,g))(x,y)~ds.\eeq
Since $\P^{(2)}_t$ is Feller and $C(f,g)$ is continuous with compact
support,
\beq \label{DL} \E[M_{0,h}f(x)M_{0,h}g(y)]=h~C(f,g)(x,y)+o(h),\eeq
uniformly in $(x,y)\in M^2$.

Therefore $\E[(M_{s,t+h}f(x)-M_{s,t}f(x))(M_{s,t+h}g(y)-M_{s,t}g(y)) 
|{\cal F}_{s,t}]$ is equivalent as $h$ tends to 0 to
$h~K_{s,t}^{\otimes 2}C(f,g)(x,y)$. This proves the lemma. \qed

\brem  In the case of Arratia's coalescing flow
$(\p_{s,t},~s\leq t)$, $C=0$ but $\frac{d}{dt}\langle
M^{f,x}_s,M^{g,y}_s\rangle_t =
1_{\{\p_{s,t}(x)=\p_{s,t}(y)\}}$. In this case, $C_K^2(M)\otimes
C_K^2(M)$ is not included in ${\cal D}(A^{(2)})$. This property also
fails for the coalescing flow associated with Tanaka's SDE. \erem

For all $s<t$, $n\geq 1$ and $0\leq k\leq 2^n-1$, let
$t_k^n=s+k2^{-n}(t-s)$ and
\beq W^n_{s,t}f=\sum_{k=0}^{2^n-1} M_{t_k^n,t_{k+1}^n}f, \eeq
where $f\in C^2_K(M)$.
Note that $(M_{t_k^n,t_{k+1}^n})_{0\leq k\leq 2^n -1}$ are independent
equidistributed random variables.

\subsubsection{Convergence in law.}
\blem\label{tcl}
For all $s<t$ and all $((x_i,f_i),~1\leq i\leq m)\in (M\times
C^2_K(M))^m$, then $\sum_{i=1}^m W^n_{s,t}f_i(x_i)$ converges
in law towards $\sum_{i=1}^m W_{s,t}f_i(x_i)$ as $n$ tends to
$\infty$, where $W$ is a vector field valued white noise of covariance
$C$. \elem
\prf Using lemma \ref{lemmart}, we have for all $f$, $g$ in $C^2_K(M)$
and all $x$, $y$ in $M$,
\beqar \E[M_{t_k^n,t_{k+1}^n}f(x)M_{t_k^n,t_{k+1}^n}g(y)] &=&
\int_0^{2^{-n}(t-s)}\P^{(2)}_u C(f,g)(x,y)~du\non\\
&=&2^{-n}(t-s)C(f,g)(x,y)+o(2^{-n})\qquad\label{petit}
\eeqar
and this developement is uniform in $x$ and $y$ in $M$.

\smallskip We will only prove the proposition when $m=1$ (the proof being
the same for $m>1$). The proposition is just an application of the central limit
theorem for arrays (see \cite{bilingsley}), which we can apply
since (\ref{petit}) is satisfied provided the Lyapounov condition
\beq\label{lyapcond}
\lim_{n\to\infty}\sum_{k=0}^{2^n-1}\E[|M_{t_k^n,t_{k+1}^n}f(x)|^{2+\delta}]=0\eeq
for some positive $\delta$, is satisfied.

\medskip
Using Burkholder-Davies-Gundy's inequality and lemma \ref{lemmart},
\beqarr
\E[|M_{t_k^n,t_{k+1}^n}f(x)|^{{2+\delta}}] &\leq& C
\E\left[\left(\int_0^{2^{-n}(t-s)} K_{0,u}^{\otimes 2}
(C(f,f))(x,x)~du\right)^{\frac{2+\delta}{2}}\right]\\
&\leq& C2^{-\frac{(2+\delta)n}{2}},
\eeqarr
where $C$ is a constant (changing every line) depending only on $f$, $(t-s)$
and $\delta$. This implies
\beq\sum_{k=0}^{2^n-1}\E[|M_{t_k^n,t_{k+1}^n}f(x)|^{2+\delta}]
\leq C2^n2^{-\frac{(2+\delta)n}{2}}\leq C2^{-\frac{n\delta}{2}}.\quad \qed \eeq

\brem  For Arratia's coalescing flow, one can show the convergence in
law as $n$ goes to $\infty$ of $(W^n_{s,t}(x_1),\dots,W^n_{s,t}(x_k))$ towards
$(B^1_{s,t},\dots,B^k_{s,t})$, where $(B^1,\dots,B^k)$ is a $k$-dimensional
white noise. \erem

\subsubsection{Convergence in $L^2(\P)$.}
In the preceding section, we proved the convergence in law of $W^n$
towards a vector field valued white noise of covariance $C$. In this
section, we prove that this convergence holds in $L^2(\P)$.
\blem\label{convl2} For all $s<t$ and all $(x,f)\in M\times C_K^2(M)$,
$W^n_{s,t}f(x)$ converges in $L^2(\P)$.\elem
\prf For all $f\in C^2_K(M)$, all $x\in M$ and all $s<t$,
\beqar \E[(W_{s,t}^nf(x)-W^{n+k}_{s,t}f(x))^2] &=& \E[(W_{s,t}^nf(x))^2] +
\E[(W^{n+k}_{s,t}f(x))^2]\non\\
&& -\quad 2\E[W_{s,t}^nf(x)W^{n+k}_{s,t}f(x)]\label{lastterm} \eeqar
Elementary computations using equation (\ref{DL}) implies
\beqar \E[(W_{s,t}^nf(x))^2] &=& (t-s)~C(f,f)(x,x)+o(1)\\
\E[(W_{s,t}^{n+k}f(x))^2] &=& (t-s)~C(f,f)(x,x)+o(1) \eeqar
as $n$ goes to $\infty$ and this uniformly in $k\in\NN$. Using the
independence of the increments, the last term (\ref{lastterm})
rewrites
\beq \E[W_{s,t}^nf(x)W^{n+k}_{s,t}f(x)] =
\sum_{i=0}^{2^n-1}\sum_{j=i2^k}^{(i+1)2^k-1}
\E[M_{t_i^n,t_{i+1}^n}f(x)M_{t_j^{n+k},t_{j+1}^{n+k}}f(x)]. \eeq
Note that for $s\leq u\leq v\leq t$, using first the martingale
property, then equation (\ref{DL}) and the uniform continuity of
$C(f,f)$, we have
\beqar \E[M_{s,t}f(x)M_{u,v}f(x)] &=& \E[M_{s,v}f(x)M_{u,v}f(x)]\non\\
&=& \E[(K_{s,u}\otimes I)(M_{u,v}f\otimes M_{u,v}f)(x,x)]\non\\
&=& \E[(K_{s,u}\otimes I)(\E[M_{u,v}f\otimes M_{u,v}f])(x,x)]\non\\
&=& (v-u)~C(f,f)(x,x)+o(v-u), \label{DLuv} \eeqar
uniformly in $x\in M$. This implies 
\beq \E[W_{s,t}^nf(x)W^{n+k}_{s,t}f(x)] = (t-s)~C(f,f)(x,x)+o(1) \eeq
as $n$ tends to $\infty$ and uniformly in $k\in\NN$. We therefore have 
\beq \lim_{n\to\infty}\sup_{k\in\NN}
\E[(W_{s,t}^nf(x)-W^{n+k}_{s,t}f(x))^2] = 0,\eeq
i.e. $(W^n_{s,t}f(x),~n\in\NN)$ is a Cauchy sequence in $L^2(\P)$. This
proves the lemma. \qed

\brem  For Arratia's coalescing flow, this lemma is not
satisfied since $(W^n_{s,t}f(x),~n\in\NN)$ fails to be a Cauchy sequence
in $L^2(\P)$. \erem

Thus, for all $s<t$, we have defined the vector field valued random
variable $W_{s,t}$ such that $W_{s,t}f(x)$ is the $L^2(\P)$-limit of
$W^n_{s,t}f(x)$ for all  $x\in M$ and all $f\in C(M)$. Then, using
lemma \ref{tcl}, it is easy to see that $W=(W_{s,t},~s\leq t)$ is a
vector field valued white noise of covariance $C$.

\subsection{The stochastic flow of kernels solves a SDE.}\label{solveSDE}

In \cite{ljr}, it is shown that a vector field valued white noise $W$
of covariance $C$ can be constructed with a sequence of independent
real white noises $(W^\alpha)_\alpha$ by the formula $W=\sum_\alpha
V_\alpha W^\alpha$, where $(V^\alpha)_\alpha$ is an orthonormal basis 
of $H_C$, the self-reproducing space associated with $C$.

For all predictable (with respect to the filtration $({\cal
F}_{-\infty,t},~t\in\RR)$) process $(H_t)_{t\in\RR}$ taking its values
in the dual of $H_C$, we define the stochastic integral of $H$ with
respect to $W$ by the formula 
\beq \int_s^t H_u(W(du)) = \sum_\alpha \int_s^t 
\langle H_u,V_\alpha\rangle~W^\alpha(du),\eeq 
for $s<t$. Note that the above definition is independent of the choice
of the orthonormal basis $(V^\alpha)_\alpha$.

In particular this applies to $H_u(V)=K_{s,u}(Vf)(x)1_{s\leq u<t}$ for
$f\in C_K(M)$ and $x\in M$. Then the stochastic integral $\sum_\alpha\int_s^t
K_{s,u}(V^\alpha f)W^\alpha(du)$ is denoted
\beq \label{SI} \int_s^t K_{s,u}(Wf(du))(x). \eeq

\brem \label{approxSI} {\em The stochastic integral
(\ref{SI}) is equal to the limit in $L^2(\P)$ of 
$$\sum_{k=0}^{2^n-1}
K_{s,t_k^n}(W_{t_k^nt_{k+1}^n}f)(x)$$ as $n$ tends to $\infty$, where
$t_k^n=s+k2^{-n}(t-s)$. Indeed,
\beqarr \E\left[\left(\int_s^t K_{s,u}(Wf(du))(x) - \sum_{k=0}^{2^n-1}
K_{s,t_k^n}(W_{t_k^nt_{k+1}^n}f)(x)\right)^2\right] &=& \\
&&\hskip-250pt =\quad \sum_{k=0}^{2^n-1} \int_{t_k^n}^{t_{k+1}^n}
\P^{(2)}_{t_k^n-s} (I+\P^{(2)}_{u-t_k^n}-2I\otimes \P^{(1)}_{u-t_k^n})
C(f,f)(x,x)~du\eeqarr
which tends to 0 as $n$ tends to $\infty$.}\erem

\bprop\label{solforte}
$W$ is the unique vector field valued white noise
such that for all $s<t$, $x\in M$ and $f\in C_K^2(M)$, $\P$-almost
surely, 
\beq\label{SDE}
K_{s,t}f(x)=f(x)+\int_s^t K_{s,u}(Wf(du))(x)+\int_s^t K_{s,u}(Af)(x)~du.\eeq
{\em Note that giving the local characteristics of the flow is equivalent
to giving this SDE. This SDE will be called the $(A,C)$-SDE.}
\eprop
\prf For all $s<t$, from remark \ref{approxSI},
\beq\int_s^t K_{s,u}(Wf(du))(x)=\lim_{n\to\infty}\sum_{k=0}^{2^n-1}
K_{s,t_k^n}(W_{t_k^n,t_{k+1}^n}f)(x)\eeq
in $L^2(\P)$, where $t_k^n=s+k2^{-n}(t-s)$.

\smallskip
For all integers $i,l,k$ and $n$ such that $l\geq n$ and
$k2^{l-n}\leq i\leq (k+1)2^{l-n}-1$, the development (\ref{DLuv})
implies 
\beq \E[M_{t_i^l,t_{i+1}^l}f(x)M_{t_k^n,t_{k+1}^n}f(x)] =
2^{-l}(t-s)C(f,f)(x,x)+o(2^{-l}),\eeq 
uniformly in $x\in M$. This implies that for $l\geq n$,
\beq \E\left[\left(\sum_{i=k2^{l-n}}^{(k+1)2^{l-n}-1}
M_{t_i^l,t_{i+1}^l}f(x) - M_{t_k^n,t_{k+1}^n}f(x)\right)^2\right]=o(2^{-n}),\eeq
uniformly in $x\in M$. Taking the limit as $l$ goes to $\infty$, we
get 
\beq \E[(W_{t_k^n,t_{k+1}^n}f(x)-M_{t_k^n,t_{k+1}^n}f(x))^2]=o(2^{-n}),\eeq
uniformly in $x\in M$. We use this estimate to prove that
\beq\label{l2lim} \int_s^t K_{s,u}(W(du)f)(x) =
\lim_{n\to\infty}
\sum_{k=0}^{2^n-1}K_{s,t_k^n}(M_{t_k^n,t_{k+1}^n}f)(x) \eeq  
in $L^2(\P)$. This holds since
\beqarr
\E\left[\left(\sum_{k=0}^{2^n-1}K_{s,t_k^n}(W_{t_k^n,t_{k+1}^n}f) -
\sum_{k=0}^{2^n-1}K_{s,t_k^n}(M_{t_k^n,t_{k+1}^n}f)\right)^2(x)\right]&=&\\ 
&&\hskip-180pt=\quad \sum_{k=0}^{2^n-1}
\E[(K_{s,t_k^n}(W_{t_k^n,t_{k+1}^n}f - M_{t_k^n,t_{k+1}^n}f))^2(x)]\\
&&\hskip-180pt \leq\quad\sum_{k=0}^{2^n-1}
\P^{(1)}_{t_k^n-s}\left(\E[(W_{t_k^n,t_{k+1}^n}f -
M_{t_k^n,t_{k+1}^n}f)^2]\right)(x)\\
&&\hskip-180pt \leq\quad 2^n~o(2^{-n})=o(1).
\eeqarr  

Note now that 

\beqarr 
\sum_{k=0}^{2^n-1}K_{s,t_k^n}(M_{t_k^n,t_{k+1}^n}f)(x)
&=& \sum_{k=0}^{2^n-1}
K_{s,t_k^n}\left(K_{t_k^n,t_{k+1}^n}f - f -
\int_{t_k^n}^{t_{k+1}^n}K_{t_k^n,u}(Af)~du\right)(x)\\
&=& K_{s,t}f(x)-f(x)-\int_s^t K_{s,u}(Af)(x)~du.
\eeqarr
This proves that $K$ solves the $(A,C)$-SDE driven
by $W$. Finally, note that if $K$ solves the $(A,C)$-SDE
driven by a vector field valued white noise $W'$ then we
must have $W'=W$. \qed

\medskip Let $X=(X_t,~t\geq 0)$ be the Markov process defined in
section \ref{cadlag} on $(\Omega\times C(\RR^+,M),
\cA\otimes {\cal B}(C(\RR^+,M)), \P(d\omega)\otimes
\P_{x,\omega}(d\omega'))$ by $X(\omega,\omega')=\omega'$.
\bprop \label{solfaible} Assume there is no pure diffusion
(i.e. for all $f\in C_K^2(M)$ and all $x\in M$,
$\Gamma(f)(x)=C(f,f)(x,x)$). Then, for all $t\geq 0$, $x\in M$ and
$f\in C_K^2(M)$, $\P(d\omega)\otimes \P_{x,\omega}(d\omega')$-almost
surely,
\beq\label{SDEf}
f(X_t)=f(x)+\int_0^t W(du)f(X_u)+\int_0^t Af(X_u)~du,\eeq
i.e. $X$ is a weak solution of this SDE (in the sense given in \cite{R-W}).
\eprop
\prf Like in the proof of (\ref{l2lim}) in proposition \ref{solforte},
we show that 
\beq \label{approxIS}\int_0^t W(du)f(X_u) = \lim_{n\to\infty} \sum_{k=0}^{2^n-1}
M_{t_k^n,t_{k+1}^n}f(X_{t_k^n}) \eeq
in $L^2(\P_x)$, with $\P_x=\P(d\omega)\otimes \P_{x,\omega}(d\omega'))$.
Let $M^f_t=f(X_t)-f(x)-\int_0^t Af(X_u)~du$, then $(M^f_t,~t\geq 0)$
is a martingale relative to the filtration $({\cal F}^X_t,~t\geq 0)$
generated by the Markov process $X$. We now prove that
$\E_x[(M^f_t-\int_0^t W(du)f(X_u))^2]=0$, where $\E_x$ denotes the
expectation with respect to $\P_x$. It is easy to see that, since there
is no pure diffusion,
\beq \E_x[(M^f_t)^2]=\E_x[(\int_0^t W(du)f(X_u))^2]=\E_x[\int_0^t
C(f,f)(X_u,X_u)~du]. \eeq
Equation (\ref{approxIS}) and the martingale property of $M^f_t$
implies that 
\beqar \E_x[M^f_t\int_0^t W(du)f(X_u)] &=& \lim_{n\to\infty} \sum_{k=0}^{2^n-1}
\E_x[M^f_{t_{k+1}^n}\times M_{t_k^n,t_{k+1}^n}f(X_{t_k^n})].\\
&=& \lim_{n\to\infty} \sum_{k=0}^{2^n-1}
\E_x[(M^f_{t_{k+1}^n}-M^f_{t_k^n})\times
M_{t_k^n,t_{k+1}^n}f(X_{t_k^n})].\non\eeqar 
Since for all $0\leq s<t$, $\E_x[M^f_t-M^f_s|\cA\vee{\cal F}^X_s]
= M_{s,t}f(X_s)$, we get 
\beqar \E_x[M^f_t\int_0^t W(du)f(X_u)] &=& \lim_{n\to\infty}
\sum_{k=0}^{2^n-1}
\E_x[(M_{t_k^n,t_{k+1}^n}f(X_{t_k^n}))^2]\non\\
&=& \E_x[(\int_0^t W(du)f(X_u))^2]. \eeqar
Therefore $\E_x[(M^f_t-\int_0^t W(du)f(X_u))^2]=0$. \qed

\subsection{The $(A,C)$-SDE.}\label{34}
In this section and in the following, let $A$ be a second order
differential operator mapping $C^2_K(M)$ in $C_K(M)$ and $C$ a
continuous covariance on vector fields.
\begin{definition} Let $K$ be a stochastic flow of kernels and $W$ a
  vector field valued white noise, defined on a probability space
  $(\Omega,\cA,\P)$.
\bdes \iti $(K,W)$ is a solution of the $(A,C)$-SDE if the covariance
of $W$ is $C$ and $(K,W)$ satisfies (\ref{SDE}) for all $s<t$, $x\in
M$ and $f\in C_K^2(M)$.
\itii $(K,W)$ is called a strong solution of the $(A,C)$-SDE if
moreover for all $s\leq t$, $K_{s,t}$ is $\cF^W_{s,t}$-measurable,
where $ \cF^W_{s,t}$ is the completion by all $\P$-negligible sets of
  $\cA$ of the $\sigma$-field $\sigma(W_{u,v},~s\leq u\leq v\leq
  t)$.
\itiii When a solution $(K,W)$ of the $(A,C)$-SDE is not a strong
  solution, we say it is a weak solution. \edes
\end{definition}

\brem  Let $(K,W)$ be a solution of the $(A,C)$-SDE and $\nu$ the
Feller convolution semigroup associated with $K$. Then $\nu$ is a
diffusion convolution semigroup with local characteristics $(A,C)$.\erem
The proof of this remark is left to the reader.

\brem {\em The fact that $(K,W)$ is a strong (respectively a weak)
  solution of the $(A,C)$-SDE only depends on the law of $K$. So that
  we can say shortly that $K$ is a strong (respectively a weak)
  solution of the $(A,C)$-SDE.} \erem 

\begin{definition} We will say that $(\P^{(n)}_t,~n\geq 1)$, a
  compatible family of Feller semigroups, or $\nu=(\nu_t)$, a Feller
  convolution semigroups, defines a strong (respectively a weak)
  solution of the $(A,C)$-SDE if $\P_\nu$ is the law of a stochastic
  flow of kernels, which is a strong (respectively a weak) solution of
  the $(A,C)$-SDE. \end{definition}

Under some additional assumptions, we will give in section \ref{noise}
a representation of all solutions of the $(A,C)$-SDE.

\begin{definition} We say that (strong) uniqueness holds for the
$(A,C)$-SDE when there is only one diffusion convolution semigroup
with local characteristics $(A,C)$ defining a (strong)
solution. \end{definition}

\subsection{Strong solution and filtering.}
Let us now consider the canonical flow associated with $\nu$, a
diffusion convolution semigroup, with local characteristics
$(A,C)$. Let $N^W_\nu:=
(\Omega^0,\cA^0,({\cal F}^W_{s,t})_{s\leq t},
\P_\nu,(T_h)_{h\in\RR})$ be the noise generated by the vector field  
valued white noise $W$. Note that $N^W_\nu$ is a linear or Gaussian
\footnote{The noise $({\cal G}_{s,t})_{s\leq t}$ is Gaussian if and
only if there exists a countable family of independent real white
noises $\{W^\alpha\}$ such that, up to negligible sets, ${\cal
G}_{s,t}$ is generated by the random variables $W^\alpha_{u,v}$
for all $s\leq u\leq v\leq t$ and all $\alpha$.} sub-noise of
$N_{\nu}$, the noise generated by the canonical flow.

\medskip Let $\bar{K}=(\bar{K}_{s,t},~s\leq t)$ be the stochastic flow
of kernels obtained by filtering the canonical flow with respect to
the sub-noise $N^W$ (see section \ref{filtering}). It is easy to see
that $\bar{K}$ also solves the $(A,C)$-SDE (see the proof of lemma 3.9 in
\cite{ljr}) and has the same local characteristics as the canonical
flow. Since, for all $s\leq t$, $\bar{K}_{s,t}$ is ${\cal
F}^W_{s,t}$-measurable, $(\bar{K},W)$ is a strong solution of the
$(A,C)$-SDE. Let $\nu^s$ denote the associated diffusion convolution
semigroup.

For any $f\in C_0(M)$ and $x\in M$, $\bar{K}_{s,t}f(x)$ can be
expanded into a sum of Wiener chaos elements, i.e. iterated Wiener
integrals of the form 
\beq \sum_{\alpha_1,\dots,\alpha_n}\int
C^{\alpha_1,\dots,\alpha_n}(s_1,\dots,s_n) ~dW_{s_n}^{\alpha_n}\cdots
dW^{\alpha_1}_{s_1}.\eeq
 Since $W$ was constructed from the flow, it is
clear that the functions $C^{\alpha_1,\dots,\alpha_n}$ are determined
by the law of the flow (we will give, under some additional
assumptions, an explicit form of them in the following section).

\subsection{The Krylov-Veretennikov expansion.}\label{Veretennikov}
We  still assume we are given $\nu=(\nu_t)_{t\geq 0}$ a diffusion
convolution semigroup, in the sense of section \ref{hypo}, associated
with a set of local characteristics $(A,C)$.

We suppose in this section the existence of a Radon measure
$m$ on $M$ such that $A$ is symmetric with respect to $m$.

Moreover, we assume that $\hbox{Im}(I-A)$ is dense in $C_0(M)$ (it
implies that $\P^{(1)}_t$ is symmetric with respect to $m$ and is the
unique Feller semigroup whose generator extends $A$).

\smallskip
Following \cite{ljr}, starting from the vector field valued white noise $W$,
one can define $(S_{s,t},~s\leq t)$ a stochastic flow of Markovian
operators (acting on $L^2(m)$) such that for all $s\leq t$, $S_{s,t}$
is $\sigma(W)$-measurable and for $f\in L^2(m)$ and $s\leq u\leq t$,
\beqarr S_{s,t}f&=&S_{s,u}S_{u,t}f, \\
S_{s,t}f&=&\P^{(1)}_{t-s}f +\int_s^t S_{s,u} W(du) \P^{(1)}_{t-u} f,
\eeqarr
where both equalities hold in $L^2(m\otimes \P)$. These operators are given
by the Wiener chaos expansion (called Krylov-Veretennikov expansion)
\beq \label{chaos} S_{s,t}f = \P^{(1)}_tf + \sum_{n\geq 1}J^n_{s,t}f, \eeq
with
\beq \label{Jn} J_{s,t}^nf = 
\int_{s\leq s_1\leq \cdots\leq s_n\leq t}
\P^{(1)}_{s_1-s}W(ds_1) \P^{(1)}_{s_2-s_1}\cdots
\P^{(1)}_{s_n-s_{n-1}}W(ds_n)\P^{(1)}_{t-s_n}f. \eeq
They can be characterized (theorem 3-2 in \cite{ljr}) as the unique
flow of random operators on $L^2(m)$, $\sigma(W)$-measurable, such that
$\E[(S_{s,t}f)^2]\leq \P^{(1)}_{t-s}f^2$ and 
\beq \label{caractS}
 S_{s,t}f-f=\int_s^t S_{s,u} W(du) f +
\frac{1}{2}\int_s^t S_{s,u}\bar{A}f~du \quad \hbox{ in } L^2(m\otimes \P)
\eeq
for all $f$ in the domain of the $L^2$-generator $\bar{A}$, denoted
$\cD(\bar{A})$. It implies the following 
\bprop \bdes \ita If $\nu$ defines a strong solution $(K,W)$ of the
$(A,C)$-SDE, then  for all $s\leq t$, $m\otimes \P$-a.e., for all
$f\in C_K(M)$, 
\beq\label{mpae} K_{s,t}f=S_{s,t}f \eeq
\itb Strong uniqueness holds. \edes  \eprop
\prf {\bf (a)}
It is clear that $K$ induces a flow of Markovian operators on
$L^2(m)$ which verifies (\ref{caractS}) for $f\in C^2_K(m)$. Then
(\ref{caractS}) extends to functions in the domain of the Feller
generator and finally to $\cD(\bar{A})$.

{\bf (b)} From {\bf (a)}, it is clear that $m^{\otimes n}$-a.e.,
$\P^{(n)}_t=\E[S_{0,t}^{\otimes n}]$.  Since it is a
Feller semigroup, it is uniquely determined. \qed

\setcounter{equation}{0}
\section{Noise and classification.}\label{noise} 
\subsection{Assumptions.}
In this section, as before $M$ denotes a smooth locally compact
manifold. We fix a pair of local characteristics $(A,C)$ on $M$.
$A$ is a second order differential operator mapping $C^2_K(M)$ in
$C_K(M)$ and $C$ a continuous covariance on vector fields.
The associated differential operators $A^{(n)}$ on $C^2_K(M)^{\otimes
  n}$ are defined by equation (\ref{genen}).

\medskip
Let $\cM(n,x)$ be the following martingale problem associated with
$A^{(n)}$ and $x\in M^n$~:
\bdes \item[] There exists a probability space on which is
constructed a $M^n$-valued stochastic process
$X^{(n)}=(X^{(n)}_t,~t\geq 0)$ such that 
\beq f(X^{(n)}_t)-f(x)-\int_0^t A^{(n)}f(X^{(n)}_s)~ds \eeq
is a martingale for all test function $f$ in
$C_K^2(M)\otimes\cdots\otimes C_K^2(M)$.\edes

\medskip
We suppose that the local characteristics $(A,C)$ verify the following
assumption
\bdes\item[(U)] For all $n\geq 1$, the martingale problem ${\cal M}(n,x)$ has a
unique solution in law on the set of continuous trajectories stopped
at $\Delta_n$.\edes

\brem  {\em Condition {\bf (U)} is satisfied when the
coefficients of the local characteristics are $C^2$ outside of
$\Delta_n$ (see theorem 12.12 and section V.19 in \cite{R-W}) or when
$A^{(n)}$ is elliptic outside of $\Delta_n$ (see section V.24 in
\cite{R-W}).} \erem

Our purpose is to classify Feller convolution semigroups associated with
these local characteristics. We will treat two cases
\bdes \item[(A)] The non coalescing case where the solution of the
martingale problem ${\cal M}(2,x)$ does not hit the diagonal when
$x=(x_1,x_2)$ with $x_1\neq x_2$.
\item[(B)] The coalescing case where
\bdes \item[] there is no pure diffusion
(i.e. $(\frac{1}{2}Af^2-fAf)(x)=C(f,f)(x,x)$ for all $f\in C^2_K(M)$ and $x\in
M$) \edes 
and where assumption
\bdes \item[(C)] For all $t>0$, $\eps>0$ and $x\in M$, $\lim_{y\to x}
\P^{(2)}_{(x,y)}[\{T_\Delta>t\}\cap\{d(X_t,Y_t)>\eps\}]=0$ and for some
$(x,y)\in M^2$, $\P^{(2)}_{(x,y)}[T_\Delta<\infty]>0$
\edes
holds for $X^{(2)}_t=(X_t,Y_t)$ a solution of ${\cal M}(2,x)$.
\edes

When the local characteristics are non coalescing (case {\bf (A)}),
these local characacteristics are associated with at most a unique
convolution semigroup and a unique canonical flow (which is not always
a flow of maps). From section \ref{34} we know
the latter has to be a strong solution of the SDE (otherwise
uniqueness would be violated). Assumption {\bf (F)} (see section
\ref{Lipschitz}) is a sufficient
(but not necessary) condition for existence. The family of semigroups given in the
example of Lipschitz SDE's (see section \ref{Lipschitz}) satisfies
these assumptions.

\medskip
In sections \ref{classi}, \ref{MR} and \ref{linear}, we assume 
{\bf (B)} is satisfied.

\subsection{The coalescing case~: classification.}\label{classi}
Following Harris \cite{Harris}, $\cM(n,x)$ has a unique
solution in law on the set of coalescing trajectories,
i.e. $X^{(n)}(\omega)\in C^{(n)}$ where $C^{(n)}$ is the set of
continuous functions $f~:~\RR^+\to M^n$ such that if $f_i(s)=f_j(s)$
for $1\leq i,j\leq n$ and $s\geq 0$ then for all $t\geq s$, 
$f_i(t)=f_j(t)$ (In \cite{Harris}, this martingale problem is solved
when $M=\RR$, but the proof can obviously be adapted to our
framework). Since {\bf (C)} holds, remark \ref{remfel} implies that
the associated semigroups are Feller.

Hence all coalescing flows with these local characteristics have the
same law $\P_{\nu^c}$. They induce the same family of semigroups
$(\P^{(n),c}_t,~n\geq 1)$ and the same diffusion convolution semigroup
$\nu^c$. 

\medskip
Let $N_{\nu^c}$ be the noise generated by the canonical coalescing
flow asociated with the local characteristics $(A,C)$.

\medskip
Let $W$ be the vector field valued white noise defined on
$(\Omega^0,\cA^0,\P_{\nu^c})$ in section \ref{secW} and
$N^W_{\nu^c}$ the sub-noise of $N_{\nu^c}$ generated by $W$. Then
$N^W_{\nu^c}$ is a Gaussian sub-noise of $N$ and it is possible to 
represent it by a countable family of independent real white noises
$\{W^\alpha\}$ such that $W=\sum_\alpha V_\alpha W^\alpha$, where
$\{V_\alpha\}$ is a countable family of vector fields on $M$.

We denote by $\nu^{s}$ the diffusion convolution semigroup
associated with the flow obtained by filtering the canonical coalescing
flow of law $\P_{\nu^c}$ with respect to $N^W_{\nu^c}$.

\medskip
The following theorem gives a representation of all flows with
the same local characteristics. They lie ``between'' the strong
solution and the coalescing solution of the SDE which are distinct
when the coalescing solution is not a strong solution of the SDE.

\bthm\label{cor} Suppose we are given a set of local characteristics
$(A,C)$ and that assumption {\bf (B)} is verified.
\bdes \ita $\nu^c$ is the unique diffusion convolution semigroup
associated with $(A,C)$ and defining a flow of maps (which is
coalescing).
\itb $\nu^{s}$ is the unique diffusion convolution semigroup
associated with $(A,C)$ and defining a strong solution of the
$(A,C)$-SDE.
\itc The diffusion convolution semigroups associated with $(A,C)$ are
all the Feller convolution semigroups weakly dominated by $\nu^c$ and
dominating $\nu^{s}$. \edes  \ethm 

Note that $\nu^c$ and $\nu^s$ are not necessarily distinct.

\medskip
\prf We have already proved {\bf (a)} at the begining of this
section. Theorem \ref{filtrage} implies that every diffusion
convolution semigroup $\bar{\nu}$ with local characteristics $(A,C)$
is weakly dominated by $\nu^c$ so that a stochastic flow $\bar{K}$ of
law $\P_{\bar{\nu}}$ can be obtained by filtering on an extension
$(N,\p)$ of $N_{\nu^c}$ the coalescing flow $\p$ with respect to a
sub-noise $\bar{N}$ of $N$.

Let $\bar{W}$ be the velocity field associated with $\bar{K}$.
Proposition \ref{solforte} shows that $(\bar{K},\bar{W})$ solves the
$(A,C)$-SDE. Notice that $\bar{W}$ can be obtained by filtering $W$
with respect to $\bar{N}$. Indeed, section \ref{constrW} shows that
$\bar{W}^n_{s,t}$ (defined from $\bar{K}$) converges (in $L^2$)
towards $\bar{W}_{s,t}$ and we have that for all $s\leq t$, $f\in
C^2_K(M)$ and $x\in M$,
$\bar{W}^n_{s,t}f(x) = \E[W^n_{s,t}f(x) | \bar{\cal F}_{s,t}]$
a.s. and therefore that $\bar{W}_{s,t}f(x) = \E[W_{s,t}f(x) | \bar{\cal
F}_{s,t}]$ a.s. Since $\bar{W}$ and $W$ have the same law, we must
have $W_{s,t} = \bar{W}_{s,t}$ a.s. This proves that $\bar{\nu}$
dominates $\nu^{s}$.

\smallskip
Let us now suppose that $(\bar{K},\bar{W})$ is a strong solution of the
$(A,C)$-SDE. Then, since $\bar{W}=W$, we must have $N^W_{\nu^c}=\bar{N}$
(since $\bar{K}_{s,t}$ is ${\cal F}^W_{s,t}$-measurable) and thus
$\nu^{s}=\bar{\nu}$. This proves the strong uniqueness for the
$(A,C)$-SDE.

\smallskip
Finally let $\bar{\nu}$ be a Feller convolution semigroup weakly
dominated by $\nu^c$ and dominating $\nu^{s}$. The fact that
$\bar{\nu}\preceq^w \nu^c$ implies that a stochastic flow $\bar{K}$ of
law $\P_{\bar{\nu}}$ can be obtained by filtering on an extension
$(N,\p)$ of $N_{\nu^c}$ the coalescing flow $\p$ with respect to a
sub-noise $\bar{N}$ of $N$. Then section \ref{constrW} shows that
$\bar{W}^n_{s,t}$ (defined from $\bar{K}$) converges (in $L^2$)
towards $\bar{W}_{s,t} = \E[W_{s,t}|\bar{\cal F}_{s,t}]$. Now,
since $\bar{\nu}\succeq \nu^{s}$, there exists (see lemma
\ref{propfiltr}) a sub-noise $\bar{\bar{N}}$ of $\bar{N}$ such that
the flow obtained by filtering $\bar{K}$ or equivalently, the
coalescing flow, with respect to $\bar{\bar{N}}$ has law
$\P_{\nu^{s}}$. The associated white noise $\bar{\bar{W}}$ verifies
for all $s\leq t$, $x\in M$ and $f\in C^2_K(M)$ 
\beq \bar{\bar{W}}_{s,t}f(x) = \E[\bar{W}_{s,t}f(x) | \bar{\bar{\cal
F}}_{s,t}] = \E[W_{s,t}f(x) | \bar{\bar{\cal F}}_{s,t}].\eeq
Since $\bar{\bar{W}}$ has covariance $C$, it has to coincide with $W$
and $\bar{W}=W$. 

Thus, $(\bar{K},\bar{W})$ solves the $(A,C)$-SDE so that $\bar{\nu}$
is a diffusion convolution semigroup whose local characteristics are
$(A,C)$. \qed

\subsection{The coalescing case~: martingale representation.} \label{MR}
On the probability space $(\Omega^0,\mathcal{A}^0,\P_{\nu^c})$, let
${\cal F}^{\nu^c}$ be the filtration $({\cal F}^{\nu^c}_{0,t})_{t\geq
0}$ and ${\cal M}({\cal F}^{\nu^c})$ be the space of locally square
integrable ${\cal F}^{\nu^c}$-martingales.
\bprop\label{repprev} For all ${\cal F}^{\nu^c}$-martingale
$M=(M_t)_{t\in\RR^+}$, there exist predictable processes
$\Phi^\alpha=(\Phi^\alpha_s)_{s\geq 0}$ such that 
\beq\label{PRP} M_t=\sum_\alpha\int_0^t\Phi^\alpha_s~W^\alpha(ds). \eeq
\eprop
\brem {\em Of course, this does not imply that $\cF^{\nu^c}$ is
  generated by $W$.}\erem
\prf We follow an argument by Dellacherie (see Rogers-Williams
(V-25) \cite{R-W}). Suppose there exists $F\in L^2({\cal F}^{\nu^c}_{0,\infty})$
orthogonal in $L^2({\cal F}^{\nu^c}_{0,\infty})$ to all stochastic integrals of
$(W^\alpha)_\alpha$ of the form (\ref{PRP}), then 
$M_t=\E[F|{\cal F}^{\nu^c}_{0,t}]$ is orthogonal to $W^\alpha$
for all $\alpha$, i.e. $\langle M,W_{0,\cdot}^\alpha\rangle_t=0$.

Let $\tau=\inf\{t,~|M_t|=1/2\}$ and $\hat{\P}_{\nu^c}=(1+M_\tau)\cdot
\P_{\nu^c}$. Since $M$ is a uniformly integrable martingale and $\tau$ a
stopping time (with $1+M_\tau\geq 1/2$), $\hat{\P}_{\nu^c}$ is a probability
measure on  $(\Omega^0,\cA^0)$. 
Since $\langle M,W_{0,\cdot}^\alpha\rangle_t=0$, we get that under $\hat{\P}_{\nu^c}$,
$(W^\alpha_{0,t})_\alpha$ is a family of independent Brownian motions.

\medskip
We are now going to prove that since {\bf (U)} is satisfied, we must
have $\P_{\nu^c}=\hat{\P}_{\nu^c}$, which implies $M_t=0$ and a contradiction.

\medskip
Let $F=\prod_{i=1}^nf_i(\p_{0,t_i}(x_i))$, for $f_1,\dots,f_n$ in 
$C_K^2(M)$, $t_1,\dots,t_n$ in $\RR^+$ and $x_1,\dots,x_n$ in $M$. We
know that under $\P_{\nu^c}$, for all $1\leq i\leq n$, $(\p_{0,t}(x_i),~t\geq
0$) is a solution of the SDE
\beq\label{SDEi} dg_i(\p_{0,t}(x_i)) = \sum_\alpha V_\alpha
g_i(\p_{0,t}(x_i))W^\alpha(dt) + Af(\p_{0,t}(x_i))dt,\eeq
for all $g_1,\dots,g_n$ in $C_K^2(M)$. Note that under $\hat{\P}_{\nu^c}$, these
SDEs are also satisfied. Since under $\hat{\P}_{\nu^c}$, $(W^\alpha)_\alpha$
is a family of independent Brownian motions, 
$((\p_{0,t}(x_i),~t\geq 0),~1\leq i\leq n)$ is a 
coalescing solution of the martingale problem associated with $A^{(n)}$
and {\bf (U)} implies that the law of $((\p_{0,t}(x_i),~t\geq
0),~1\leq i\leq n)$ is the same under $\P_{\nu^c}$ and under
$\hat{\P}_{\nu^c}$. Therefore $\hat{\E}[F]=\E[F]$, where $\hat{\E}$ denotes
the expectation with respect to $\hat{\P}_{\nu^c}$. 

\smallskip
To conclude that $\hat{\P}_{\nu^c}=\P_{\nu^c}$, we need to prove $\hat{\E}[F]=\E[F]$
with $F=\prod_{i=1}^nf_i(\p_{s_i,t_i}(x_i))$ for all $f_1,\dots,f_n$
in $C_K^2(M)$, $0\leq s_i<t_i$ in $\RR^+$ and $x_1,\dots,x_n$ in
$M$. This can be proved the same way but using the kernel
$\tilde{K}_t$ introduced in section \ref{cadlag}. In this case
$\tilde{K}_t=\delta_{\tilde{\p}_t}$, where $\tilde{\p}_t~:~\RR^+\times
M\to \RR^+\times M$ is measurable. Then
$F=\prod_{i=1}^n\tilde{f_i}(\tilde{\p}_{t_i}(s_i,x_i))$ and
$(\tilde{\p}_t(s_i,x_i),~t\geq 0)$ is a solution of an SDE on 
$\RR^+\times M$. \qed

\subsection{The coalescing case~: the linear noise.}\label{linear}
Let us remark that if $\nu$ is a diffusion convolution semigroup, then
$N_\nu$ is a predictable noise (see proposition \ref{martcont}),
i.e. $\mathcal{M}(\mathcal{F}^\nu)$ is formed of continuous martingales
(in particular, a Gaussian noise is predictable). Following Tsirelson
\cite{Tsirelson}, a linear representation of a predictable noise
$N=(\Omega,\cA, ({\cal F}_{s,t})_{s\leq t},\P,(T_h)_{h\in\RR})$
is a family of real random variables $X=(X_{s,t},~s\leq t)$ such that  
\bdes
\ita $X_{s,t}\circ T_h=X_{s+h,t+h}$ for all $s\leq t$ and all $h\in\RR$,
\itb $X_{s,t}$ is ${\cal F}_{s,t}$-measurable for all $s\leq t$,
\itc $X_{r,s}+X_{s,t}=X_{r,t}$ a.s., for all $r\leq s\leq t$.
\edes
The space of linear representations is a vector space. Equipped with
the norm $\|X\|=(\E[|X_{0,1}|^2])^{\frac{1}{2}}$, it is a
Hilbert space we denote by $H_{\hbox{{\small lin}}}$. Let $H^0_{\hbox{{\small lin}}}$ be
the orthogonal in $H_{\hbox{{\small lin}}}$ of the one-dimensional vector
space constituted of the representation $X_{s,t}=v(t-s)$ for
$v\in\RR$, then $H^0_{\hbox{{\small lin}}}$ is constituted with the centered
linear representations. Note that if $X\in
H^0_{\hbox{{\small lin}}}$ with $\|X\|=1$, then
$(X_{0,t})_{t\geq 0}$ is a standard Brownian motion. The Hilbert space
$H^0_{\hbox{{\small lin}}}$ is a Gaussian system and every
$X\in H^0_{\hbox{{\small lin}}}$ is a real white
noise.

Note that if $X$ and $Y$ are orthogonal linear representations then
$X$ and $Y$ are independent.

\medskip
For all $-\infty\leq s\leq t\leq\infty$, let 
${\cal F}_{s,t}^{\hbox{{\small lin}}}$ be the $\sigma$-field generated
by the random variables $X_{u,v}$ for all $X\in H^0_{\hbox{{\small
lin}}}$ and $s\leq u\leq v\leq t$, and completed by all
$\P$-negligible sets of $\mathcal{F}_{-\infty,+\infty}$. Then
$N_{\hbox{\small{lin}}} :=(\Omega,\cA, ({\cal F}^{\hbox{{\small lin}}}_{s,t})_{s\leq t},
\P,(T_h)_{h\in\RR})$ is a noise. It is called the linearizable part of
the noise $N$. The noise $N_{\hbox{{\small lin}}}$ is a maximal
Gaussian sub-noise of $N$, hence $N$ is Gaussian if and only if
$N_{\hbox{{\small lin}}}=N$. When $N_{\hbox{{\small lin}}}$ is trivial
(i.e. constituted of trivial $\sigma$-fields), one says that $N$ is a
black noise (when $N$ is not trivial).  

\bthm\label{Wlin} $N^W_{\nu^c}=N^{\hbox{\em{\small lin}}}_{\nu^c}$.\ethm
\prf Let $H^W$ be the space of centered linear representations of the
noise $N^W_{\nu^c}$. Then $H^W$ is an Hilbert space (an orthonormal
basis of $H^W$ is given by $\{(W^\alpha_{s,t})_{s\leq t}\}$) and we
have $H^W\subset H^0_{\hbox{{\small lin}}}$. This implies that
$N^W_{\nu^c}$ is a Gaussian sub-noise of $N^{\hbox{{\small lin}}}_{\nu^c}$.

\medskip
If $N^W_{\nu^c}\neq N^{\hbox{{\small lin}}}_{\nu^c}$ then there exists a linear
representation $X\neq 0\in H^0_{\hbox{{\small lin}}}$ orthogonal to
$H^W$ and therefore independent of $\{W^\alpha\}$. Since
$(X_{0,t})_{t\geq 0}\in {\cal M}({\cal F})$, proposition \ref{repprev}
implies that the martingale bracket of $X_{0,t}$ equals 0. This is a
contradiction. \qed

\medskip In section \ref{isotrope}, we give an example of a stochastic
coalescing flow whose noise is predictable but not Gaussian. It is an
example of non-uniqueness of the diffusion convolution semigroup
associated with a set of local characteristics.

\brem  \label{sgnrem} {\em In example \ref{sgn}, although the
covariance function $C$ is not continuous, it is still possible to
construct a white noise $W$
from the coalescing flow $(\p_{s,t},~s\leq t)$. For all $s<t$, we set
$W_{s,t} = \int_s^t \hbox{sgn}(\p_{s,u}(0))~d\p_{s,u}(0)$. Then we
have $W_{s,t} = \int_s^t \hbox{sgn}(\p_{s,u}(x))~d\p_{s,u}(x)$ for all
$x\in\RR$. Therefore one can check that $W=(W_{s,t},~s\leq t)$ is a
real white noise.

The coalescing flow $(\p_{s,t},~s\leq t)$ solves the SDE
\beq \p_{s,t}(x) = \int_s^t \hbox{sgn}(\p_{s,u}(x))~dW_u,\hbox{ for }
s<t \hbox{ and } x\in\RR. \label{TanakaSDE} \eeq
The results of this subsection apply since proposition \ref{repprev}
is also satisfied if we only assume the uniqueness in law of the
coalescing solutions \footnote{i.e. such that if $(X^1,\dots,X^n)$
solves the SDE then if for $i\neq j$ and $s\geq 0$, $X^i_s=X^j_s$ then
$X^i_t=X^j_t$ for all $t\geq s$.} of the SDE satisfied by the
$n$-point motion (i.e. the SDE (\ref{SDEi})), which here is almost
obvious. Therefore, the linear part of the noise generated by the
coalescing flow is given by the noise generated by $W$. But since the
strong solution of the SDE (\ref{TanakaSDE}) is not a flow of
mappings, the coalescing flow is not a strong solution. Therefore, we
recover the result of Warren \cite{warren} and Watanabe
\cite{Watanabe} that the noise of this stochastic coalescing flow is
predictable but not Gaussian.

The strong solution given in section \ref{sgn} can be recovered by
filtering the coalescing solution with respect to the noise generated
by $W$.}\erem

\setcounter{equation}{0}
\section{Isotropic Brownian flows.}\label{isotrope}
In this section, we give examples of compatible families of Feller
semigroups. They are constructed on $M$, a two-point symmetric
space, with $C$ an isotropic covariance function on the space of
vector fields and the semigroup of a Brownian motion on $M$. 

\subsection{Isotropic covariance functions.}
Let $M=G/K$ be a two-point symmetric space. This class of spaces
includes euclidean spaces, hyperbolic spaces and spheres, see
\cite{helgason}, chapter III. $G$ is the group of isometries on $M$.
A covariance function $C$ is said isotropic if
\beq\label{isocond}C(g\cdot\xi,g\cdot\xi')=C(\xi,\xi')\eeq
for all $g\in G$ and $(\xi,\xi')\in (T^*M)^2$ and where $g\cdot \xi = Tg(\xi)$
(or $g\cdot(x,u)=(gx,Tg_xu)$ for $(x,u)\in T^*M$).

Examples of isotropic covariances are given by Monin and Yaglom  in
\cite{yaglom} on $\RR^d$ and by Raimond \cite{ra,ra'} on the sphere
and on the hyperbolic plane. In these examples, the group $G$ of
isometries on $\RR^d$ (making $\RR^d$ homogeneous) is generated by
$O(d)$ and by the translations. For the sphere $\SS^d$, this group is
$O(d+1)$ and for the hyperbolic space, it is $O(d,1)$.

\subsection{A compatible family of Markovian semigroups.}
Let $C$ be an isotropic covariance on ${\cal X}(M)$, the space of vector
fields on the two-point symmetric space $M=G/K$. To this isotropic covariance
function is associated a Brownian vector field on $M$ (i.e. a ${\cal
X}(M)$-valued Brownian motion $W$ such that 
$\E[\langle W_t,\xi\rangle\langle W_s,\xi'\rangle] = t\wedge s~
C(\xi,\xi')$).
Let $\P$ be the associated Wiener measure, constructed
on the canonical space $\Omega=\{\omega:\RR^+\to {\cal X}(M)\}$, equipped
with the $\sigma$-field $\cA$ generated by the coordinate
functions. 

We denote by $W$ the random variable $W(\omega)=\omega$. $W$ is a
Brownian vector field of covariance $C$ which is isotropic in the sense
that for all $g\in G$, $(Tg_x^{-1}W_t(gx),~t\in\RR^+,~x\in M)$ is a
Brownian vector field of covariance $C$.

\medskip Let $\P_t$ be the heat semigroup on $M$, $m$ the volume
element and $\Delta$ the Laplacian.

Let $(S_t,~t\geq 0)$ be the family of random operators defined in
\cite{ljr}, associated with $W$ and to the heat semigroup
$\P_t$. Following \cite{ljr}, we define the associated semigroups of
the $n$-point motion, $\P^{(n)}_t=\E[S_t^{\otimes n}]$ (with
$\P^{(1)}_t=\P_t$). Then, it is obvious that $(\P^{(n)}_t,~n\geq 1)$
is a compatible family of Markovian semigroups of operators acting on
$L^2(m^{\otimes n})$. We now prove that these semigroups are induced by Feller
semigroups (the question was raised in \cite{Ma}).

\medskip One can extend $(W_t)_{t\geq 0}$ into a vector field valued
white noise $(W_{s,t},~s\leq t)$ of covariance $C$ such that
$W_t=W_{0,t}$ for $t\geq 0$ and associate to it a stationary cocycle
of random operators $(S_{s,t},~s\leq t)$ such that $S_{0,t}=S_t$ for
$t\geq 0$. 

\subsection{Verification of the Feller property.}\label{sectverif}
For all $g\in G$, let $L_g:\Omega\to\Omega$ defined by 
$L_g\omega_t (\cdot)=Tg^{-1}(\omega_t(g\cdot))$, for all $t\in \RR$ and $x\in
M$. Then $L_g$ is linear and for all $g_1$ and $g_2$ in $G$, 
$L_{g_1g_2}=L_{g_1}L_{g_2}$ (i.e. $g\mapsto L_g$ is a representation
of $G$). It is easy to check that for all $g\in G$,
$(L_g)^*\P=\P$. Note that this last condition is also a
characterization that $C$ is isotropic.

\medskip
For all $g\in G$, $L_g$ induces a linear transformation on
$L^2(\Omega,\cA,\P)$ we will also denote by $L_g$. Then for all
$f\in L^2(\Omega,\cA,\P)$, we have
$L_gf(\omega)=f(L_g\omega)$. This transformation is unitary since 
$$\|L_gf\|^2=\int f^2(L_g\omega)~\P(d\omega)=\int
f^2(\omega)~((L_g)^*\P)(d\omega)=\|f\|^2,$$
(where $\|\cdot\|$ denotes the $L^2(\P)$-norm).

\bprop For all $v\in L^2(\Omega,\cA,\P)$, the mapping $g\mapsto
L_g v$ is continuous. \eprop
\prf Note that, since $L$ is a representation, it is enough to prove
the continuity at $e$, the identity element in $G$.
\brem \label{gcont} {\em Let $(v_n,~n\in\NN)$ be a sequence in
$L^2(\Omega,\cA,\P)$ converging towards $v\in L^2(\Omega,{\cal
A},\P)$ as $n\to \infty$ such that $\lim_{g\to e}L_gv_n=v_n$ for all
integer $n$, then $\lim_{g\to e} L_gv=v$. Indeed, since for all $g\in
G$, $L_g$ is unitary, $\|L_gv-v\|\leq 2\|v_n-v\|+\|L_gv_n-v_n\|$.
Hence $\limsup_{g\to e}\|L_gv-v\|\leq 2\|v_n-v\|$ for all
integer $n$.} \erem

We first prove that $\lim_{g\to e}L_g v=v$ for every $v$
of the form $\sum_iW_{t_i}(\xi_i)$ (with $W_t(x,u)=\langle
W_t(x),u\rangle$, where $\langle \cdot,\cdot\rangle$ denotes the
Riemannian metric)~: 
$$\|L_g(\sum_iW_{t_i}(\xi_i))-\sum_iW_{t_i}(\xi_i)\|^2 =
2\sum_{i,j}t_i\wedge t_j(C(\xi_i,\xi_j)-C(g\cdot\xi_i,\xi_j))$$
which converges towards 0 as $g$ tends to $e$. 

Let $H$ denote the closure (in $L^2(\Omega,\cA,\P)$) of the class
of all $v$ of the form $\sum_iW_{t_i}(\xi_i)$. Remark \ref{gcont}
implies that $\lim_{g\to e}L_gv=v$ holds for all $v\in H$. 

\smallskip
It is well known that $L^2(\Omega,\cA,\P)$ is the orthogonal sum
of the Wick powers $H^n$ of $H$ (See \cite{Simon}), also called the
$n$-th Wiener chaos (see \cite{Neveu}), $H^0$ is constituted by the
constants. The space $H^n$ is isometric to the symmetric tensor
product Hilbert space $H^{\otimes^s n}$. We now prove that $\lim_{g\to
e}L_gv=v$ holds for all $v\in H^n$. For all 
$v=v_1\otimes^s\cdots\otimes^sv_n\in H^n$ (or $:v_1v_2\cdots v_n:$ in
wick notation), with $v_1,\dots,v_n$ in $H$, 
\beqarr \|L_gv - v\|
&\leq& \sum_j\|L_gv_1\otimes^s\cdots\otimes^sL_gv_{j-1}\otimes^s
(L_gv_j-v_j)\otimes^sv_{j+1}\otimes^s\cdots\otimes^sv_n\|\\
&\leq&\sqrt{n!}\sum_j\|L_gv_j-v_j\|\times\prod_{i\neq j}\|v_i\|
\eeqarr
which converges towards 0 as $g$ tends to $e$. Since the class of
linear combinaisons of elements of the form $v_1\otimes^s\cdots\otimes^sv_n$
is dense in $H^n$, we have $\lim_{g\to e}L_gv=v$ for all $v$ in
$H^n$. And we conclude since
$L^2(\Omega,\cA,\P)=\oplus_{n\geq 0}H^n$. \qed

\medskip
For all $x\in M$, $s\leq t$ and $f\in C_0(M)$,
since $\P^{(1)}_\eps$ is absolutely continuous with respect to $m$, we
have
\beq \P^{(1)}_\eps S_{s+\eps,t}f(x) = 
\E[\P^{(1)}_{\eps'} S_{s+\eps',t}f(x) |{\cal F}_{s+\eps',t}], \eeq 
for $0<\eps'\leq \eps$. Thus, for all $s<t$, $\P^{(1)}_\eps
S_{s+\eps,t}f(x)$ is a martingale as $\eps$ decreases.
This martingale converges and we denote its limit by $K_{s,t}f(x)$.
Then $S_{s,t}f=K_{s,t}f$ in $L^2(m\otimes \P)$ and $\P^{(n)}_t =
\tilde{P}_t^{(n)}$ $m^{\otimes n}$-a.e., where $\tilde{P}_t^{(n)}$ denotes
$\E[K_{s,t}^{\otimes n}]$.

\blem\label{Scont} The mapping $x\mapsto K_{s,t}f(x)$ is continuous for all
Lipschitz function $f$ and all $s\leq t$.\elem
\prf Note that for all $g\in G$ and all $x\in M$,
\beq L_gK_{s,t}f(x)=K_{s,t}f^{g^{-1}}(gx)\eeq
where $f^{g^{-1}}(x)=f(g^{-1}x)$. We then have
\beqarr \|K_{s,t}f(gx)-K_{s,t}f(x)\| &\leq&
\|K_{s,t}f(gx)-K_{s,t}f^{g^{-1}}(gx)\|\\
&& +\quad\|L_gK_{s,t}f(x)-K_{s,t}f(x)\|. \eeqarr

Hence $\lim_{g\to e}K_{s,t}f(gx)=K_{s,t}f(x)$ since
$\lim_{g\to e}L_gK_{s,t}f(x)=K_{s,t}f(x)$ and
$\|K_{s,t}f(gx)-K_{s,t}f^{g^{-1}}(gx)\|\leq
\|f-f^{g^{-1}}\|_\infty$ which converges towards 0 (since
$|f(x)-f^{g^{-1}}(x)|\leq Cd(x,g^{-1}x)$, which converges towards $0$
as $g\to e$). This implies the lemma. \qed

\bprop \bdes \ita $(\tilde{\P}^{(n)}_t,~n\geq 1)$ is a compatible
family of Feller semigroups.
\itb The associated convolution semigroup $\nu^s=(\nu^s_t)_{t\geq 0}$ is
a diffusion convolution semigroup with local characteristics
$(\frac{1}{2}\Delta,C)$. \edes  \eprop
\prf For all bounded Lipschitz functions $f_1,\dots,f_n$,
lemma \ref{Scont} implies that
$(x_1,\dots,x_n)\mapsto\tilde{\P}^{(n)}_tf_1\otimes\cdots\otimes 
f_n(x_1,\dots,x_n)=\E[\prod_{i=1}^nK_{s,t}f_i(x_i)]$ is
continuous. This suffices to prove {\bf (a)} (the proof that
$\lim_{t\to 0}\P^{(n)}_th(x)=h(x)$ for all $h\in C(M^n)$ is the same
as in lemma \ref{Fel2n}).

\medskip
To prove {\bf (b)}, notice that It\^o's formula for $(S_{s,t},~s\leq
t)$ (see theorem 3.2 in \cite{ljr}) implies that for all $f\in
C^2_K(M)$ and $s\leq t$,
\beq K_{s,t}f(x) = f(x) + \int_s^t K_{s,u}(Wf(du))(x) +
\frac{1}{2}\int_s^t K_{s,u}(\Delta f)(x)~du, \eeq
i.e. $(K,W)$ solves the $(\frac{1}{2}\Delta,C)$-SDE. \qed

\subsection{Classification.}

Let $\nu^s$ be the diffusion convolution semigroup constructed
above. It defines a strong solution of the
$(\frac{1}{2}\Delta,C)$-SDE. Note that there is no pure diffusion.

Let $(d_t)_{t\geq 0}$ denote the distance process induced by the
2-point motion $X^{(2)}_t=(X_t,Y_t)$ (then $d_t=d(X_t,Y_t)$). The
isotropy condition and the fact that in two point homogeneous spaces,
pairs of equidistant points, can be exchanged by an isometry imply that
$d_t$ is a real diffusion. We denote in the following the law of this
diffusion starting from $r\geq 0$ by $\P_r$. Let
$H_r=\inf\{t>0,~d_t=r\}$.

\bprop \label{classprop} \bdes \item[(1)] $\nu^s$ defines a non-coalescing flow of maps
(i.e. such that the 2-point motion starting outside of the diagonal
never hits the diagonal)  if and only if 0 is a natural boundary
point, i.e. if 
\beq \forall r>0,~\P_r[H_0<\infty]=0 \hbox{ and }
\P_0[H_r<\infty]=0. \eeq 
\item[(2)] $\nu^s$ defines a coalescing flow of maps if and only if 0 is a closed
exit boundary point, i.e. if
\beq \exists r>0,~\P_r[H_0<\infty]>0 \hbox{ and } \forall r>0,
~\P_0[H_r<\infty]=0. \eeq
\item[(3)] $\nu^s$ defines a turbulent flow \footnote{We recall that a
turbulent flow was defined as a stochastic flow of kernels which is
not a flow of maps and without pure diffusion.} without hitting
(i.e. such that the 2-point motion starting outside of the diagonal
never hits the diagonal) if and only if 0 is an open entrance boundary
point, i.e. if
\beq \forall r>0,~\P_r[H_0<\infty]=0 \hbox{ and } \exists r>0,
~\P_0[H_r<\infty]>0. \eeq
\item[(4)] $\nu^s$ defines a turbulent flow with hitting (i.e. such that the
2-point motion starting outside of the diagonal hits the diagonal with
a positive probability) if and only if 0 is a reflecting regular boundary point, i.e. if
\beq \exists r>0,~\P_r[H_0<\infty]>0 \hbox{ and } \exists r>0,
~\P_0[H_r<\infty]>0. \eeq \edes

In all cases except {\bf (4)}, $\nu^s$ is the unique diffusion
convolution semigroup with local characteristics
$(\frac{1}{2}\Delta,C)$.

In case {\bf (4)}, called the intermediate phase, $\nu^c\neq\nu^s$ and
theorems \ref{cor} and \ref{Wlin} apply. Thus $N_{\nu^c}$ is a predictable
non-Gaussian noise. \eprop

\prf The proof of {\bf (1)}, {\bf (2)}, {\bf (3)} and {\bf (4)} is
straightforward. Notice that the local characteristics satisfy {\bf
(U)}. In all cases, $\nu^s$ defines a strong solution of the
$(\frac{1}{2}\Delta,C)$-SDE. This with theorem \ref{cor}
implies that in the coalescing case {\bf (2)}, since 
$\nu^s=\nu^c$, $\nu^s$ is the unique diffusion convolution
semigroup whose local characteristics are $(\frac{1}{2}\Delta,C)$.

In the non-coalescing case {\bf (1)} and in the turbulent case without
hitting {\bf (3)}, the fact that $\nu^s$ is the unique diffusion
convolution semigroup whose local characteristics are
$(\frac{1}{2}\Delta,C)$ follows directly from {\bf (U)}.

In the intermediate phase {\bf (4)}, we must have $\nu^c\neq \nu^s$
since $\nu^s$ defines a turbulent flow and $\nu^c$ a flow of
maps. Moreover, condition {\bf (B)} holds so that we can
conclude using theorems \ref{cor} and \ref{Wlin}. \qed

\brem  The $(\frac{1}{2}\Delta,C)$-SDE has a solution, unique
in law except in the intermediate phase, in which case all solutions
are obtained by filtering, on an extension $(N,\p)$ of the noise of the
coalescing solution, this coalescing solution $\p$ with respect to a
sub-noise of $N$ containing $W$. \erem

\brem  The conditions involving the distance process can be verified
using the speed and scale measures of this process which are
explicitly determined by the spectral measures of the isotropic fields
(cf \cite{ljr} for $\RR^d$ and for $\SS^d$). \erem

\subsection{Sobolev flows.}
In \cite{ljr}, Sobolev flows $(S_{s,t},~s\leq t)$
on $\RR^d$ and on $\SS^d$ are studied. The Sobolev covariances are
described with two parameters $\alpha>0$ and $\eta\in [0,1]$.
The associated self-reproducing spaces are Sobolev spaces of vector fields of
order $(d+\alpha)/2$. The incompressible and gradient subspaces are
orthogonal and respectively weighted by factors $\eta$ and $1-\eta$.

\medskip
Let us apply the results obtained in \cite{ljr}. We will call the
stochastic flow associated with $(S_{s,t},~s\leq t)$ (see section
\ref{Veretennikov} and \ref{sectverif}) Sobolev flow as well.
When $\alpha>2$, we are in case {\bf (1)} and Sobolev flows are flows
of diffeomorphisms. More interestingly, when $0<\alpha<2$ then 
\bdes \item[] If $d\in\{2,3\}$ and $\eta<1-\frac{d}{\alpha^2}$, we are 
in case {\bf (2)} and the Sobolev flow is a coalescing flow.
\item[] If $d\geq 4$ or if $d\in\{2,3\}$ and
$\eta>\frac{1}{2}-\frac{(d-2)}{2\alpha}$, we are in case {\bf (3)} and
the Sobolev flow is turbulent without hitting.
\item[] if $d\in\{2,3\}$ and $1-\frac{d}{\alpha^2} <
\eta<\frac{1}{2}-\frac{(d-2)}{2\alpha}$, we are in case {\bf (4)}
(i.e. the intermediate phase) and the Sobolev flow is turbulent with
hitting. \edes
In dimension 1, the parameter $\eta$ vanishes. The critical case was
studied in \cite{air,fang,mal}. There is a strong coalescing solution
for $\alpha\in [1,2[$ and an intermediate phase for $\alpha\in]0,1[$.

By construction, in all these cases, the noise generated by the
Sobolev flows are Gaussian noises. For the intermediate phase, in
which there exist two different solutions to the
$(\frac{1}{2}\Delta,C)$-SDE (namely the coalescing one and the
turbulent one), the noise of the associated coalescing flow is
predictable but not Gaussian.

\medskip
These different cases are represented by the phase diagram below, for
the homogeneous space $\SS^3$. Recall that a flow of diffeomorphisms
is called stable (respectively unstable) when the first Lyapounov
exponent is negative (respectively positive). These exponents actually
converge actually towards $-\infty$ or to $+\infty$ as $\alpha$
approches the critical value $2$.

\begin{center}
\includegraphics[scale=0.9]{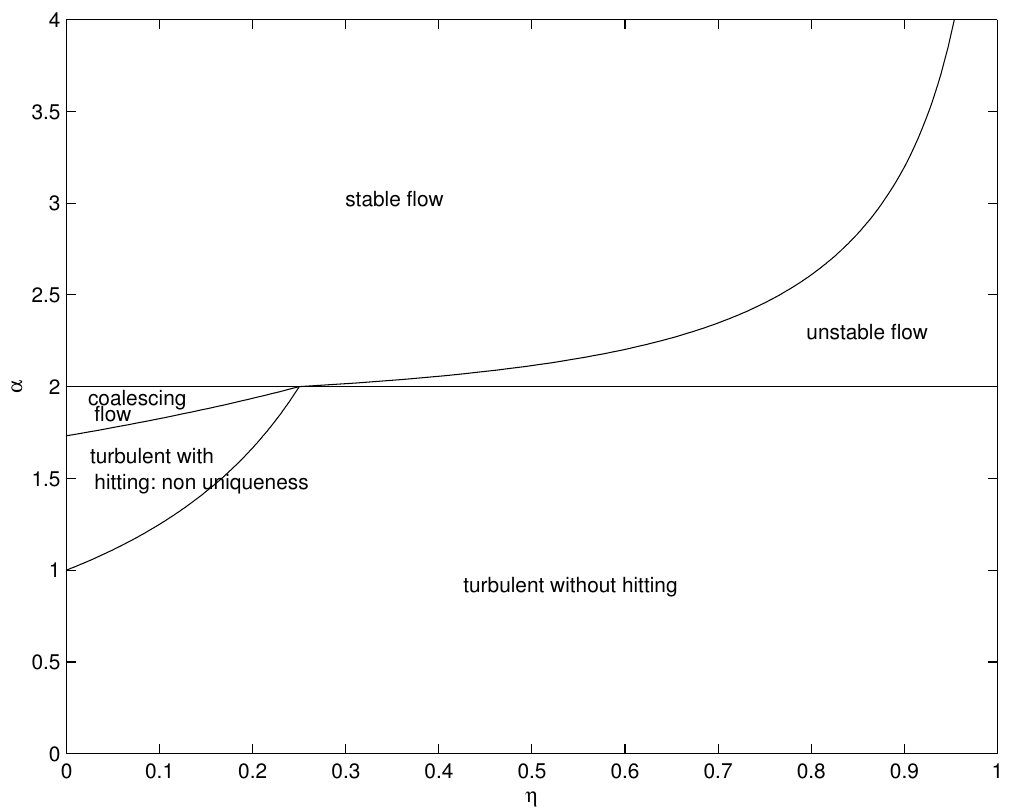}
\end{center}

%

\setcounter{equation}{0}
\section{Conclusion.}
Looking at the phase diagram above, it looks as this case has been
fully analysed. 

The three different types of motion which can be defined by a consistent system 
of Feller semigroups appear in this picture: Flows of non coalescing maps occur 
when, for the two point motion, the diagonal and the complement of the diagonal 
are absorbing.
 
When the first condition fails, i.e. when the diagonal is not
absorbing, we get a diffusive flow, i.e. a flow of non trivial Markov
kernels. We see in this exemple that this can happen without pure
diffusion, i.e. when the evolution equation has no dissipative
term. In that case we say that the flow is turbulent. It can be viewed
as an effect of extreme  unstability due to the importance of very
high frequency divergence free components in
the velocity field near the diagonal. 

When the second condition fails, i.e. when the complement of the
diagonal is not absorbing, we get flows of coalescing maps.
We see, in the intermediate phase, that a turbulent and a coalescing flow can 
have the same local characteristics. This happens when both conditions fail for 
the two point motion associated with the turbulent flow.

Moreover, it is likely that at least in the other
isotropic situations, a very similar picture will occur, the
parameters being the singularity of the covariance on the diagonal and
the balance between gradient and incompressible velocity fields.

Yet there is still some important work to do about the intermediate
phase. We know there exists two remarkable distinct solutions in that
case for the SDE: the coalescing flow, the noise of which is not
linear but for which the linear part has been identified as the
velocity white noise $W$, and the unique strong solution which is a
flow of non trivial kernels obtained by averaging the coalescing flow
with respect to $W$. Other solutions do exist and we have shown that
their associated convolution semigroups are weakly dominated by the
``coalescing'' convolution semigroup and dominate the ``strong'' or
``linear'' one. But this classification should be made analytically
precise and one can conjecture it involves a ``gluing'' parameter on
the diagonal (see section \ref{exemple}, arXiv math.PR/0203221 and
math.PR/0212269 for first steps in this direction.) Moreover, the non
linear part of the relevant noises remains to be fully
analysed. Finally, one can expect that more complex phenomena occur
for SDEs in which a multiplicity of weak solutions with different
one-point motions do exist. Hence this paper can only be a step in the
understanding of the multiplicity of flows with given velocity field,
or given local characteristics.

\end{document}